\documentclass[10pt,a4paper]{amsart}
\usepackage{macro}
\usepackage{amsxtra}
\usepackage{xcolor}
\date{\today}
\begin{document}
\title{Quantum Lefschetz theorem Revisited} \author{Jun Wang }
\maketitle

\begin{abstract}
  Let $X$ be any smooth Deligne-Mumford stack with projective coarse moduli, and
  $Y$ be a smooth complete intersection in $X$ associated with a direct sum of
  semi-positive line bundles. We will introduce a useful and broad class known as admissible
  series for discussing quantum Lefschetz theorem. For any admissible series on the Givental's Lagrangian cone of
  $X$, we will show
  that a hypergeometric modification of the series lies on the Lagrangian cone
  of $Y$. This confirms a prediction from Coates-Corti-Iritani-Tseng about the
  genus zero quantum Lefschetz theorem beyond convexity. In our quantum
  Lefschetz theorem, we use extended variables to formulate the hypergeometric
  modification, which may be of self-independent interest.
\end{abstract}

\renewcommand{\thefootnote}{}
\tableofcontents

\renewcommand{\thefootnote}{\arabic{footnote}}
\section{Introduction}

\emph{Gromov-Witten (GW) invariants} count the (virtual) number of stable maps from
curves to a target variety/stack $X$ with prescribed conditions. \emph{Quantum Lefschetz}
is one of the central topics in Gromov-Witten theory and it compares the GW invariants of a complete
intersection and its ambient space. On the other hand, in Givental's formalism\cite{Givental_2004, coates07_quant_rieman_roch_lefsc_serre,
  tseng10_orbif_quant_rieman_roch_lefsc_serre}, all information of
genus-zero Gromov-Witten invariants is encoded in an \emph{overruled Lagrangian
  cone} $\mathcal L_{X}$ sitting inside
an infinite dimensional symplectic vector space. Consequently, a natural approach to the
quantum Lefschetz problem in genus zero is to find an explicit slice on
Givental's Lagrangian cone of the complete intersection using a given slice on
Givental's Lagrangian cone of the ambient space, This approach is often referred
to as a \emph{mirror theorem} in the literature. The earliest quantum Lefschetz
theorem can be traced back to the verification of the famous mirror conjecture
for quintic threefolds\cite{Candelas_1991, Givental_1996, Lian_1999}. Since
then, more cases related to quantum Lefschetz have been proved. In genus zero,
Givental-style quantum Lefschetz theorems can be divided into two categories:
\begin{enumerate}
\item The complete intersection $Y$ is associated with a direct sum of
  \emph{convex} line bundles over the ambient space $X$, see e.g.,\cite{Givental_1998, coates07_quant_rieman_roch_lefsc_serre, coates2019some}. The reason for
  requiring the convexity condition is that we need to apply the so-called
  Quantum Lefschetz principle proved in~\cite{kim03_funct_inter_theor_conjec_cox_katz_lee}. This principle
  expresses the Gromov-Witten invariants of the complete intersection as an
  Euler class of a certain bundle over the moduli stack of genus-zero stable
  maps to $X$. See also \cite{Coates_2012} for some examples where this
  principle can fail if we drop the convexity condition, where the examples are
  associated with semi-positive line bundles in the sense of this paper. We note
  that convexity is a very rare condition when $X$ is a not a variety.
\item When convexity fails, proving a quantum Lefschetz theorem is much more
  difficult. There has been significant progress on this case recently by
  solving the genus zero quasimap wall-crossing
  conjecture~\cite{wang19_mirror_theor_gromov_witten_theor_without_convex,
    zhou2022quasimap}. However, it's important to note that we require the
  ambient space to be a GIT quotient, which is a prerequisite for applying
  quasimap theory\cite{Ciocan_Fontanine_2014, Cheong_2015}.
\end{enumerate}
The main objective of this paper is to prove a genus-zero\footnote[1]{Recently,
  there have been significant advance in the high genus case, where the
  failure of convexity provides one main obstacle to the computation of high
  genus GW
  invariants, see e.g., \cite{ zinger08_reduc_genus_gromov_witten_invar,
    guo17_mirror_theor_genus_two_gromov, guo2018structure, chang2021polynomial,
    chen2022punctured, liu22_castel_bound_higher_genus_gromov,
    chen2021logarithmic} and their references therein.}quantum Lefschetz theorem
that goes beyond the scope of the two categories mentioned previously, where the
ambient space can be a non-GIT target and the convexity condition may not hold.
Actually our scope of quantum Lefschetz theorem is quite broad, there are not many restrictions on the ambient space and we only assume
that the complete intersection is associated with semi-positive line bundles,
which is much weaker than convexity (see\cite{Coates_2012}). The proof relies on
(and generalizes) the \emph{recursive relation} of a point on the Givental's Lagrangian cone
discovered by the author in~\cite{wang19_mirror_theor_gromov_witten_theor_without_convex},
which was used to demonstrate the genus-zero quasimap wall-crossing conjecture
for abelian GIT quotients. However, we need to
carefully pick a space carried with a $\C^{*}-$action which could provide the recursive
relation we want; in our case we will use a \emph{root-stack modification of the space
of deformation to the normal cone}. This is \emph{different} from the spaces used
in loc.cit and it has the advantage of generalizing the quantum Lefschetz
theorem proved
in loc.cit to non-GIT targets.  



\subsection{Main theorem}
Let $X$ be a smooth Deligne-Mumford stack over $\mathbb C$ with projective
coarse moduli. Let $\bar{I}_{\mu}X$ be the rigidified (cyclotomic) inertia stack
of $X$, where a $\mathbb C-$point of $\bar{I}_{\mu}X $ can be written as a pair
$(x,g)$ where $x$ is a $\mathbb C$-point of $X$ and $g$ is an element in the
isotropy group $Aut(x)$. Denote by $C$ the finite set of connected components of
$\bar{I}_{\mu}X$ with $ \bar{I}_{c}X$ being the component corresponding to the
element $c\in C$. The involution on $\bar{I}_{\mu}X$ by sending $(x,g)$ to
$(x,g^{-1})$ also induces an involution on $C$, we write $c^{-1}$ to be image of
$c$ under the involution.

The \emph{overruled Lagrangian cone} $\mathcal L_{X}$ introduced by A.Givental
is comprised by the so-called (big) $J$-function:
\begin{equation}
  \begin{split}
    &J^{X}(q,\mathbf{t}(z),-z):=-z\mathbbm{1}_{X}+\mathbf{t}(z)+\sum_{\beta\in
      \mathrm{Eff}(X)}\sum_{m\geq
      0} \\
    & \frac{q^{\beta}}{m!}\phi^{\alpha}\langle
    \mathbf{t}(\bar{\psi}_{1}),\cdots,
    \mathbf{t}(\bar{\psi}_{m}),\frac{\phi_{\alpha}}{-z-\bar{\psi}_{\star}}\rangle^{X}_{0,[m]\cup\star,\beta}\
    ,
  \end{split}
\end{equation}
where the input\footnote{We will also use the notation $\mu(z)$ in this paper.}
$\mathbf{t}(z)\in H^{*}_{CR}(X,\C)[z]$ is a polynomial in $z$ with coefficients
in the Chen-Ruan cohomology $H^{*}_{CR}(X,\C)=H^{*}(\bar{I}_{\mu}X,\mathbb C)$,
the notation $q^{\beta}$ stands for the Novikov variable corresponding to the
degree $\beta$ in the cone $\mathrm{Eff}(X)$ of effective curve classes of $X$.
We note that such a choice of the input $\mathbf{t}(z)$ usually lead to some
(rather mild) convergence issue about $J-$function. One way to handle this is
using formal geometry; let $N$ be any positive integer, we will choose our input
$\mathbf{t}$ in the \emph{$\mathbb Z_{2}-$graded space}
$$ (q,t_{1},\cdots,t_{N})H^{*}(\bar{I}_{\mu}X,\mathbb
C)[z][\![t_{1},\cdots,t_{N}]\!][\![\mathrm{Eff}(X)]\!]\ ,$$ where each variable
$t_{i}$ is associated with a grading in $\mathbb Z_{2}=\mathbb Z/2 \mathbb Z$.
More precisely, $\mathbf{t}(z)$ can be written in the form of
\begin{equation}\label{eq:exseries}
  \sum_{\substack{\beta \in \mathrm{Eff}(X)\\ \vec{k}=(k_{1},\cdots, k_{N})\in
      \mathbb Z_{\geq 0}^{N}}}q^{\beta}t_{1}^{k_{1}}\cdots, t_{N}^{k_{N}}f_{\beta,\vec{k}}(z)\ ,
\end{equation} 
where $f_{\beta,\vec{k}}(z)\in H^{*}(\bar{I}_{\mu}X,\mathbb C )[z]$, $f_{0,
  \overrightarrow{0}}(z)=0$ and $t_{1}^{k_{1}}\cdots
t_{N}^{k_{N}}f_{\beta,\vec{k}}(z)$ is of even grading. See \S \ref{sec:back} for
more details.

Let $Y\subset X$ be a smooth complete intersection associated with a direct sum
of semi-positive line bundles $\oplus_{j=1}^{r}L_{j}$, i.e., the pairing
$\beta(L_{j}):=(c_{1}(L_{j}),\beta)\geq 0$ for any degree $\beta\in
\mathrm{Eff}(X)$. Then the age function $age_{g}(L_{j}|_{x})$ for each $(x,g)\in
\bar{I}_{c}X $ is constant on each connected component $\bar{I}_{c}X $ of
$\bar{I}_{\mu}X$ and we will denote the constant value to be $age_{c}(L_{j})$.
To state the main theorem in this paper, we will introduce a useful category for
discussing the quantum Lefschetz known as \emph{admissible series}. This concept
summarizes a common feature of all previous proved $J-$functions about quantum
Lefschetz. Here we will describe the content about admissible series needed in
the statement of the main theorem, see \S \ref{sec:back} for a more general
definition of admissible series.

Endow each variable $t_{i}$ and line bundle $L_{j}$ with a weight $w_{ij}\in
\mathbb Q_{\geq 0}$. We will call a tuple
$(\beta,\vec{k}=(k_{1},\cdots,k_{N}),c)\in \mathrm{Eff}(X)\times \mathbb Z_{\geq
  0}^{N}\times C$ an \emph{admissible pair} if
$$age_{c}(L_{j})\equiv \beta(L_{j})+\sum_{i=1}^{N} w_{ij}k_{i}\; \mathrm{mod}\;
\mathbb Z $$ for all line bundles $L_{j}$. Denote by $\mathrm{Adm}$ the set of
all admissible pairs. Write each component $f_{\beta,\vec{k}}(\mathbf{t},z)$ in
\ref{eq:exseries} as a sum
$$ f_{\beta,\vec{k}}(z)= \sum _{c\in C}f_{\beta,\vec{k},c}(z)$$
where $f_{\beta,\vec{k},c}(z)$ belongs to the space
$H^{*}(\bar{I}_{c^{-1}}X,\mathbb C)[z]$. Then we call the input $\mathbf{t}(z)$
an \emph{admissible (power) series} if $\mathbf{t}(z)$ can be written as in
\ref{eq:exseries} and $f_{\beta,\vec{k},c}(z)=0$ whenever $(\beta,\vec{k},c)$ is
\emph{not} an admissible pair.

Let $J^{X}(q,\mathbf{t},-z)$ be a point on the Lagrangian cone $\mathcal L_{X}$
of $X$ with $\mathbf{t}(z)$ being an admissible (power) series. Following the
tradition in the literature, to prove a quantum Lefschetz theorem, it's more
convenient to flip the sign of $z$ in the function $J(q,\mathbf{t},-z)$; in our
case, $J^{X}(q,\mathbf{t},z)$ will always have an asymptotic expansion in
variables $q^{\beta}t_{1}^{k_{1}}\cdots t_{N}^{k_{N}}$ as:
$$z+ \sum_{(\beta,\vec{k},c)\in \mathrm{Adm} }q^{\beta}t_{1}^{k_{1}}\cdots t_{N}^{k_{N}}J^{X}_{\beta,\vec{k},c}(\mathbf{t},z)\ ,$$ 
where $J^{X}_{\beta,\vec{k},c}(\mathbf{t},z)\in H^{*}(\bar{I}_{c^{-1}}X,\mathbb
C )[z,z^{-1}]]$. Define $J^{X,tw}(q,\mathbf{t},z)$ to be the hypergeometric
modification of $J^{X}(q,\mathbf{t},z)$:
\begin{equation}\label{eq:tw}
  \begin{split}
    &J^{X,tw}(q,\mathbf{t},z):=z+ \sum_{(\beta, \vec{k},c)\in \mathrm{Adm} }q^{\beta}t_{1}^{k_{1}}\cdots t_{N}^{k_{N}}J^{X}_{\beta,\vec{k},c}(\mathbf{t},z)\\
    &\cdot \prod_{j=1}^{r}\prod_{0\leq m<
      \beta(L_{j})+\sum_{i}w_{ij}k_{i}}\big(c_{1}(L_{j})
    +(\beta(L_{j})+\sum_{i}w_{ij}k_{i}-m)z \big)\ .
  \end{split}
\end{equation}
Our main theorem in this paper will be the following:
\begin{theorem}[=Theorem \ref{thm:main1}]
  Let $i:\bar{I}_{\mu}Y \rightarrow \bar{I}_{\mu}X$ be the inclusion of
  rigidified inertia stacks, then the series $i^{*} J^{X,tw}(q,\mathbf{t},-z)$
  lies on the Lagrangian cone $\mathcal L_{Y}$ of $Y$.
\end{theorem}
Here we are suppressing a change about the $J-$function defining the Lagrangian
cone $\mathcal L_{Y}$ of $Y$ by using Novikov variables from $\mathrm{Eff}(X)$
rather than $\mathrm{Eff}(Y)$, see \S \ref{subsec:spe-deg} for more details.

Even when the convexity holds, our theorem still have some new results compared
to previous quantum Lefschetz theorems. Our selection of the variables $t_{i}$
with their respective weights $w_{ij}$ plays the same role as the Novikov
variables $q^{\beta}$ with the corresponding numbers $\beta(L_{j})$, which
resembles much similarity with the usage of extended degrees for $S$-extended
$I-$functions of toric stacks in~\cite{Coates_2015, coates2019some} and we will
call the pair $(\beta,\vec{k})$ an \emph{extended degree} and call $t_{i}$
extended variables. However the appearance of extended variables beyond toric
cases is new. We note that introducing extended variables in loc.cits is a vital
aspect to establish a Laudau-Ginzburg mirror for the big (equivariant) quantum
ring of toric stacks (see~\cite{coates2020hodge}). We expect our theorem will
have the similar Mirror result (at least for ambient cohomology) and we will
return to this question in the future. It's worth noting that the method of introducing extended variables in loc.cits
cannot generalized to general cases as it relies on the $S-$extended fan
structure for toric stacks in loc.cits (see also~\cite{Jiang_2008}). This
reflects a novel aspect of our method.

\subsection{Relation to twisted Gromov-Witten invariants}
We will discuss a relationship of our quantum Lefschetz theorem to twisted
Gromov-Witten invariants, in particular a conjecture of
Coates-Corti-Iritani-Tseng.

Consider the $\C^{*}-$action on the vector bundle $E:=\oplus_{j=1}^{r}L_{j}$ via
scaling the fibers. The $\C^{*}-$equivariant Euler class of $E$ can be written
in term of the Chern roots $c_{1}(L_{j})$ as
$$e^{\C^{*}}(E):=\prod_{j}(\kappa+c_{1}(L_{j}))\ ,$$
where $\kappa$ be the equivariant parameter corresponding to the standard
representation of $\C^{*}$. Then one can define a \emph{twisted Lagrangian cone}
$\mathcal L^{tw}_{X}$ using $(e^{\C^{*}},E)$-\emph{twisted Gromov-Witten
  invariants} (See~\cite{coates07_quant_rieman_roch_lefsc_serre}
and~\cite{tseng10_orbif_quant_rieman_roch_lefsc_serre} for more details).

The series $J^{X}(q,\mathbf{t},z)$ in our main theorems are related to points on
the twisted Lagrangian cone $\mathcal L^{tw}_{X}$ in the following way: we can
associate $J^{X}$ a series defined by
\begin{equation}\label{eq:twi}
  \begin{split}
    &I^{X,tw}:= z+ \sum_{(\beta, \vec{k},c)\in \mathrm{Adm} }q^{\beta}t_{1}^{k_{1}}\cdots t_{N}^{k_{N}}J^{X}_{\beta,\vec{k},c}(\mathbf{t},z)\\
    &\cdot \prod_{j=1}^{r}\prod_{0\leq m<
      \beta(L_{j})+\sum_{i}w_{ij}k_{i}}\kappa+\big(c_{1}(L_{j})
    +(\beta(L_{j})+\sum_{i}w_{ij}k_{i}-m)z \big)\ .
  \end{split}
\end{equation}
Then we have $$\lim_{\kappa\rightarrow 0}I^{X,tw}= J^{X,tw}\ .$$ Moreover, using
the same idea in~\cite[Theorem 22]{coates2019some} where we view our $t_{i}$
here as their $x_{i}$ in loc.cit, one can prove that $I^{X,tw}$ lies on the
twisted Lagrangian cone $\mathcal L^{tw}_{X}$. In general,
Coates-Corti-Iritani-Tseng make the following conjecture (see~\cite[Conejcture
5.2]{oneto2018quantum}):
\begin{conjecture}\label{cor:ccit}
  Let $X$ be a smooth proper Deligne-Mumford stack with projective coarse moduli
  and $Y$ is a smooth complete intersection\footnote{Note that the conjecture
    doesn't require the line bundles defining the complete intersection to be
    semi-positive.} of $X$. Let $I^{tw}$ be a point on the twisted Lagrangian
  cone $\mathcal L^{tw}_{X}$ of $X$, if the limit $\lim_{\kappa\rightarrow
    0}i^{*}I^{tw}$ exits, then $\lim_{\kappa\rightarrow 0}i^{*}I^{tw}$ is a
  point on the Lagrangian cone $\mathcal L_{Y}$ of $Y$.
\end{conjecture}

This conjecture is known to be false in general, where a counterexample is found
in~\cite{sultani22_subtl_quant_lefsh_without_convex}. However we can still take
this as a guiding principle towards genus zero quantum Lefschetz theorem in
general case and we are able to verify this conjecture under mild hypotheses as
in Theorem \ref{thm:main1}.

\subsection{Sketch of the main idea of the proof}
The proof will first reduce to the hypersurface case. Set
$\mu^{X,tw}:=[J^{X,tw}(q,\mathbf{t},z)-z]_{+}$ to be the truncation in
nonnegative $z-$powers. Note that the $J-$function $J^{Y}(q,\mathbf{t},z)$ is of
the form $z+\mathbf{t}+\mathcal O(z^{-1})$, we observe that, to prove the main
theorem, we need to show that
$$ J^{Y}(q,i^{*}\mu^{X,tw},z)=i^{*}J^{X,tw}(q,\mathbf{t}, z)\ ,$$
for which we only need to consider their \emph{negative $z-$powers}. We aim to
demonstrate that both sides of the negative $z-$powers of the above equation
satisfy the \emph{same recursive relations} in each degree corresponding to
$q^{\beta}t_{1}^{k_{1}}\cdots t_{N}^{k_{N}}$ except a special consideration in
low degrees (see Lemma \ref{lem:deg0main}). First we will consider two auxiliary
spaces carried with $\C^{*}-$actions, which are root-stack modifications of a
certain orbi-$\P^{1}$ bundle over Y (see \S\ref{subsec:space1}) and the space of
deformation to the normal cone (see \S\ref{subsec:space2}). We then apply
virtual localization to express two auxiliary cycles (see
\eqref{eq:mainintegral1} and \eqref{eq:mainintegralvar}) corresponding to the
two spaces in graph sums and extract $\lambda^{-1}$ coefficient and
$\lambda^{-2}$ coefficient respectively (where $\lambda$ is an equivariant
parameter). Finally, the polynomiality of the two auxiliary cycles implies that
the coefficients must vanish, leading to the same type of recursive relations
(see also Proposition \ref{prop:chara-general} and Theorem \ref{prop:charavar}).

\subsection{Outline}
The rest of this paper is organized as follows. In \S\ref{sec:back}, we
introduce the concept of admissible series in more details than the in the
introduction and collect
some background on Gromov-Witten theory. In \S\ref{sec:prop-Lagrangian}, we will
give a recursive relation of an admissible series on the Lagrangian cone. Based
on this recursive relation, in \S\ref{sec:proofmain}, we will prove our main theorem. In the Appendix, we
carry out a detailed computation about edge contributions needed in the
localization formula in our proof.



\subsection{Notation and convention}
We work over the field $\mathbb C$. Let $D$ be an effective divisor of a smooth
stack $X$, we will denote the $\mathcal O(D)$ or $\mathcal O_{X}(D)$ to be the
line bundle associated with $D$.

For any positive integer $i$, we will use the notation $\boldsymbol{\mu}_{i}$ to
mean the finite cyclic subgroup of $\C^{*}$ of order $i$, and use the notation
$[i]$ to mean the set $\{1,2,\cdots,i\}$. For any rational number $a$, we will
use the notation $\mathrm{e}^{\;a}$ for to mean the exponential
$\mathrm{exp}(\frac{\sqrt{-1}a}{2\pi})$.

Let $\lambda$ be the equivariant class in
$H^{2}_{\C^{*}}(\{\mathrm{Spec}(\mathbb C)\})$ corresponding to the stand
representation of $\C^{*}$. For any rational number $q$, we will write $\mathbb
C_{q\lambda}$ to be the trivial line bundle with a $\C^{*}-$action of weight
$q$.

\section{Background on Gromov-Witten theory}\label{sec:back}
\subsection{Admissible series}\label{subsec:top-type}
Let $X$ be a smooth DM stack over $\mathbb C$ with projective coarse moduli and
$\bar{I}_{\mu}X:=\sqcup_{a\geq 1}\bar{I}_{\mu_{a}}X$ be the associated
rigidified (cyclotomic) inertia stack. For our purpose, a finer decomposition of
the rigidified inertia stack $\bar{I}_{\mu}X$ is needed. Given a set of line
bundles $L_{1},\cdots,L_{r}$ over $X$, we will choose a decomposition of the
rigidified inertia stack into open-closed components
$$\bar{I}_{\mu}X:= \bigsqcup_{c\in C} \bar{I}_{c}X$$
where $C$ is a \emph{finite} index set for the decomposition which satisfies the
following assumption:
\begin{assumption}\label{assump:inertiaindex}
  For each $c\in C$ and every $\mathbb C$-point $(x,g)\in \bar{I}_{c}X$, where
  $x$ is a $\mathbb C$-point of $X$ and $g$ is an element in the isotropy group
  $Aut(x)$, we will assume that the age function\footnote{Recall that if $g$
    acts on the fiber of $L$ over $x$ by multiplication by the number $\mathrm{e}^{q}$,
    where $q\in \mathbb Q$, then we define the age $age_{g}(L_{j}|_{x})$ to be
    $\langle{q}\rangle$, which is the fractional part of $q$.}
  $age_{g}L_{j}|_{x}$ and the order function $ord(g)$ are all constant functions
  on $\bar{I}_{c} X$. Then we can define $age_{c}L_{j}:=age_{g}L_{j}|_{x}$ and
  $a(c):=ord(g)$, which are well-defined as they are independent of the choice
  of the point $(x,g)\in \bar{I}_{c}X$.
  
  Moreover, we require that the index set $C$ hass an involution which is
  compatible with the involution $\iota: \bar{I}_{\mu}X \rightarrow
  \bar{I}_{\mu}X $ sending $(x,g)$ to $(x,g^{-1})$ for any $\mathbb C$-point
  $(x,g)$ of $\bar{I}_{\mu}X$; namely, for each index $c\in C$, there exists a
  unique index in $ C$ and we will denote it to be $c^{-1}$ such that
  $\iota(\bar{I}_{c}X )=\bar{I}_{c^{-1}}X $.
\end{assumption}

We will call an index set satisfying the above assumption a \emph{good index
  set}.

\begin{remark}
  We note that a good index set $C$ always exists (but may not unique), e.g., we
  can decompose $\bar{I}_{\mu}X$ into connected components, which corresponds to
  the biggest good index set. When $X$ is a quotient stack $[W/G]$ where $W$ is
  a connected affine scheme and $G$ is a finite group, then we can choose the conjugacy class $\mathrm{Conj}(G)$ of $G$ to
  be the index set $C$.
\end{remark}

When $Y$ is a substack of $X$, we will usually use $C$ to index
$\bar{I}_{\mu}Y$ as well, which is also a good index set with respect to the restriction
of line bundles $L_{1},\dots,L_{r}$ to $Y$.

Let $m$ be a nonnegative integer, for an (ordered) tuple
$\vec{m}=\{c_{1},\cdots,c_{m}\}\in C^{m}$, we will denote
$$\mathcal K_{g,\vec{m}}(X,\beta):=\mathcal K_{g,m}(X,\beta)\cap
\bigcap_{i=1}^{m} ev_{i}^{-1}(\bar{I}_{c_{i}}X)\ .$$ Here $\mathcal
K_{g,m}(X,\beta)$ is the moduli stack of genus $g$ twisted stable maps to $X$
with $m$ (not necessarily trivialized) gerby-marked points of degree $\beta$ as
defined in~\cite{abramovich2002compactifying} and $ev_{i}: \mathcal K_{g,m}(X,\beta)\rightarrow \bar{I}_{\mu}X$
is the evaluation map at the $i$th marking. Then $\mathcal K_{g,m}(X,\beta)$ can be
written as a disjoint union
$$ \mathcal
K_{g,m}(X,\beta):= \bigsqcup_{\vec{m}\in C^{m}} \mathcal K_{g,\vec{m}}(X,\beta)\
. $$

Sometimes, we will use the notation $[m]\cup \star$ (resp. $\vec{m}\cup \star$)
to mean $[m+1]$ (resp. $\overrightarrow{m+1}$) to distinguish the $m+1$-th
marking, which we denote to be $\star$.

Now we introduce \emph{admissible pairs}:
\begin{definition}\label{adaptedpair}
  Given a set of line bundles $L_{1},\cdots,L_{r}$ over $X$, a positive integer
  $N$ and $r$ weights $\overrightarrow{w_{j}}=(w_{1j},\cdots,w_{Nj})\in \mathbb
  Q^{N}_{\geq 0}$ corresponding to each $L_{j}$. Let $C$ be a good index set.
  Let $c\in C$, $\vec{k}=(k_{1},\cdots,k_{N})\in \mathbb Z_{\geq 0}^{N}$ and
  $\beta\in \mathrm{Eff}(X)$, where $\mathrm{Eff}(X)$ is the cone of effective
  curve classes of $X$. we will call $(\beta,\vec{k},c)$ an admissible pair if
  and only if
$$age_{c}(L_{j})\equiv \beta(L_{j})+(\overrightarrow{w_{j}},\vec{k})\; \mathrm{mod}\;
\mathbb Z $$ for all line bundles $L_{j}$, where
$(\overrightarrow{w_{j}},\vec{k})=\sum_{i}w_{ij}k_{i}$. We will denote
$\mathrm{Adm}$ to be the set of all admissible pairs.
\end{definition}





\begin{definition}\label{def:adm-series}
  Let $$f(z)=\sum_{(\beta,\vec{k},c)\in \mathrm{Eff}(X)\times \mathbb Z_{\geq
      0}^{N}\times C}q^{\beta}t_{1}^{k_{1}}\cdots
  t_{N}^{k_{N}}f_{\beta,\vec{k},c}(z)\ ,$$ be a formal series in
  $H^{*}(\bar{I}_{\mu}X,\mathbb
  C)[z,z^{-1}]][\![t_{1},\cdots,t_{N}]\!][\![\mathrm{Eff}(X)]\!]$ in
  which $$f_{\beta,\vec{k},c}\in H^{*}(\bar{I}_{c^{-1}}X,\mathbb C)[z,z^{-1}]]\
  .$$ We further put a (not unique!) $\mathbb Z_{2}(=\mathbb Z/ 2 \mathbb Z)$-grading on each
  variable $t_{i}$, which we denote the grading to be $\bar{t}_{i}$. Note that we also require the super-commutativity:$$t_{i}t_{j}=(-1)^{\bar{t}_{i}\cdot
    \bar{t}_{j}}t_{j}t_{i}\ .$$
  We will put a $Z_{2}$-grading on $H^{*}(\bar{I}_{\mu}X,\mathbb C)$, which is
  induced from its ordinary cohomological grading, and put even grading on Novikov variables $q^{\beta}$ and $z$.

  We will call $f(z)$ an admissible series if $f(z)$ satisfies the following
  condition.
  \begin{enumerate}
  \item $f_{\beta,\vec{k},c}(z)=0$ whenever $(\beta,\overrightarrow{k},c)$ is
    not an admissible pair;
  \item $f_{0,\overrightarrow{0},c}(z)=0$ for all $c\in C$.
  \item $f_{\beta,\vec{k},c}(z)$ belongs to the space
    $H^{\sum_{i}k_{i}\bar{t}_{i}}(\bar{I}_{c^{-1}}X,\mathbb C)[z]$, i.e., when
    $\sum_{i}k_{i}\bar{t}_{i}=0\in \mathbb Z_{2}$,
    $H^{\sum_{i}k_{i}\bar{t}_{i}}(\bar{I}_{c^{-1}}X,\mathbb C)$ means the even
    degree part of the
    cohomology group of $\bar{I}_{c^{-1}}X$; when
    $\sum_{i}k_{i}\bar{t}_{i}=1\in \mathbb Z_{2}$,
    $H^{\sum_{i}k_{i}\bar{t}_{i}}(\bar{I}_{c^{-1}}X,\mathbb C)$ means the odd
    degree part of the cohomology group of $\bar{I}_{c^{-1}}X$.       
  \end{enumerate}
  Furthermore, let $g\in \mathbb C[z]$, we say a formal series $g+f(z)$ is an
  admissible series near $g$ if $f(z)$ is an admissible series as above.

  For any pair $(\beta,\overrightarrow{k})$ in $\mathrm{Eff}(X)\times \mathbb
  Z_{\geq 0}^{N}$, we will
write $q^{\beta}\mathbf{t}^{\overrightarrow{k}}:=q^{\beta}t_{1}^{k_{1}}\cdots
t_{N}^{k_{N}}$ and call $(\beta,\overrightarrow{k})$ an extended degree.
\end{definition}

\begin{remark}
  Examples of admissible series include $J-$functions\cite{Givental_2004} and $I-$functions in
  quasimap theory\cite{Cheong_2015,
    ciocan-fontanine14_wall_cross_genus_zero_quasim,
    webb18_abelian_nonab_corres_i,
    webb21_abelian_quant_lefsc_orbif_quasim_i_funct,
    wang19_mirror_theor_gromov_witten_theor_without_convex}. Therefore
  admissible series consists of a large class of objects studied in
  Gromov-Witten theory.
\end{remark}

\subsection{Background on orbifold Gromov-Witten theory}
Now assume that $X$ carries a algebraic torus $T$ action (can be trivial), we
can define the so-called {\it Chen-Ruan cohomology} of $X$, Given any two
elements $\alpha_{1},\alpha_{2}$ in the $T-$equivariant
$$H^*_{\mathrm{CR}, T} (X,  \mathbb C):=H ^*_T(\bar{I}_{\mu}X , \mathbb C ) \ ,$$
We can define the Poincar\'e pairing in the {\it non-rigidified } inertia stack
$I_\mu X$ of $X$:
$$\langle \alpha_{1} , \alpha_{2} \rangle _{\mathrm{orb}} := \int _{\sum  _{c\in
    C} a(c)^{-1} [\bar{I}_{c}X]} \alpha_{1} \cdot \iota ^* \alpha_{2} \ .$$ Here
$\iota$ is the involution of $\bar{I}_{\mu}X $ obtained from the inversion
automorphisms of the band. Therefore, the diagonal class $[\Delta _{\bar{I}_{c}
  X}] $ obtained via push-forward of the fundamental class by $(\mathrm{id},
\iota): \bar{I}_{c}X \rightarrow \bar{I}_{c}X \times \bar{I}_{c^{-1}}X$ can be
written as
$$\sum _{c\in C} a(c)[ \Delta _{\bar{I}_{c}X}]
= \sum _{\alpha} \phi_{\alpha} \otimes \phi^{\alpha} \text{ in }
H^*_{T}(\bar{I}_{\mu}X \times \bar{I}_{\mu}X, \mathbb C), $$ where $\{
\phi_{\alpha}\}$ is a basis of $H^*_{\mathrm{CR}, T} (X, \mathbb C)$ with
$\{\phi^{\alpha}\}$ the dual basis with respect to the Poincar\'e pairing
defined above. Set $g_{\alpha\beta} = \langle { \phi_\alpha,\phi_\beta
}\rangle_{\mathrm orb}$ and $g^{\alpha\beta} = \langle
{\phi^\alpha,\phi^\beta}\rangle_{\mathrm orb}$.

Denote by $\bar{\psi}_i$ the first Chern class of the universal cotangent line
whose fiber at $((C, q_1, ..., q_m), [x])$ is the cotangent space of the coarse
moduli $\underline{C}$ of $C$ at $i$-th marking $\underline{q} _i$. For
non-negative integers $a_i$ and classes $\alpha _i \in H ^*_T (\bar{I}_{\mu}X ,
\mathbb C)$, using the virtual cycle $ [\mathcal
K_{0,\vec{m}}(X,\beta)]^{\mathrm{vir} }$ defined in~\cite{li1998virtual,
  behrend97_intrin_normal_cone, Dan_Abramovich_2008}, we write the Gromov-Witten
invariant:
\begin{align*}
  \langle \alpha _1\bar{\psi} ^{a_1}, ..., \alpha _m
  \bar{\psi}^{a_m} \rangle^{X} _{g, \vec{m}, \beta} & :=
                                                      \int _{[\mathcal K_{g,\vec{m}}(X,\beta)]^{\mathrm{vir} }}  \prod _i  ev _i ^* (\alpha _i) \bar{\psi} ^{a_i}_i \ .
\end{align*}
When the tuple $(g,\vec{m},\beta)$ gives rise to an empty stack $\mathcal
K_{g,\vec{m}}(X,\beta)$, we define the above integral is zero.

We will also need the stablemap Chen-Ruan classes
\begin{equation} \label{tilde-ev}(\widetilde {ev _j})_* = \iota
  _*(\boldsymbol{r} _{j}(ev_{j})_*),\end{equation} where $\boldsymbol{r} _{j}$
is the order function of the band of the gerbe structure at the marking $q_{j}$.
Define a class in $H_*^T(\bar{I}_{\mu}X )\cong H^*_T(\bar{I}_{\mu}X )$ by
\begin{align*}
  \langle \alpha _1,..., \alpha _m, - \rangle
  ^{X}_{0, \beta} &:=(\widetilde{ev_{m+1}})_* \left((\prod ev _i ^* \alpha _i) \cap [Q^{\epsilon}_{0, m}(X, \beta)]^{\mathrm{vir} }\right)\\
                  &=\sum _{\alpha}\phi^{\alpha} \langle\alpha_1,...,\alpha_m,\phi_{\alpha}\rangle^{X}_{0,m+1,\beta}\ .
\end{align*}

\section{A recursive relation about the Lagrangian
  cone}\label{sec:prop-Lagrangian}
\subsection{A root-stack modification of the twisted graph
  space}\label{subsec:space1}



Let $Y$ be any smooth Deligne-Mumford stack with projective coarse moduli. Let
$L$ be a semi-positive line bundle over $Y$, i.e., $\beta(L)\geq 0$ for any
degree $\beta\in \mathrm{Eff}(Y) $. Denote by $\P Y_{r,s}$ the root stack of the
projective bundle (also named twisted graph space\footnote{The terminology of
  ``twisted graph space'' is taken from
  \cite{clader17_higher_genus_quasim_wall_cross_via_local,
    clader17_higher_genus_wall_cross_gauged}, see loc.cits for some other
  applications of the twisted graph space in Gromov-Witten theory. }) $\P_{Y}(L^{\vee}\oplus \mathbb C)$
over $Y$ by taking $s-$th root of the zero section $D_{0}:=\P(0 \oplus \mathbb
C)$ and $r-$th root of the infinity section $D_{\infty}:=\P(L^{\vee}\oplus 0)$.
Let $\mathcal D_{0}$ and $\mathcal D_{\infty}$ be the corresponding root
divisors of $D_{0}$ and $D_{\infty}$ respectively, then the zero section
$\mathcal D_{0}\subset \P Y_{r,s}$ is isomorphic to the root stack
$\sqrt[s]{L^{\vee}/Y}$ with normal bundle isomorphic to the root bundle
$(L^{\vee})^{\frac{1}{s}}$, and the infinity section $\mathcal D_{\infty}\subset
\P Y_{r,s}$ is isomorphic to the root stack $\sqrt[r]{L/Y}$ with normal bundle
isomorphic to the root bundle $L^{\frac{1}{r}}$.

There are two morphisms
$$q_{0},q_{\infty}: \mathbb PY_{r,s} \rightarrow \mathbb B \C^{*}  \ ,$$
associated to the line bundles $\mathcal O(\mathcal D_{0})$ and $\mathcal
O(\mathcal D_{\infty})$ respectively such that $\mathcal O(\mathcal D_{0})\cong
q_{0}^{*}(\mathbb L)$ and $\mathcal O(\mathcal D_{\infty})\cong
q_{\infty}^{*}(\mathbb L)$, where $\mathbb L$ is the universal line bundle over
$\mathbb B \C^{*}$. This will induce morphisms on their rigidified inertia stack
counterparts:
$$ q_{0},q_{\infty}: \bar{I}_{\mu}\mathbb PY_{r,s} \rightarrow \bar{I}_{\mu} \mathbb B \C^{*}  \ ,$$
Note the rigidified inertia stack $ \bar{I}_{\mu} \mathbb B \C^{*}$ of $\mathbb
B \C^{*}$ can be written as the disjoint union $$ \bigsqcup_{c\in \C^{*}}
\bar{I}_{c} \mathbb B \C^{*} \ .$$Let $c$ be a complex number, any $\mathbb
C$-point of $ \bar{I}_{c} \mathbb B \C^{*}$ is isomorphic to the pair
$(1_{\mathbb C},c)$, where $1_{\mathbb C}$ is the trivial principal $\mathbb
C$-bundle over $\mathbb C$ and $c$ is an element in the automorphism group
$\mathrm{Aut}_{\mathbb C}(1_{\mathbb C})\cong \C^{*}$. Let $(y,g)$ be a $\mathbb
C-$point of $\bar{I}_{\mu} \P Y_{r,s}$, where $y\in Ob(\P Y_{r,s}(\mathbb C) )$
and $g\in \mathrm{Aut}(y)$. Assume that $q_{0}((y,g))=(1_{\mathbb C}, c)$, then
$g$ acts on the fiber $\mathcal O(\mathcal D_{0})|_{y}$ via the multiplication
by $c$. We have a similar relation for $q_{\infty}$.

Let $\mathrm{pr}_{r,s}: \P Y_{r,s} \rightarrow Y$ be the projection to the base,
which is a composition of deroot-stackification and projection from
$\P_{Y}(L^{\vee}\oplus \mathbb C) $ to $Y$.

Now choose a good index set $C$ for $\bar{I}_{\mu}Y$ satisfying the assumption
\ref{assump:inertiaindex}. For each $c\in C$, $c_{0},c_{\infty}\in \C^{*}$, we
will define the rigidified inertia component $ \bar{I}_{(c,c_{0},c_{\infty})} \P
Y_{r,s}$ to be
$$
\mathrm {pr}_{r,s}^{-1}(\bar{I}_{c}Y)\cap q_{0}^{-1} (\bar{I}_{c_{0}}\mathbb B
\C^{*} ) \cap q_{\infty}^{-1} (\bar{I}_{c_{\infty}}\mathbb B \C^{*}) \ .$$ We
also denote $(c,c_{0},c_{\infty})^{-1}=(c^{-1},c^{-1}_{0},c^{-1}_{\infty})$.
Notice that $\bar{I}_{(c,c_{0},c_{\infty})}\P Y_{r,s}$ is nonempty if and only
if $s\cdot age_{c}(L)\equiv -c_{0}$ mod $\mathbb Z$ and $r\cdot
age_{c}(L)\equiv c_{\infty}$ mod $\mathbb Z$.

\begin{definition}\label{def:deg}
  For any degree $\beta\in \mathrm{Eff}(Y)$ and rational number $\delta\in
  \mathbb Q_{\geq 0}$, we say a stable map $f:C\rightarrow \P Y_{r,s}$ is of degree
  $(\beta,\frac{\delta}{r})$ if $(\mathrm{pr}_{r,s}\circ f)_{*}[C]=\beta$ and
  $deg(f^{*}\mathcal O(\mathcal D_{\infty}))=\frac{\delta}{r}$. We will denote
  $\mathcal K_{0,m}(\P Y_{r,s},(\beta,\frac{\delta}{r}))$ to be the
  corresponding moduli stack of stable maps to $\P Y_{r,s}$ of degree
  $(\beta,\frac{\delta}{r})$. Note that the line bundle $\mathcal O(\mathcal
  D_{\infty})$ is semi-positive.
\end{definition}

\subsection{Localization analysis}\label{sec:lc1}
Consider the $\C^{*}-$action on the projection bundle $\mathbb P_{Y}(L^{\vee}\oplus
\mathbb C)$ by scaling the fiber such that the normal bundle of $D_{0}$ (resp.
$D_{\infty}$) in $\mathbb P_{Y}(L^{\vee}\oplus \mathbb C)$ is of $\C^{*}-$weight 1
(resp. $-1$). This $\C^{*}-$action induces a $\C^{*}-$action on $\P Y_{r,s}$
such that the normal bundle of $\mathcal D_{0}$ (resp. $\mathcal D_{\infty}$) in $\P
Y_{r,s}$ is of $\C^{*}-$weight $\frac{1}{s}$ (resp. $-\frac{1}{r}$). Then we
have an induced $\C^{*}-$action on the moduli $\mathcal
K_{0,\overrightarrow{m}}(\P Y_{r,s}, (\beta,\frac{\delta}{r}))$ of twisted stable maps to $\P
Y_{r,s}$. We will apply the virtual localization formula of
Graber--Pandharipande~\cite{Graber_1999} to $\mathcal
K_{0,\overrightarrow{m}}(\P Y_{r,s}, (\beta,\frac{\delta}{r}))$, for which
introduce the notation of decorated graph so that we index the components of
$\C^{*}-$fixed loci of $\mathcal K_{0,\overrightarrow{m}}(\P
Y_{r,s},(\beta,\frac{\delta}{r}))$ by decorated graphs (trees). First, a decorated graph
$\Gamma$ consists of vertexes, edges and legs with the following decorations:
\begin{itemize}
\item Each vertex $v$ is associated with an index $j(v) \in \{0, \infty\}$, and
  a degree $\beta(v) \in \mathrm{Eff}(Y)$.
\item Each edge $e$ consists of a pair of half-edges $\{h_{0},h_{\infty}\}$, and
  $e$ is equipped with a degree $\delta(e) \in \mathbb{Q}_{>0}$. Here we call
  $h_{0}$ and $h_{\infty}$ half edges, and $h_{0}$ (resp. $h_{\infty}$) is
  attached to a vertex labeled by $0$ (resp. $\infty$).
\item Each half-edge $h$ and each leg $l$ has a multiplicity $m(h)$ or $m(l)$ in
  $C \!\times\! \C^{*} \!\times\! \C^{*} $.
\item The legs are labeled with the numbers $\{1, \ldots, m\}$ and each leg is
  incident to a unique vertex.
\end{itemize}

By the ``valence" of a vertex $v$, denoted $\text{val}(v)$, we mean the total
number of incident half-edges and legs.

For each $\C^{*}-$fixed stable map $f:(C; q_1, \ldots, q_m)\rightarrow \P Y_{r,s}$ in $\mathcal K_{0,\overrightarrow{m}}(\P
Y_{r,s},(\beta,\frac{\delta}{r}))$, we can associate a decorated graph $\Gamma$
in the following way.
\begin{itemize}
\item Each edge $e$ corresponds to a genus-zero component $C_e$ which maps
  constantly to the base $Y$ and intersects $\mathcal D_{0}$ and $\mathcal
  D_{\infty}$ properly. The restriction $f|_{C_{e}}$ satisfies that $r\cdot
  deg(f|_{C_{e}}^{*}\mathcal O(\mathcal D_{\infty}) )= s\cdot
  deg(f|_{C_{e}}^{*}\mathcal O(\mathcal D_{0}) )= \delta(e)$.
\item Each vertex $v$ for which $j(v) = 0$ (with unstable exceptional cases
  noted below) corresponds to a maximal sub-curve $C_v$ of $C$ which maps
  totally into $\mathcal D_{0}$, then the restriction of $f$ to $C_{v}$ defines a twisted stable map in
  $$\mathcal
  K_{0,val(v)}(\sqrt[s]{L^{\vee}/Y},\beta(v))\ .$$ Each vertex $v$ for which
  $j(v) = \infty$ (again with unstable exceptions) corresponds to a maximal
  sub-curve which maps totally into $\mathcal D_{\infty}$, then the restriction
  of $f$ to $C_{v}$
  defines a twisted stable map in
  $$\mathcal
  K_{0,val(v)}(\sqrt[r]{L/Y},\beta(v))\ .$$ The label $\beta(v)$ denotes the
  degree coming from the restriction $\mathrm{pr}_{r,s}\circ f|_{C_{v}}:C_{v}\rightarrow Y$.
\item A vertex $v$ is {\it unstable} if stable twisted maps of the type
  described above do not exist. In this case, we have $\beta(v)=0$ and $v$ corresponds to a single point
  of the component $C_e$ for each incident edge $e$, which may be a node at
  which $C_e$ meets another edge curve $C_{e'}$, a marked point of $C_e$, or an
  unmarked point.
\item Each leg $l$ labeled by $i$ corresponds to the $i$th marking $q_{i}$, and it's
  incident to the vertex $v$ if $q_{i}$ lies on $C_{v}$. The index $m(l)$ on a leg $l$ indicates the rigidified inertia stack
  component $\bar{I}_{m(l)}\P Y_{r,s}$ of $\P Y_{r,s}$ on which the gerby marked
  point corresponding to the leg $l$ is evaluated.
\item A half-edge $h$ of an edge $e$ corresponds a ramification point $q \in
  C_{e}$, i.e., there are two distinguished points(we call them ramification
  points) $q_{0}$ and $q_{\infty}$ on $C_{e}$ satisfying that $q_{0}$ maps to
  $\mathcal D_{0}$ and $q_{\infty}$ maps to $\mathcal D_{\infty}$, respectively.
  We associate half-edges $h_{0}$ and $h_{\infty}$ to $q_{0}$ and $q_{\infty}$
  respectively. Then $m(h)$ indicates the rigidified inertia component
  $\bar{I}_{m(h)}\P Y_{r,s}$ of $\P Y_{r,s}$ on which the ramification point $q$
  associated with $h$ is evaluated. We will write the edge $e$ as
  $\{h_{0},h_{\infty}\}$.
\end{itemize}

Moreover, the decorated graph from a $\C^{*}-$fixed stable map should satisfy
the following:
\begin{enumerate}\label{prop:edge1}
\item(Monodromy constraint for edge) For each edge $e=\{h_{0},h_{\infty}\}$, we
  have that $m(h_{0})=(c^{-1}, \mathrm{e}^{\frac{\delta(e)}{s}},1)$ and
  $m(h_{\infty})=(c, 1, \mathrm{e}^{\frac{\delta(e)}{r}})$ for some $c\in C$.
  Moreover the choice of $c$ should satisfy that $$\delta(e)= age_{c}(L)\;
  \text{mod}\; \mathbb Z\ ,$$ which follows from the fact that the ages of line
  bundles $f|^{*}_{C_{e}}\mathcal O(s \mathcal D_{0})$ and
  $f|^{*}_{C_{e}}(\mathcal O(r\mathcal D_{\infty})\otimes
  \mathrm{pr}_{r,s}^{*}L^{\vee})$ at $q_{\infty}$ should be equal, as the two
  line bundles on $C_{e}$ are isomorphic.
\item(Monodromy constraint for vertex) For each vertex $v$, let $I_{v}\subset
  [m]$ be the set of incident legs to $v$, and $H_{v}$ be the set of incident
  half-edges to $v$. For each $i\in I_{v}$, write
  $m(l_{i})=(c_{i},a_{i},b_{i})$, and for each $h\in H_{v}$, write
  $m(h)=(c_{h},a_{h},b_{h})$. If $v$ is labeled by $0$, we have that all $b_{i}$
  and $b_{h}$ are equal to $1$, and
  \begin{align}\label{eq:monver0}
    \mathrm{e}^{\frac{\beta(v)(L)}{s}}\times \prod_{i\in I_{v}} a_{i} \times \prod_{h\in H_{v}}
    a_{h}^{-1}=1\ .
  \end{align}
  If $v$ is labeled by $\infty$, then we have all $a_{i}$ and $a_{h}$ are all
  equal to $1$ and
  \begin{align}\label{eq:monverinfty}
    \mathrm{e}^{-\frac{\beta(v)(L)}{r}}\times \prod_{i\in I_{v}} b_{i} \times
    \prod_{h\in E_{v}}b^{-1}_{h}=1\ .
  \end{align}
\end{enumerate}

We will only consider decorated graph coming from a $\C^{*}-$fixed stable map. In particular, we note that the decorations at each stable vertex $v$ yield a
tuple
$$\overrightarrow{val}(v) \in (C \!\times\! \C^{*} \!\times\! \C^{*})^{\text{val}(v)}$$
recording the multiplicities at every special point of $C_v$.\footnote{For each
  node of $C_{v}$, let $h$ be the incident half-edge, then we define the
  multiplicity at the branch of node at $C_{v}$ to be $m(h)^{-1}$.} Then the
restriction gives a stable map in
$$ \mathcal
K_{0,\overrightarrow{val}(v)}(\mathcal D_{j(v)},\beta(v))\ .$$

For each decorated graph $\Gamma$, we will associate each vertex $v$ (resp. edge
$e$) a moduli space $\mathcal M_{v}$(resp. $\mathcal M_{e}$) over which there is
a family of $\C^{*}-$stable maps to $\P Y_{r,s}$ with decorated degree(there are
exceptions for unstable vertexes $v$, see section \ref{subsubsec:ver-cntr1} for more
details). Denote by $F_{\Gamma}$ the space
$$\prod_{v:j(v)=0}\mathcal M_{v}\times_{\bar{I}_{\mu}\mathcal D_{0}}\prod_{e\in E}\mathcal
M_{e}\times_{\bar{I}_{\mu}\mathcal D_{\infty}} \prod_{v:j(v)=\infty} \mathcal
M_{v} \ ,$$ where the fiber product is taken by gluing the two branches at each
node. By virtual localization formula~\cite{Graber_1999}, we can write
$$[\mathcal K_{0,\overrightarrow{m}}(\P Y_{r,s},\beta )]^{ \mathrm{vir}}\ ,$$
in terms of contributions from each decorated graph $\Gamma$:
\begin{equation}
  \label{localization}
  [\mathcal K_{0,\overrightarrow{m}}(\P Y_{r,s},\beta )]^{\mathrm{vir}} = \sum_{\Gamma} \frac{1}{\mathbb A_{\Gamma}}\iota_{\Gamma*} \left(\frac{[F_{\Gamma}]^{\mathrm{vir}}}{e^{\C^{*}}(N_{\Gamma}^{\mathrm{vir}})}\right)\ .
\end{equation}
Here, for each graph $\Gamma$, $[F_{\Gamma}]^{\mathrm{vir}}$ is obtained via the
$\C^*$-fixed part of the restriction to the fixed loci of the obstruction theory
of $\mathcal K_{0,\overrightarrow{m}\cup \star}(\P Y_{r,s},\beta )$, and
$N_{\Gamma}^{\mathrm{vir} }$ is the equivariant Euler class of the $\C^*$-moving
part of this restriction. Besides, $\mathbb A_{\Gamma}$ is the automorphism
factor for the graph $\Gamma$, which represents the degree of $F_{\Gamma}$ into
the corresponding open and closed $\C^{*}$-fixed substack
$i_{\Gamma}(F_{\Gamma})$ in $\mathcal K_{0,\overrightarrow{m}\cup \star}(\P
Y_{r,s},\beta )$. In our case, $\mathbb A_{\Gamma}$ is the product of the size
of the automorphism group $\mathrm{Aut}(\Gamma)$ of $\Gamma$ and the degrees
from each edge moduli
$\mathcal M_{e}$ into the corresponding fixed loci.

\emph{Assume that $r,s$ are sufficiently large primes.} We will do an explicit computation for the contributions of each graph $\Gamma$
in the following. As for the contribution of a graph $\Gamma$ to
\eqref{localization}, one can first apply the normalization exact sequence to
the obstruction theory, which decomposes the contribution from $\Gamma$ to
\eqref{localization} into contributions from vertex, edge, and node factors.

\subsubsection{Vertex contribution}\label{subsubsec:ver-cntr1}
Assume that $v$ is a stable vertex. If the vertex is labeled by $\infty$, we
define the vertex moduli $\mathcal M_{v}$ to be $\mathcal
K_{0,\overrightarrow{val}(v)}(\mathcal D_{\infty},\beta(v))$ and the fixed part
of perfect obstruction theory gives rise to $[\mathcal
K_{0,\overrightarrow{val}(v)}(\mathcal D_{\infty},\beta(v))]^{\mathrm{vir}}$.
Let $\pi:\mathcal C\rightarrow \mathcal K_{0,\overrightarrow{val}(v)}(\mathcal
D_{\infty},\beta(v))$ be the universal curve and $f:\mathcal C\rightarrow
\mathcal D_{\infty}$ be the universal map, then the movable part of the perfect
obstruction theory yields the \emph{inverse of the Euler class} of the virtual
normal bundle which is equal to
$$e^{\C^{*}}((-R^{\bullet}\pi_{*}f^{*} L^{\frac{1}{r}})\otimes
\C_{-\frac{\lambda}{r}} )\ .$$ When $r$ is a sufficiently large prime and the multiplicity
$m(l)$ corresponding to each leg $l$ incident to $v$ is equal to $(c_{l},1,\frac{\delta_{l}}{r})$
for some prefixed $\delta_{l}\in \mathbb Q_{\geq 0}$(note this implies
$\delta_{l}\ll r$) and $c_{l}\in C$, following
a generalization of~\cite{janda2020double} to the orbifold case, the above Euler
class has a representation
$$\sum_{d\geq 0}c_{d}(-R^{\bullet}\pi_{*}f^{*}L^{\frac{1}{r}})(\frac{-\lambda}{r})^{|E(v)|-1-d} \ .$$ Here the
virtual bundle $-R^{\bullet}\pi_{*}f^{*}L^{\frac{1}{r}}$ has virtual rank
$|E(v)|-1$, where $|E(v)|$ is the number of edges incident to the vertex $v$.

If the vertex $v$ is labeled by $0$, we define the vertex moduli $\mathcal
M_{v}$ to be $\mathcal K_{0,\overrightarrow{val}(v)}(\mathcal D_{0},\beta(v))$
and the fixed part of perfect obstruction theory gives rise to $[\mathcal
K_{0,\overrightarrow{val}(v)}(\mathcal D_{0},\beta(v))]^{\mathrm{vir}}$. Let
$\pi:\mathcal C\rightarrow \mathcal K_{0,\overrightarrow{val}(v)}(\mathcal
D_{0},\beta(v))$ be the universal curve and $f:\mathcal C\rightarrow \mathcal
D_{0}$ be the universal map, then the movable part of the perfect obstruction
theory yields the \emph{inverse of the Euler class} of the virtual normal bundle
which is equal to
$$e^{\C^{*}}((-R^{\bullet}\pi_{*}f^{*} (L^{\vee})^{\frac{1}{s}})\otimes
\C_{\frac{\lambda}{s}})\ .$$ Now assume that there is only one edge incident to
$v$ and the multiplicity
$m(l)$ corresponding to each leg $l$ incident to $v$ is equal to $(c_{l},\frac{\delta_{l}}{s},1)$
for some prefixed $\delta_{l}\in \mathbb Q_{\geq 0}$(note this implies $\delta_{l}\ll r$) and $c_{l}\in G$. By~\cite[Lemma 5.2, Remark
5.3]{wang19_mirror_theor_gromov_witten_theor_without_convex}, the above Euler
class is equal to $1$ when $\beta(v)(L)>0$, we note here the semi-positivity of
$L$ is essentially used in the proof of\cite[Lemma
5.2]{wang19_mirror_theor_gromov_witten_theor_without_convex}. On the other hand,
when $\beta(v)(L)=0$, by the argument in the proof of~\cite[Lemma
6.5]{wang19_mirror_theor_gromov_witten_theor_without_convex}, the above Euler
class is still equal to $1$ there is a special point (we will choose a node in
our case) on $C_{v}$ whose multiplicity $m$ satisfies that
$age_{m}\mathcal O(\mathcal D_{0})\neq 0$.

When $v$ is an unstable vertex
over $0$ (resp. $\infty$), let $h$ be the half-edge incident to $v$, the vertex
moduli $\mathcal M_{v}$ is defined to be $\bar{I}_{m(h)^{-1}}\mathcal
D_{0}$(resp. $\bar{I}_{m(h)^{-1}}\mathcal D_{\infty}$) with $[\mathcal
M_{v}]^{\mathrm{vir}}=[\mathcal M_{v}]$ and zero virtual normal bundle.

\subsubsection{Edge contribution}\label{subsubsec:edge-cntr1}
Let $e=\{h_{0},h_{\infty}\}$ be an edge in $\Gamma$ with decorated degree
$\delta(e)\in \mathbb Q$(we will write $\delta(e)$ as $\delta$ for simplicity if
no confusion occurs) and decorated multiplicities $m(h_{0})$ and
$m(h_{\infty})$. Assume that $m(h_{\infty})=(c, 1,
\mathrm{e}^{\frac{\delta(e)}{r}})$, denote $a_{e}:=a(c)$ be the number as
Assumption \ref{assump:inertiaindex}. We define the edge moduli $\mathcal M_{e}$
to be the root gerbe $\sqrt[as\delta(e)]{L^{\vee}/I_{c}Y}$ over $I_{c}Y $, then the
virtual cycle $[\mathcal M_{e}]^{vir}$ coming from the fix part of the
obstruction theory is equal to the fundamental class $[\mathcal M_{e}]$ of
$\mathcal M_{e}$ and
\emph{inverse of the Euler class} of virtual normal bundle is equal to $1$ when
$r$ is a sufficiently large prime. We
note that $\mathcal M_{e}$ allows a finite \'etale map into the corresponding fixed-loci in $\mathcal
K_{0,2}(\P Y_{r,s}, (0,\frac{\delta(e)}{r}))$ of degree $\frac{1}{as}$. See
appendix \ref{sec:app-edge} for more details.

\subsubsection{Node contributions}\label{subsubsec:node-cntr1}
The deformations in $\mathcal K_{0,\overrightarrow{m}}(\P
Y_{r,s},(\beta,\frac{\delta}{r})$ smoothing a node contribute to the Euler class
of the virtual normal bundle as the first Chern class of the tensor product of
the two cotangent line bundles at the branches of the node. For nodes at which a
component $C_e$ meets a component $C_v$ over the vertex labeled by $X_{0}$, this
contribution is
\begin{equation*}
  \frac{\lambda - c_{1}(L)}{a_{e}s\delta(e)} - \frac{\bar{\psi}_v}{a_{e}s}.
\end{equation*}

For nodes at which a component $C_e$ meets a component $C_v$ at the vertex over
$\mathcal D_{\infty}$, this contribution is
\begin{equation*}
  \frac{-\lambda + c_{1}(L)}{a_{e}r\delta(e)} -\frac{\bar{ \psi}_v}{a_{e}r}.
\end{equation*}

We will not need node contributions from other types of nodes as the above types
suffice the need in this paper. 

\subsubsection{Total contribution}\label{subsec:total1}
For any decorated graph $\Gamma$, we define $F_{\Gamma}$ to be the fiber product
$$\prod_{v:j(v)=0}\mathcal M_{v}\times_{\bar{I}_{\mu}\mathcal D_{0}}\prod_{e\in E}\mathcal
M_{e}\times_{\bar{I}_{\mu}\mathcal D_{\infty}} \prod_{v:j(v)=\infty} \mathcal
M_{v}
$$ of the following
diagram:
$$\xymatrix{
  F_{\Gamma}\ar[r]\ar[d] &\prod\limits_{v:j(v)=0} \mathcal M_{v}\times
  \prod\limits_{e\in E}\mathcal M_{e}\times \prod\limits_{v:j(v)=\infty}
  \mathcal M_{v}
  \ar[d]^-{ev_{nodes}}\\
  \prod\limits_{E}(\bar{I}_{\mu}\mathcal D_{0}) \times \bar{I}_{\mu}\mathcal
  D_{\infty}\ar[r]^-{(\Delta^{0}\times \Delta^{\infty})^{|E|}}
  &\prod\limits_{E}(\bar{I}_{\mu}\mathcal D_{0})^{2}\times
  (\bar{I}_{\mu}\mathcal D_{\infty})^{2}\ , }
$$
where $\Delta^{0}=(id,\iota)$(resp. $\Delta^{\infty}=(id,\iota)$) is the
diagonal map of $\bar{I}_{\mu}\mathcal D_{0}$ (resp. $\bar{I}_{\mu}\mathcal
D_{\infty}$). Here the right-hand
vertical map $ev_{nodes}$ is the product of the evaluation maps at the two
branches of each node. We note that there are two nodes corresponding to $h_{0}$
and $h_{\infty}$ for each edge $e=\{h_{0},h_{\infty}\}$.

We define $[F_{\Gamma}]^{\mathrm{vir}}$ to be:
\begin{equation*}
  \begin{split}
    \prod_{v:j(v)=0}[\mathcal M_{v}]^{\mathrm{vir}}
    \times_{\bar{I}_{\mu}\mathcal D_{0}}\prod_{e\in E}[\mathcal
    M_{e}]^{\mathrm{vir}}\times_{\bar{I}_{\mu}\mathcal D_{\infty}}
    \prod_{v:j(v)=\infty}[\mathcal M_{v}]^{\mathrm{vir}}\ .
  \end{split}
\end{equation*}
Then the contribution of decorated graph $\Gamma$ to the virtual localization
is:
\begin{equation}\label{eq:auto-loc1}
  Cont_{\Gamma}=\frac{\prod_{e\in
      E}sa_{e}}{|\text{Aut}(\Gamma)|}(\iota_{\Gamma})_{*}\left(\frac{[F_{\Gamma}]^{\mathrm{vir}}}{e^{\C^{*}}(N^{\mathrm{vir}}_{\Gamma})}\right)\ .
\end{equation}
Here $\iota_{F}:F_{\Gamma}\rightarrow \mathcal K_{0,\overrightarrow{m}}(\P
Y_{r,s},(\beta,\frac{\delta}{r}))$ is a finite etale map of degree
$\frac{|\text{Aut}(\Gamma)|}{\prod_{e\in E}sa_{e}}$ into the corresponding
$\C^{*}$-fixed loci. The virtual normal bundle
$e^{\C^{*}}(N^{\mathrm{vir}}_{\Gamma})$ is the product of virtual normal bundles
from vertex contributions, edge contributions and node contributions.

\subsection{Recursive relation}


Take $L_{1}=L$ and $r=1$ in \ref{def:adm-series} and write
$\vec{w}=\overrightarrow{w_{1}}$ to be the weight. Let
$\mu(z)=\sum_{\beta,\vec{k},c}q^{\beta}\mathbf{t}^{\vec{k}}\mu_{\beta,\vec{k},c}(z)$ be an admissible (power)
series as in \ref{def:adm-series}. We further require that
$\mu_{\beta,\vec{k},c}\in H^{*}(\bar{I}_{c^{-1}}X,\mathbb C)[z]$ is a polynomial
in $z$. For any admissible pair $(\beta,\vec{k},c)$ in $\mathrm{Eff}(Y)\times
\mathbb Z_{\geq 0}^{N} \times C$, define $J^{Y}_{\beta,\vec{k},c}(\mu,z)$ to be
\begin{equation}\label{eq:jcomp}
  \begin{split}
    &\mu_{\beta,\vec{k},c}(z)+\mathbf{Coeff}_{\mathbf{t}^{\vec{k}}}\bigg[\sum_{m\geq 0} \sum_{\substack{(\beta_{j},
        \overrightarrow{k_{j}},c_{j})_{j=1}^{m}\in \mathrm{Adm}^{m},
        \beta_{\star}\in \mathrm{Eff}(Y)\\ \beta_{1}+\cdots+\beta_{m}+\beta_{\star}=\beta\\ \overrightarrow{k_{1}}+\cdots+\overrightarrow{k_{m}}=\overrightarrow{k}  }}\\
    & \frac{1}{m!}\phi^{\alpha}\langle
    \mathbf{t^{\overrightarrow{k_{1}}}\mu_{\beta_{1},\overrightarrow{k_{1}},c_{1}}}(-\bar{\psi}_{1}),\cdots,\mathbf{t^{\overrightarrow{k_{m}}}\mu_{\beta_{m},\overrightarrow{k_{m}},c_{m}}}(-\bar{\psi}_{m}),\frac{\phi_{\alpha}}{z-\bar{\psi}_{\star}}\rangle^{Y}_{0,\overrightarrow{m}\cup\star,\beta_{\star}}\bigg]\
    ,
  \end{split}
\end{equation}
where $\overrightarrow{m}\cup\star =\big(c_{1}^{-1}\cdots,c_{m}^{-1},c\big)\in
C^{m+1}$, $\overrightarrow{k_{j}}=(k_{j1},\cdots,k_{jN})$ and
$\overrightarrow{k_{1}}+\cdots+\overrightarrow{k_{m}}=\overrightarrow{k}$ means
the component-wise addition in $\mathbb Z_{\geq 0}^{N}$.

\begin{remark}
  We can show that when $(\beta,\vec{k},c)$ is not an admissible pair and we
  define $J^{Y}_{\beta,\overrightarrow{k},c}(\mu,z)$ in the same way as above,
  we have $J^{Y}_{\beta,\overrightarrow{k},c}(\mu,z)=0$. Moreover
  $J_{0,\vec{0},c}=0$. This implies that
  $J-$function $$J^{Y}(q,\mu,z)=z+\sum_{\beta,\overrightarrow{k},
    c}q^{\beta}\mathbf{t}^{\vec{k}}
  J_{\beta,\overrightarrow{k},c}^{Y}(\mu,z)$$ associated with the input $\mu(z)$
  above is an admissible series near $z$.
\end{remark}

\begin{definition}\label{def:star-graph}
  Let $m,n$ be two nonnegative integers and $(\beta,\vec{k},c)$ be an admissible pair. We denote
  $\Lambda_{\beta,\vec{k},c,m,n}$ to be the set of tuples
$$ \big(c,
\beta_{\star},((\beta_{1},\overrightarrow{k_{1}},c_{1}),\cdots,(\beta_{m+n},\overrightarrow{k_{m+n}},c_{m+n}))\big)\in
C\times \mathrm{Eff}(Y) \times \mathrm{Adm}^{m+n}\ ,$$ where we require that
$\beta_{\star}+\sum_{i=1}^{m+n}\beta_{i}=\beta$,
$\sum_{i=1}^{m+n}\overrightarrow{k_{i}}=\overrightarrow{k}$, $\beta_{i}(L)+ (\vec{w},
\overrightarrow{k_{i}}) >0 $ for $ 1\leq i\leq m$ and $\beta_{i}(L)+
(\vec{w},\overrightarrow{k_{i}})=0 $ for $ m+1\leq i\leq m+n$. We call an
element of $\Lambda_{\beta,\vec{k},c,m,n}$ stable if $\beta_{\star}\neq 0$ or
$m+n\geq 2$ when $\beta_{\star}=0$.
\end{definition}
We note that $\Lambda_{\beta,\vec{k},c,m,n}$ is a finite set as $\mathcal
K_{0,m+n}(X,\beta)$ is finite type over $\mathbb C$, hence Noetherian.

For any admissible pair $(\beta,\vec{k},c)$ with $\beta(L)+(\vec{w},\vec{k})>0$. We have the following recursive characterization about
$J^{Y}_{\beta,\overrightarrow{k},c}(\mu,z)$.
\begin{proposition}\label{prop:chara-general}
For any admissible pair $(\beta,\vec{k},c)$ with $\beta(L)+(\vec{w},\vec{k})>0$.
Assume that $r$ is a sufficiently large prime. Then for any nonnegative integer $b$, we have the following recursive relation:
  \begin{equation}\label{eq:rec1}
    \begin{split}
      & [J^{Y}_{\beta,\vec{k},c}(\mu,z)]_{z^{-b-1}}\\
      &=\bigg[\sum_{m=0}^{\infty}\sum_{n=0}^{\infty} \sum_{\substack{\Gamma\in
          \Lambda_{\beta,\overrightarrow{k},c,m,n}\\ \Gamma\; \textrm{is
            stable}}}\frac{1}{m!n!}(\widetilde {ev_{\star}})_{*}\bigg(\sum
      _{d=0}^{\infty}\epsilon_{*}\big(c_{d}(-R^{\bullet}\pi_{*}f^{*}L^{\frac{1}{r}})(\frac{\lambda}{r})^{-1+m-d}(-1)^{d}\\
      &\cap [\mathcal K_{0,\overrightarrow{m+n}\cup
        \star}(\sqrt[r]{L/Y},\beta_{\star})]^{\mathrm{vir}}\big) \cap
      \prod_{i=1}^{m}\frac{ev_{i}^{*}(\frac{1}{\delta_i} \mathbf{t}^{\overrightarrow{k_{i}}}J^{Y}_{\beta_{i},\vec{k_i},c_{i}}(\mu,z)|_{z=\frac{\lambda-c_1(L)}{\delta_{i}}})}{\frac{\lambda-ev^{*}_{i}c_1(L)}{r\delta_{i}}+\frac{\bar{\psi_{i}}}{r}}\\
      &\cap
      \prod_{i=m+1}^{m+n}ev_{i}^{*}(\mathbf{t}^{\overrightarrow{k_{i}}}\mu_{\beta_{i},\vec{k_i},c_{i}}(-\bar{\psi}_{i}))
      \cap \bar{\psi}_{\star}^{b}\bigg)\bigg]_{\mathbf{t}^{\overrightarrow{k}}\lambda^{-1}}\ .
    \end{split}
  \end{equation}
  Here $\delta_{i}=\beta_{i}(L)+ (\vec{w}, \overrightarrow{k_{i}})$, and
  $\epsilon:\mathcal K_{0,\overrightarrow{m+n}\cup
    \star}(\sqrt[r]{L/Y},\beta_{\star})\rightarrow \mathcal
  K_{0,\overrightarrow{m+n}\cup \star}(Y,\beta_{\star})$ is the natural
  structural morphism by forgetting the root structure of $\sqrt[r]{L/Y}$
  (c.f.\cite{tang2021quantum}), where we choose the
  tuple $\overrightarrow{m+n}\cup\star$ for $\mathcal
  K_{0,\overrightarrow{m+n}\cup \star}(\sqrt[r]{L/Y},\beta_{\star})$ to be
  $$\big((c_{1}^{-1},1,\mathrm{e}^{\frac{-\delta_1}{r}}),\cdots,
  (c_{m+n}^{-1},1,\mathrm{e}^{\frac{-\delta_{m+n}}{r}}), (c,1,
  \mathrm{e}^{\frac{\delta}{r}})\big)\in (C \!\times\! \C^{*} \!\times\!
  \C^{*})^{m+n+1}\ ,$$ and choose the tuple $\overrightarrow{m+n}\cup\star$ for
  $\mathcal K_{0,\overrightarrow{m+n}\cup \star}(Y,\beta_{\star})$ to be
    $$(c_{1}^{-1},\cdots, c_{m+n}^{-1},c)\in C^{m+n+1}\ .$$   
  \end{proposition}

  The proof of the above proposition is based on applying virtual localization
  to the following integral.
  \begin{equation}\label{eq:mainintegral1}
    \begin{split}
      &\sum_{m=0}^{\infty}\sum_{\substack{ (\beta_{j},
          \overrightarrow{k_{i}},c_{i})_{i=1}^{m}\in \mathrm{Adm}^{m},
          \beta_{\star}\in \mathrm{Eff}(Y)\\
          \beta_{\star}+\sum_{i=1}^{m}\beta_{i}=\beta \\
          \overrightarrow{k_{1}}+\cdots+\overrightarrow{k_{m}}=\overrightarrow{k}
        }}\frac{1}{m!}(\widetilde{EV_{\star}})_{*}\bigg(
      \prod_{i=1}^{m}ev_{i}^{*}\big(\mathrm{pr}_{r,s}^{*}(\mathbf{t}^{\overrightarrow{k_{i}}}\mu_{\beta_{i},\vec{k_i},c_{i}}(-\bar{\psi}_{i}))\big)\\
      &\cap \bar{\psi}_{\star}^{b}\cap [\mathcal K_{0,\overrightarrow{m}\cup
        \star}(\P
      Y_{r,s},(\beta_{\star},\frac{\delta}{r}))]^{\mathrm{vir}}\bigg)\ .
    \end{split}
  \end{equation}
  Here an explanation of the notations is in order:
  \begin{enumerate}
 
  \item For degrees $\beta_{\star},\beta_{1},\cdots, \beta_{m}$ in
    $\mathrm{Eff}(X)$ and tuples $ \overrightarrow{k_{1}},\cdots,
    \overrightarrow{k_{m}}$ and with
    $\sum_{i=1}^{m}\beta_{i}+\beta_{\star}=\beta$ and
    $\overrightarrow{k_{1}}+\cdots+\overrightarrow{k_{m}}=\overrightarrow{k}$,
    let $(c_{1},\cdots,c_{m})\in C^{m}$ such that $(\beta_{i},c_{i})$ are
    admissible pairs. Write $\delta_{i}=
    \beta_{i}(L)+(\vec{w},\overrightarrow{k_{i}})$ and $\delta=\beta(L)+(\vec{w},\overrightarrow{k})$,  we define $\vec{m}\cup
    \star$ to be the $m+1$tuple
$$\big((c^{-1}_{1},\mathrm{e}^{\frac{\delta_i}{s}},1),\cdots,(c^{-1}_{m},\mathrm{e}^{\frac{\delta_m}{s}},1),(c,1,\mathrm{e}^{\frac{\delta}{r}})\big)\ . $$ 

\item the morphism $EV_{\star}$ is a composition of the following maps:
  $$\xymatrix{ \mathcal K_{0,\overrightarrow{m}\cup \star}(\P
    Y_{r,s},(\beta_{\star},\frac{\delta}{r}))\ar[r]^-{ev_{\star}}
    &\bar{I}_{\mu}\P Y_{r,s} \ar[r]^{\mathrm{pr}_{r,s}} &\bar{I}_{\mu}Y\ ,}$$
  and $(\widetilde{EV_{\star}})_{*}$ is defined by
  $$\iota_{*}(r_{\star}(EV_{\star})_{*})\ .$$
  Note here $r_{\star}$ is the order of the band from the gerbe structure of
  $\bar{I}_{\mu}Y$ but not $\bar{I}_{\mu}\P Y_{r,s}$.
\end{enumerate}

First we state a vanishing lemma regarding applying localization formula to
\ref{eq:mainintegral1}.
\begin{lemma}\label{lem:vanish1}
  Assume $r,s$ is sufficiently large. If the localization graph $\Gamma$ has
  more than one vertex labeled by $\infty$, then the corresponding fixed loci
  moduli $F_{\Gamma}$ is empty, therefore it will contribute zero to
  \eqref{eq:mainintegral1}.
\end{lemma}
\begin{proof}
  For any twisted stable map $f:C\rightarrow \mathfrak \P Y_{r,s} $ in $\mathcal
  K_{0,\overrightarrow{m}\cup \star}(\P Y_{r,s}, (\beta,\frac{\delta}{r}))$,
  denote $N:=f^{*}(\mathcal O(-\mathcal D_{\infty}))$, first we show that
  $H^{1}(C, N )=0$. Indeed, using orbifold Riemann-Roch, we have
    $$\chi(N)=
    1+deg(N)-age(N |_{q_{\star}})=0 \ ,$$ as $deg(N)=-\frac{\delta}{r}$, and
    $age(N|_{q_{\star}})= 1-\frac{\delta}{r}$, then showing $H^{1}(C,N)=0$ is
    equivalent to show $H^{0}(C,N)=0$. As
    the degree of $N$ is negative on $C$, it remains to show that the degree of
    the restriction of the line bundle $N$ to every irreducible component $E$ of
    $C$ is non-positive. Observe that the degree $deg (N|_{E})$ is equal to
    intersection number $-([E], [\mathcal D_{\infty}])$. If the image of an
    irreducible component of $C$ via $f$ isn't contained in $\mathcal
    D_{\infty}$, the restricted degree of $N$ to $E$ is obviously non-positive.
    Otherwise, observe that $N$ is isomorphic to $(L^{\frac{1}{r}})^{\vee}$ over
    $\mathcal D_{\infty}$. As $L$ is semi-positive, the claim follows. This
    finishes the proof of the part $H^{1}(C,N )=0$.

    Now assume by contradiction that the moduli of fixed-loci $F_{\Gamma}$ is
    nonempty, by the connectedness of the graph $\Gamma$, there is at least one
    vertex of the graph $\Gamma$ labeled by $0$ with at least two edges
    attached. Suppose $f:C\rightarrow \P Y_{r,s}$ belongs to the moduli
    $F_{\Gamma}$ of $\C^{*}-$fixed loci. Assume that $C_{0}\cap C_{1}\cap C_{2}$
    is part of curve $C$, where $C_{0}$ is mapped by $f$ to $\mathcal D_{0}$ and
    $C_{1},C_{2}$ are edges meeting with $C_{0}$ at $b_{1}$ and $b_{2}$. Then in
    the normalization sequence for $R^{\bullet}\pi_{*}N$, it contains the part
    \begin{align*}\label{cd:vanishing}
      &H^{0}\left(C_{0}, N\right) \oplus H^{0}\left(C_{1}, N\right) \oplus H^{0}\left(C_{2}, N\right)\\
      \rightarrow &H^{0}\left(b_{1}, N\right) \oplus H^{0}\left(b_{2},
                    N\right)\\
      \rightarrow &H^{1}\left(C, N\right).
    \end{align*}
    Hence there is one of the weight-0 pieces in $H^{0}\left(b_{1}, N\right)
    \oplus H^{0}\left(b_{2}, N\right)$ that is canceled with a weight-0 piece of
    $H^{0}\left(C_{0}, N\right)$, and the other is mapped injectively into
    $H^{1}\left(C, N\right)$, but this contradicts that $H^{1}(C,N)=0$. So
    $F_{\Gamma}$ is empty.
  \end{proof}

  Now we prove Proposition \ref{prop:chara-general}.

\begin{proof}[Proof of Proposition \ref{prop:chara-general}]
  The proof is almost identical to the proof in~\cite[\S
  6.2]{wang19_mirror_theor_gromov_witten_theor_without_convex}, we here sketch
  the main steps for the convenience of readers. First without loss of
  generality, we will

  By Lemma \ref{lem:vanish1}, we
  only need to consider the decorated graph $\Gamma$ that has only one vertex
  labeled by $\infty$. Note the marking $q_{\star}$ corresponding to $\star$ is incident to the
  vertex $v_{\star}$ due to the choice of multiplicity at the marking
  $q_{\star}$, then the vertex $v_{\star}$ can't be a node linking two edges.
  Such decorated graph is \emph{star-shaped}, i.e., the vertex set $V$ allows a
  decomposition
  $$V= \{v_{\star}\}\bigsqcup V_{0}\ .$$
  where the vertex $v_{\star}$ is labeled by $\infty$, and all the vertexes in
  $V_{0}$ are labeled by $0$. Each vertex labeled by $0$ is linked to
  $v_{\star}$ by a unique edge. There are three types of star-shaped decorated
  graphs:
  \begin{enumerate}
  \item(Type I) The vertex $v_{\star}$ is unstable and there is only one edge
    incident to $v_{\star}$ in $\Gamma$, the unique vertex $v$ labeled by $0$ is
    unstable;
  \item(Type II) The vertex $v_{\star}$ is unstable and there is only one edge
    incident to $v_{\star}$ in $\Gamma$, the unique vertex $v$ labeled by $0$ is
    stable;
  \item(Type III) The vertex $v_{\star}$ is stable.
  \end{enumerate}

  Denote by $\beta_{\star}$ the degree of the unique vertex $v_{\star}$ labeled
  by $\infty$. Now let's compute the localization contribution from the above
  three types of graphs:
  \begin{enumerate}\label{case:graph2}
  \item If the graph $\Gamma$ is of type \textrm{I}, the vertex over $0$
    corresponds to a marked point with input $\mu_{\beta,\overrightarrow{k},c}$,
    then the graph $\Gamma$ contributes
    $$\mathbf{t}^{\overrightarrow{k}}\frac{\mu_{\beta,\vec{k},c}(\frac{\lambda-c_1(L)}{\delta})}{\delta}\cdot(\frac{\lambda-c_1(L)}{\delta})^{b}$$
    to \eqref{eq:mainintegralvar}. Note the use the fact that the restriction of
    the psi-class $\bar{\psi}_{\star}$ to $\mathcal M_{e}$ is equal to
    $\frac{\lambda-c_1(L)}{\delta}$(see Remark \ref{rmk:psi-edge}).
  \item If the graph $\Gamma$ is of type \textrm{II}, $\Gamma$ has the same
    description with type I except that the vertex over $0$ is stable. Then this
    type of graphs contributes
    \begin{equation*}
      \begin{split}
        &\sum_{\substack{(\beta_{j},\overrightarrow{k_j},c_{j})_{j=1}^{m}\in
            \mathrm{Adm}^{m}, \beta_{\star}\in \mathrm{Eff}(Y)\\
            \beta_{1}+\cdots+\beta_{m}+\beta_{\star}=\beta \\
            \overrightarrow{k_{1}}+\cdots+\overrightarrow{k_{m}}=\overrightarrow{k}
          }} \frac{1}{m!}(\widetilde{ev_{\star}})_{*}\bigg(\epsilon'_{*}\big(
        [\mathcal
        K_{0,\vec{m}_{s}\cup\star}(\sqrt[s]{L^{\vee}/Y},\beta_{\star})]^{\mathrm{vir}}\big)\\
        &\cap
        \bigcap_{i=1}^{m}ev_{i}^{*}(\mathbf{t}^{\overrightarrow{k_{i}}}\mu_{\beta_{i},\vec{k_i},c_{i}}(-\bar{\psi}_{i}))\cap
        \frac{\frac{1}{\delta}(\frac{\lambda-ev^*_\star
            c_{1}(L)}{\delta})^{b}}{\frac{\lambda-ev^*_\star
            c_{1}(L)}{s\delta}-\frac{\bar{\psi}_{\star}}{s}}\bigg)
      \end{split}
    \end{equation*} 
    to \eqref{eq:mainintegral1}, where
    $\vec{m}_{s}\cup\star=\big((c_{1}^{-1},e^{\frac{\delta_1}{s}},1),\cdots,
    (c^{-1}_{m},\mathrm{e}^{\frac{\delta_m}{s}},1),(c,e^{-\frac{\delta}{s}},1)\big)\in
    (C \!\times\! \C^{*} \!\times\! \C^{*})^{m+1}$,
    $\vec{m}\cup\star=\big(c^{-1}_{1},\cdots,c^{-1}_{m},c\big)\in C^{m+1}$
    and $$\epsilon':\mathcal
    K_{0,\vec{m}_{s}\cup\star}(\sqrt[s]{L^{\vee}/Y},\beta_{\star})\rightarrow
    \mathcal K_{0,\vec{m}\cup \star}(Y,\beta_{\star})$$ is the natural structure
    morphism by forgetting root of $\sqrt[s]{L^{\vee}/Y}$
    (c.f.\cite{tang2021quantum}). Note here we
    implicitly cancel out $a_{e}$ factors from the node contribution and
    automorphic factor everywhere and use the fact that the Euler class of
    virtual normal bundle for the vertex over $0$ is equal to 1. It's proved
    that
    $$ \epsilon'_{*}\big( [\mathcal
    K_{0,\vec{m}_{s}\cup\star}(\sqrt[s]{L^{\vee}/Y},\beta_{\star})]^{\mathrm{vir}}\big)=\frac{1}{s}
    [\mathcal K_{0,\vec{m}\cup\star}(Y,\beta_{\star})]^{\mathrm{vir}}\ , $$ in~
    \cite{tang2021quantum}, then the above formula
    is equal to
    \begin{equation*}
      \sum_{m\geq 0} \sum_{\substack{(\beta_{j},\overrightarrow{k_j},c_{j})_{j=1}^{m}\in \mathrm{Adm}^{m},
          \beta_{\star}\in \mathrm{Eff}(Y)\\ \beta_{1}+\cdots+\beta_{m}+\beta_{\star}=\beta \\ \overrightarrow{k_{1}}+\cdots+\overrightarrow{k_{m}}=\overrightarrow{k}}}
      \frac{1}{m!}\phi^{\alpha}\langle{\mu_{\beta_{1},\overrightarrow{k_{1}},c_{1}}(-\bar{\psi}_{1}),\cdots,\mu_{\beta_{m},\overrightarrow{k_{m}},c_{m}}(-\bar{\psi}_{m}),\frac{\frac{1}{\delta}(\frac{\lambda-c_1(L)}{\delta})^{b}\phi_{\alpha}}{\frac{\lambda-c_1(L)}{\delta}-\bar{\psi}_{\star}}}
      \rangle_{0,\vec{m}\cup \star,\beta_{\star}}\ .
    \end{equation*}
    
  \item If the graph $\Gamma$ is of type \textrm{III}, then $v_{\star}$ is
    incident to the distinguished leg corresponding to the marking $q_{\star}$
    and $m$ edges (m can be 0) and n legs. Each leg must be associated with the
    input $\mathbf{t}^{\overrightarrow{k'}}\mu_{\beta',\overrightarrow{k'},c'}(z)$ with
    $\beta'(L)+(\vec{w},\overrightarrow{k'})=0$, otherwise the support of the
    twisted sector $\bar{I}_{c'^{-1}}Y$ will avoid $\mathcal D_{\infty}$. Choose
    a labeling\footnote{Such a labeling is not unique, but we will divide $m!$
      to offset the labeling in the end.} of the $m$ edges attached to the
    vertex $v_{\star}$ by $[m]:=\{1,\cdots,m\}$. Let $\frac{\delta_{i}}{r}$ be
    the decorated degree associated with the $i$th edge $e_{i}$. Let $v_{i}$ be
    the vertex over $0$ incident to $e_{i}$, then $v_{i}$ can't be a unstable vertex
    of valence 1 as the corresponding ramification point $q_{0}$ is a stacky
    point, or it can't be a node linking two edges by Lemma \ref{lem:vanish1}.
    Therefore $v_{i}$ corresponds to either a marking point or a stable vertex.
    Assume that there are $l$ legs ($l$ can be zero) incident to $v$. let's label
    the legs incident to $v_{i}$ by $\{i1,\cdots,il\}\subset [\;|\{\textrm{Legs of
      $\Gamma$}\}|\;]$. Note that when $v_{i}$ is unstable, $l=1$.  

    Assume that the vertex $v_{i}$ is decorated by the degree $\beta_{i0}$.
    Assume that the insertion at the marking $q_{ij}$ on the curve\footnote{When $v$ is unstable, we just take $v$ to be $q_{i1}$.}
    $C_{v_{i}}$ corresponds to $\mathbf{t}^{\overrightarrow{k_{ij}}}\mu_{\beta_{ij},\overrightarrow{k_{ij}},c_{ij}}(-\bar{\psi}_{ij})$ in
    \eqref{eq:mainintegral1}, let's say the leg for $q_{ij}$ has \emph{virtual
      extended degree} $(\beta_{ij}, \overrightarrow{k_{ij}})$ contribution to
    the vertex $v_{i}$, denote $(\beta_{i},\overrightarrow{k_{i}})$ to be
    summation of $(\beta_{i0},0)$ and all the virtual degrees from the legs
    incident to $v_{i}$. We call $(\beta_{i},\overrightarrow{k_{i}})$ the \emph{total
      extended degree} at the vertex $v_{i}$. From \eqref{eq:mainintegral1}, one
    has
    $$\beta_{\star}+ \sum_{i=1}^{m}\beta_{i}=\beta,\quad
    \sum_{i=1}^{m}\overrightarrow{k_{i}}=\overrightarrow{k}\ . $$ Assume the
    multiplicity of the half edge incident to $C_{v_{i}}$ is equal to
    $(c^{-1}_{i},e^{\frac{\delta_{i}}{s}},1)$, here $\delta_{i}$ is the
    decorated degree associated with the edge $e_{i}$ and $\delta_{i}=
    age_{c}(L)\; \text{mod}\; \mathbb Z$ by \ref{prop:edge1}. Observe that to
    ensure such a graph $\Gamma$ exists, one must have
    \begin{equation}\label{eq:adm-graph}
      \beta_{i}(L)+ (\vec{w},\overrightarrow{k_{i}}) =\delta_{i}\ .
    \end{equation}
    Indeed, by
    orbifold Riemann-Roch Theorem, one has
    $$deg(f^{*}\mathcal O(\mathcal D_{0})|_{C_{v_{i}}})=-\frac{\beta_{i0}(L)}{s}=
    (1-\frac{\delta_{i}}{s})+\sum_{j=1}^{l}\frac{\beta_{ij}(L)+
    (\vec{w},\overrightarrow{k_{ij}}) }{s}\quad mod
    \quad \Z\ .$$ Here the first term on the right hand is the age of
    $f^{*}\mathcal O(D_{0})$ at the node of $C_{c_{i}}$, and the second term
    on
    the right is the sum of the ages of $f^{*}\mathcal O(\mathcal D_{0})$ at
    the
    marked points on $C_{c_{i}}$. As $s$ is sufficiently large, one must have
    $$\frac{\delta_{i}}{s}=\frac{\beta_{i0}(L)}{s}+\sum_{j=1}^{l}\frac{\beta_{ij}(L)+ (\vec{w},\overrightarrow{k_{ij}})
    }{s}\ ,$$
      which implies that
      $\beta_{i}(L)+(\vec{w},\overrightarrow{k_{i}})=\delta_{i}$.
    Note that this also imply that
    $(\beta_{i},\overrightarrow{k_i},c_{i})$ is an admissible pair.  

    Now we can group the edge-labeled decorated graphs by the set
    $\Lambda_{\beta,\overrightarrow{k},c,m,n}$ defined in \ref{def:star-graph}.
    For any element of $\Lambda_{\beta,\overrightarrow{k},c,m,n}$, we can
    naturally associate a group of edge-labeled star-shaped decorated graphs
    such that the vertex incident to the edge labeled by $i$ has total extended
    degree $(\beta_{i},\overrightarrow{k_{i}})$ and the multiplicity at the
    half-edge $h_{i}$ incident to $v_{i}$ over $0$ is $m(h_{i}):=
    (c^{-1}_{i},e^{\frac{\delta_{i}}{s}},1)$. We may call an element
    $\phi:=(c,\beta_{\star},((\beta_{1},\overrightarrow{k_{1}},c_{1}),\cdots,(\beta_{m+n},\overrightarrow{k_{m+n}},c_{m+n})))$
    of $\Lambda_{\beta,\overrightarrow{k},c,m,n}$ a \emph{meta graph}, and
    denote all the (star-shaped) edge-labeled decoration graph associated to $\phi$ by $\Gamma_{\phi}$.
        
    Now we use the localization formula in \S\ref{subsec:total1} to compute the
    contribution from $\Gamma_{\phi}$ to \eqref{eq:mainintegral1}. Summing over
    the localization contribution of the vertex $v_{i}$ together with branch of
    node $h_{i}$ at $v_{i}$ from all graphs in $\Gamma_{\phi}$, and pushing
    forward to $\bar{I}_{c_{i}^{-1}}Y\cong \bar{I}_{m_{h_{i}}}\mathcal D_{0}$
    along $\iota\circ (ev_{h_{i}})_{*}$(recall that $\iota$ is the inversion map
    on Chen-Ruan cohomology $H^{*}(\bar{I}_{\mu}Y)$), it yields
    \begin{equation*}
      \begin{split}
        &\mathbf{t}^{\overrightarrow{k_{i}}}\mu_{\beta_{i},\vec{k_i},c_{i}}\big(\frac{\lambda-c_1(L)}{\delta_{i}}\big)+
        \sum_{m=0}^{\infty}\sum_{\substack{
            (\beta_{j},\overrightarrow{k_j},c_{j})_{j=1}^{m}\in
            \mathrm{Adm}^{m},
            \beta_{\star}\in \mathrm{Eff}(Y)\\ \beta_{\star}+\sum_{j=1}^{m}\beta_{j}=\beta_{i}\\ \overrightarrow{k_{1}}+\cdots+\overrightarrow{k_{m}}=\overrightarrow{k}}}\frac{1}{m!}\\
        &(\widetilde{ev_{\star}})_{*}\big(\epsilon'_{*}\big( [\mathcal
        K_{0,\vec{m}_{s}\cup\star}(\sqrt[s]{L^{\vee}/Y},\beta_{\star})
        ]^{\mathrm{vir}}\big)\cap
        \bigcap_{j=1}^{m}ev_{j}^{*}(\mathbf{t}^{\overrightarrow{k_{j}}}\mu_{\beta_{j},\overrightarrow{k_{j}},c_{j}}(-\bar{\psi}_{j}))\cap
        \frac{1}{\frac{\lambda-ev^*_\star
            c_{1}(L)}{\delta_{j}s}-\frac{\bar{\psi}_{\star}}{s}}\big),
      \end{split}
    \end{equation*}
    which is equal to
    $\mathbf{t}^{\overrightarrow{k_{i}}}J^{Y}_{\beta_{i},\vec{k_i},c_{i}}(\mu,z)|_{\frac{\lambda-c_1(L)}{\delta_{i}}}$.
    Note here we use the fact the Euler class of the virtual normal bundle for
    $v_{i}$ is equal to $1$ as the node has nontrivial isotropy action on the
    normal bundle $N_{\mathcal D_{0}/\P Y_{r,s} }$. Observe that all the
    edge-labeled decorated graphs in $\Gamma_{\phi}$ have the same localization
    contributions from the unique vertex $v_{\star}$ labeled by $\infty$, the edge $e_{i}$ and the node over
    $\infty$ incident to $e_{i}$. Moreover the localization formula for any
    graph in $\Gamma_{\phi}$ depends multi-linearly on the localization
    contributions of vertexes over $0$, then the localization from all meta
    graphs in $\Lambda_{\beta,\overrightarrow{k},c,m,n}$ to
    \eqref{eq:mainintegral1} yields the summation:
    \begin{equation}
      \begin{split}
        &\sum_{m=0}^{\infty}\sum_{n=0}^{\infty} \sum_{\substack{\Gamma\in
            \Lambda_{\beta,\overrightarrow{k},c,m,n}\\ \Gamma\; \textrm{is
              stable}}}\frac{1}{m!n!}(\widetilde {ev_{\star}})_{*}\bigg(\sum
        _{d=0}^{\infty}\epsilon_{*}\big(c_{d}(-R^{\bullet}\pi_{*}f^{*}L^{\frac{1}{r}})(\frac{-\lambda}{r})^{-1+m-d}\\
        &\cap [\mathcal K_{0,\overrightarrow{m+n}\cup
          \star}(\sqrt[r]{L/Y},\beta_{\star})]^{\mathrm{vir}}\big) \cap
        \prod_{i=1}^{m}\frac{ev_{i}^{*}(\frac{1}{\delta_i} \mathbf{t}^{\overrightarrow{k_{i}}}J_{\beta_{i},\vec{k_i},c_{i}}(\mu,z)|_{z=\frac{\lambda-c_1(L)}{\delta_{i}}})}{-\frac{\lambda-ev^{*}_{i}c_1(L)}{r\delta_i}-\frac{\bar{\psi_{i}}}{r}}\\
        &\cap
        \prod_{i=m+1}^{m+n}ev_{i}^{*}(\mathbf{t}^{\overrightarrow{k_{i}}}\mu_{\beta_{i},\vec{k_i},c_{i}}(-\bar{\psi}_{i}))
        \cap \bar{\psi}_{\star}^{b}\bigg)\ .
      \end{split}
    \end{equation}

  \end{enumerate}
  By the discussion above, we can write \eqref{eq:mainintegral1} in the
  following way:
  \begin{equation}\label{eq:rec2*}
    \begin{split}
      &\mathbf{t}^{\overrightarrow{k}}\frac{\mu_{\beta,\overrightarrow{k},c}(\frac{\lambda-c_1(L)}{\delta})}{\delta}\cdot(\frac{\lambda-c_1(L)}{\delta})^{b}
      +\sum_{m=0}^{\infty}\sum_{\substack{(\beta_{j},\overrightarrow{k_j},c_{j})_{j=1}^{m}\in
          \mathrm{Adm}^{m},
          \beta_{\star}\in \mathrm{Eff}(Y)\\ \beta_{1}+\cdots+\beta_{m}+\beta_{\star}=\beta \\ \overrightarrow{k_{1}}+\cdots+\overrightarrow{k_{m}}=\overrightarrow{k}}}\frac{1}{m!}(\widetilde{ev_{\star}})_{*}\\
      &\bigg([\mathcal
      K_{0,\vec{m}\cup\star}(Y,\beta_{\star})]^{\mathrm{vir}}\cap
      \bigcap_{i=1}^{m}ev_{i}^{*}(\mathbf{t}^{\overrightarrow{k_{i}}}\mu_{\beta_{i},\overrightarrow{k_{i}},c_{i}}(-\bar{\psi}_{i}))\bigcap
      \frac{\frac{1}{\delta}(\frac{\lambda-ev^*_\star c_{1}(L)}{\delta})^{b}}{\frac{\lambda-ev^*_\star c_{1}(L)}{\delta}-\bar{\psi}_{\star}}\bigg)\\
      &+\sum_{m=0}^{\infty}\sum_{n=0}^{\infty} \sum_{\substack{\Gamma\in
          \Lambda_{\beta,\overrightarrow{k},c,m,n}\\ \Gamma\; \textrm{is
            stable}}}\frac{1}{m!n!}(\widetilde {ev_{\star}})_{*}\bigg(\sum
      _{d=0}^{\infty}\epsilon_{*}\big(c_{d}(-R^{\bullet}\pi_{*}f^{*}L^{\frac{1}{r}})(\frac{-\lambda}{r})^{-1+m-d}\\
      &\cap [\mathcal K_{0,\overrightarrow{m+n}\cup
        \star}(\sqrt[r]{L/Y},\beta_{\star})]^{\mathrm{vir}}\big) \cap
      \prod_{i=1}^{m}\frac{ev_{i}^{*}(\frac{1}{\delta_i} \mathbf{t}^{\overrightarrow{k_{i}}}J_{\beta_{i},\vec{k_i},c_{i}}(\mu,z)|_{z=\frac{\lambda-c_1(L)}{\delta_{i}}})}{-\frac{\lambda-ev^{*}_{i}c_1(L)}{r\delta_i}-\frac{\bar{\psi_{i}}}{r}}\\
      &\cap
      \prod_{i=m+1}^{m+n}ev_{i}^{*}(\mathbf{t}^{\overrightarrow{k_{i}}}\mu_{\beta_{i},\vec{k_i},c_{i}}(-\bar{\psi}_{i}))
      \cap \bar{\psi}_{\star}^{b}\bigg) \ .
    \end{split}
  \end{equation}
  As \eqref{eq:mainintegral1} lies in $\mathbf{t}^{\overrightarrow{k}}H^{*}(\bar{I}_{\mu}Y, \mathbb
  Q)[\lambda]$, the coefficient of $\mathbf{t}^{\overrightarrow{k}}\lambda^{-1}$ term in \eqref{eq:rec2*} must
  vanish. Note that the coefficients before $\lambda^{-1}$ in the first two
  terms in \eqref{eq:rec2*} yields
  \begin{equation*}
    \sum_{m=0}^{\infty}\sum_{\substack{(\beta_{j},\overrightarrow{k_j},c_{j})_{j=1}^{m}\in \mathrm{Adm}^{m},
        \beta_{\star}\in \mathrm{Eff}(Y)\\ \beta_{1}+\cdots+\beta_{m}+\beta_{\star}=\beta \\ \overrightarrow{k_{1}}+\cdots+\overrightarrow{k_{m}}=\overrightarrow{k}}}\frac{1}{m!}\phi^{\alpha}\langle{\mu_{\beta_{1},\overrightarrow{k_{1}},c_{1}}(-\bar{\psi}_{1}),\cdots,\mu_{\beta_{m},\overrightarrow{k_{m}},c_{m}}(-\bar{\psi}_{m}),\phi_{\alpha}\bar{\psi}^{b}_{\star}}
    \rangle_{0,\vec{m}\cup \star,\beta_{\star}},
  \end{equation*} 
  which is the left hand side of equality in \eqref{eq:recvar}. Then we extract
  the coefficient of the $\mathbf{t}^{\overrightarrow{k}}\lambda^{-1}$ term in the third term in
  \eqref{eq:rec2*}, this yields the term on the right hand side of
  \eqref{eq:recvar} up to a minus sign. This completes the proof of
  \eqref{eq:recvar}.
\end{proof}

\subsection{Specializing Novikov degrees}\label{subsec:spe-deg} 
When $Y$ can be embedded into another smooth DM stack $X$ with projective coarse
moduli and $L$ is a restriction of a line bundle of $X$, which we still denote
to be $l$ by an abuse of notation. For any degree $\beta\in \mathrm{Eff}(X)$, we
will denote $\mathcal K_{0,\vec{m}}(Y,\beta)$ to be the disjoint
union\footnote{It's a finite disjoint union as $\mathcal K_{0,m}(X,\beta)$ is
  finite type over $\mathbb C$, hence Noetherian. }
$$\bigsqcup_{\substack{d\in \mathrm{Eff}(Y)\\ i_{*}(d)=\beta}} \mathcal
K_{0,\vec{m}}(Y,d)\ .$$ The same rule applies to the notation when we define
Gromov-Witten invariants: for any degree $\beta\in \mathrm{Eff}(X)$, we will
denote the Gromov-witten invariants
$$\langle {\cdots} \rangle^{Y}_{0,\vec{m},\beta} := \sum _{d\in
  \mathrm{Eff}(Y):i_{*}(d)=\beta}\langle{\cdots} \rangle^{Y}_{0,\vec{m},d}\ .$$
Now we replace the effective cone $\mathrm{Eff}(Y)$ by $\mathrm{Eff}(X)$ in the
definition of $\Lambda_{\beta,\overrightarrow{k},c,m,n}$ and the recursion
relation \ref{eq:rec1} and apply the rule of specialization of degrees, we can
also apply the same strategy of proving of proposition \ref{prop:chara-general}
to show the following variation of \ref{prop:chara-general}:

\begin{proposition}\label{prop:deg-spe}
  Let $$\mu(z)=\sum_{\beta,\vec{k},c}q^{\beta}\mathbf{t}^{\overrightarrow{k}}\mu_{\beta,\vec{k},c}(z)$$ be an admissible series in
  $H^{*}(\bar{I}_{\mu}Y,\mathbb
  C)[z][\![t_{1},\cdots,t_{N}]\!][\![\mathrm{Eff}(X)]\!]$. For any integer
  $b\geq 0$ and admissible pair $(\beta,\overrightarrow{k},c)$ with
  $\beta(L)+(\vec{w},\vec{k})>0$, we have the following recursive relation:
  \begin{equation}\label{eq:rec1-spe}
    \begin{split}
      & [J^{Y}_{\beta,\vec{k},c}(\mu,z)]_{z^{-b-1}}\\
      &=\bigg[\sum_{m=0}^{\infty}\sum_{n=0}^{\infty} \sum_{\substack{\Gamma\in
          \Lambda_{\beta,\overrightarrow{k},c,m,n}\\ \Gamma\; \textrm{is
            stable}}}\frac{1}{m!n!}(\widetilde {ev_{\star}})_{*}\bigg(\sum
      _{d=0}^{\infty}\epsilon_{*}\big(c_{d}(-R^{\bullet}\pi_{*}f^{*}L^{\frac{1}{r}})(\frac{\lambda}{r})^{-1+m-d}(-1)^{d}\\
      &\cap [\mathcal K_{0,\overrightarrow{m+n}\cup
        \star}(\sqrt[r]{L/Y},\beta_{\star})]^{\mathrm{vir}}\big) \cap
      \prod_{i=1}^{m}\frac{ev_{i}^{*}(\frac{1}{\delta_i} \mathbf{t}^{\overrightarrow{k_{i}}}J^{Y}_{\beta_{i},\vec{k_i},c_{i}}(\mu,z)|_{z=\frac{\lambda-c_1(L)}{\delta_{i}}})}{\frac{\lambda-ev^{*}_{i}c_1(L)}{r\delta_i}+\frac{\bar{\psi_{i}}}{r}}\\
      &\cap
      \prod_{i=m+1}^{m+n}ev_{i}^{*}(\mathbf{t}^{\overrightarrow{k_{i}}}\mu_{\beta_{i},\vec{k_i},c_{i}}(-\bar{\psi}_{i}))
      \cap \bar{\psi}_{\star}^{b}\bigg)\bigg]_{\mathbf{t}^{\overrightarrow{k}}\lambda^{-1}}\ .
    \end{split}
  \end{equation}

\end{proposition}

\section{Proof of the main theorem}\label{sec:proofmain}
In this section, we will first assume that $X$ is a smooth Deligne-Mumford stack
with projective moduli and $Y\subset X$ is a smooth \emph{hypersurface} such
that the line bundle $L:=\mathcal O_{X}(Y)$ is \emph{semi-positive}, i.e.,
$\beta(L)\geq 0$ for any degree $\beta\in \mathrm{Eff}(X)$. Later in \S
\ref{subsec:mainthm}, $Y$ will be assumed to be a complete intersection
associated to a direct sum of semi-positive line bundles.

\subsection{A root stack modification of the space of the deformation to the
  normal cone}\label{subsec:space2}

Let $P_{X}:=\mathbb P_{X}(\mathbb C\oplus \mathbb C)$ be the trivial $\P^{1}$
bundle over $X$. It has two section: one is the zero section $X_{0}:=\mathbb
P(0\oplus \mathbb C)\subset P_{X}$, and the other is the $\infty-$section
$X_{\infty}:=\mathbb P(\mathbb C\oplus 0)$. We will introduce the notation
$P_{Y}$ to the mean the divisor $\mathbb P_{Y}(\mathbb C\oplus \mathbb C)
\subset P_{X}$.


Let $\mathfrak Q$ be the blow-up\footnote{This space is known as (a
  compactification of) the space of the deformation to the normal cone,
  c.f.,~\cite[ch5]{MR732620}.} of $P_{X}$ along $ Y_{0}:= X_{0}\cap P_{Y}$. More
explicitly, $\mathfrak Q$ can be constructed as a hypersurface in the projective
bundle $\pi:\mathfrak F:=\mathbb P_{ P_{X}}(L^{\vee}\oplus \mathcal
O_{P_{X}}(-1) )\rightarrow P_{X}$ associated to the section $z_{1} \cdot
\pi^{-1}(s_{P_{Y}})-z_{2}\cdot \pi^{-1}(s_{X_{0}})$ of the line bundle $\mathcal
O_{\mathfrak F}(1)$ , where $s_{Y}$ and $s_{X_{0}}$ are the defining equations
of $P_{Y}$ and $X_{0}$ in $P_{X}$, and $z_{1},z_{2}$ are the tautological
sections of line bundles $\pi^{*}(L^{\vee})\otimes O_{\mathfrak F}(1)$ and
$\pi^{*}\mathcal O_{P_{X}}(-1)\otimes \mathcal O_{\mathfrak F}(1)$ respectively.
The space $\mathfrak Q$ has four special smooth divisors: 1, the strict
transformation of $P_{Y}$, which we will still denote it to be $P_{Y}$; 2, the
exceptional divisor, which we will denote to be $E$ and it's isomorphic to the
projective bundle $\mathbb P_{Y}( L^{\vee}\oplus \mathbb C)$; 3, the strict
transformation of the $0-$section $X_{0}$, we will still denote to be $X_{0}$ by
an abuse of notation, and it's isomorphic to $X$ with normal bundle equal to
$L$, 4, the strict transformation of the $\infty-$section $X_{\infty}$, which we
still denote it to be $X_{\infty}$.


Let $\C^{*}$ act on $P_{X}$ by scaling the $\P^{1}-$fiber so that the weight of
the $\C^{*}-$action on the normal bundle of $X_{\infty}$ in $P_{X}$ is $-1$. It
induces a $\C^{*}-$action on $\mathfrak Q$ with the $\C^{*}-$fixed loci $X_{0}$,
$X_{\infty}$ and $D_{0}:=P_{Y}\cap E$.

We will consider the root stack $\mathfrak R$ of $\mathfrak Q$ by taking $s-$th
root of $X_{0}$ and $r-$th root of $P_{Y}$. We will also assume that $r,s$ are
distinct primes. We will still use the notation $X_{0}$, $X_{\infty}$, $E$ and
$P_{Y}$ to mean their corresponding divisors in $\mathfrak R$ after taking
roots. We note $E$ is isomorphic to $\P Y_{r,s} $ so that $\mathcal
D_{0}:=X_{0}\cap E $ is isomorphic to the root stack $\sqrt[s]{L^{\vee}/Y}$ and
$\mathcal D_{\infty}:= P_{Y}\cap E $ is isomorphic to the root stack
$\sqrt[r]{L/Y}$.

The $\C^{*}-$action on $\mathfrak Q$ induces a $\C^{*}-$action on $\mathfrak R$.
We see that the normal bundle of $X_{0}$ in $\mathfrak R$ is
($\C^{*}-$equivariantly) isomorphic to $(L^{\vee})^{\frac{1}{s}}\otimes \mathbb
C_{\frac{1}{s}}$, the normal bundle of $X_{\infty}$ in $\mathfrak R$ is
isomorphic to $\mathbb C_{-\lambda}$, and the normal bundle of $\mathcal
D_{\infty}$ in $\mathfrak R$ is isomorphic to $(L^{\frac{1}{r}}\otimes \mathbb
C_{-\frac{1}{r}})\oplus \mathbb C_{\lambda}$.

Let $\mathrm{pr}_{r,s}: \mathfrak R\rightarrow X$ be the morphism induced from
composition of the following three maps: (1) a morphism $\mathfrak R\rightarrow
\mathfrak Q$ forgetting root structure; (2) the blow-down from $\mathfrak Q$ to
$P_{X}$; (3) and the projection from $P_{X} $ to the base $X$. Note
$\mathrm{pr_{r,s}}$ induces a morphism from the rigidified inertia stack
$\bar{I}_{\mu} \mathfrak R$ to the rigidified inertia stack $\bar{I}_{\mu}X$.

There are two morphisms
$$q_{0},q_{\infty}: \mathfrak R \rightarrow \mathbb B \C^{*}  \ ,$$
associated to the line bundle $\mathcal O( X_{0})$ and $\mathcal O(P_{Y})$
respectively such that $\mathcal O(X_{0})\cong q_{0}^{*}(\mathbb L)$ and
$\mathcal O(P_{Y})\cong q_{\infty}^{*}(\mathbb L)$, here $\mathbb L$ is the
universal line bundle over $\mathbb B \C^{*}$. This will induce morphisms on
their rigidified inertia stack counterparts:
$$ q_{0},q_{\infty}: \bar{I}_{\mu}\mathfrak R \rightarrow \bar{I}_{\mu} \mathbb B \C^{*}  \ ,$$
Note the rigidified inertia stack $ \bar{I}_{\mu} \mathbb B \C^{*}$ of $\mathbb
B \C^{*}$ can be written as the disjoint union $$ \bigsqcup_{c\in \C^{*}}
\bar{I}_{c} \mathbb B \C^{*} \ ,$$ Let $c$ be a complex number, any $\mathbb
C$-point of $ \bar{I}_{c} \mathbb B \C^{*}$ is isomorphic to the pair
$(1_{\mathbb C},c)$, where $1_{\mathbb C}$ is the trivial principal $\mathbb
C$-bundle over $\mathbb C$ and $c$ is an element in the automorphism group
$\mathrm{Aut}_{\mathbb C}(1_{\mathbb C})\cong \C^{*}$. Let $(y,g)$ be a $\mathbb
C-$point of $\bar{I}_{\mu}\mathfrak R$, where $y\in Ob(\P Y_{r,s}(\mathfrak R)
)$ and $g\in \mathrm{Aut}(y)$. Assume that $q_{0}((y,g))=(1_{\mathbb C}, c)$,
then $g$ acts on the fiber $\mathcal O(X_{0})|_{y}$ via the multiplication by
$c$. We have a similar relation for $q_{\infty}$.

Let $C$ be a good index set for rigidified inertia stack $\bar{I}_{\mu}X $
satisfying the assumption \ref{assump:inertiaindex}. Now we will use the index
set $C\times \C^{*} \times \C^{*}$ to index the components of the rigidified
inertia stack $\bar{I}_{\mu}\mathfrak R$: for each $c\in C$,
$c_{0},c_{\infty}\in \C^{*}$, we will define the rigidified inertia component $
\bar{I}_{(c,c_{0},c_{\infty})} \mathfrak R$ to be
$$
\mathrm {pr}_{r,s}^{-1}(\bar{I}_{c}X)\cap q_{0}^{-1} (\bar{I}_{c_{0}}\mathbb B
\C^{*} ) \cap q_{\infty}^{-1} (\bar{I}_{c_{\infty}}\mathbb B \C^{*} )\ .$$
Notice that $\bar{I}_{(c,c_{0},c_{\infty})}\mathfrak R$ is nonempty if and
only if $s\cdot age_{c}(L)\equiv -c_{0}$ mod $\mathbb Z$ and $r\cdot age_{c}(L)\equiv c_{\infty}$ mod $\mathbb Z$.

\begin{definition}\label{def:deg2}
  For any degree $\beta\in \mathrm{Eff}(X)$ and $\delta,d\in \mathbb Q$, we say
  a stable map $f:C\rightarrow \mathfrak R$ is of degree
  $(\beta,\frac{\delta}{r},d)$ if $(\mathrm{pr}_{r,s}\circ f)_{*}[C]=\beta$,
  $deg(f^{*}\mathcal O(P_{Y}))=\frac{\delta}{r}$ and $deg(f^{*}\mathcal
  O(X_{\infty}))=d$. We will denote $\mathcal K_{0,\overrightarrow{m}}(\mathfrak
  R,(\beta,\frac{\delta}{r},d))$ to be the corresponding moduli stack of twisted
  stable maps to $\mathfrak R$ of degree $(\beta,\frac{\delta}{r},d)$. We note
  $\mathcal K_{0,\overrightarrow{m}}(\mathfrak R,(\beta,\frac{\delta}{r},d))$ is
  proper as the degree data uniquely determines a homology class in
  $H_{2}(\mathfrak R,\mathbb Q)$. Indeed, dual to $H_{2}(\mathfrak R,\mathbb
  Q)$, we have $H_{2}(\mathfrak R,\mathbb Q)^{\vee}\cong H^{2}(\mathfrak
  R,\mathbb Q)\cong H^{2}(\mathfrak Q,\mathbb Q)$, and, by the \emph{blow-up}
  construction of $\mathfrak Q$, $H^{2}(\mathfrak Q,\mathbb Q)$ is isomorphic to
  the direct sum
$$H^{2}(X,\mathbb Q)  \oplus \mathbb Q[X_{\infty}] \oplus \mathbb Q[P_{Y}] \ ,$$
where $[P_{Y}]$, $[X_{\infty}]$ (and $[E]$) are the fundamental classes of
$P_{Y}$, $X_{\infty}$ (and $E$) respectively and the first summand is embedded into $H^{*}(\mathfrak Q,\mathbb Q)$
via the pullback $\mathrm{pr}_{r,s}^{*}$. Note that we have
$r[P_{Y}]+[E]=\mathrm{pr}_{r.s}^{*}([Y])$ in $H^{2}(\mathfrak R,\mathbb Q)$.
Note that the line bundles $\mathcal O(X_{\infty})$ and $\mathcal O(P_{Y})$ are semi-positive. 
\end{definition}

\subsection{Localization analysis}\label{sec:lc2}
The $\C^{*}-$action on $\mathfrak R$ induces a $\C^{*}-$action on the moduli of
stable maps to $\mathfrak R$, we will use decorated graphs (trees) to index the
components of $\C^{*}-$fixed loci of $ \mathcal
K_{0,\overrightarrow{m}}(\mathfrak R,(\beta,\frac{\delta}{r},d))$ similar to
\ref{sec:lc1}. The decorations on a graph $\Gamma$ is given by the following data:
\begin{itemize}
\item Each vertex $v$ is associated with an index $j(v) \in \{X_{0}, \mathcal
  D_{\infty}, X_{\infty}\}$, and a degree $\beta(v) \in \mathrm{Eff}(X)$. In
  particular, we will also say the vertex $v$ is labeled by $0$ if $j(v)=X_{0}$
  and is labeled by $\infty$ if $j(v)=\mathcal D_{\infty}$.
\item Each edge $e$ consists a pair of half-edges $\{h_{j},h_{j'}\}$, and $e$ is
  equipped with two numbers $\delta(e),d(e) \in \mathbb{Q}$. Here we call
  $h_{j}$ and $h_{j'}$ half edges and the subscript $j$ or $j'$ is taken in the
  set $ \{X_{0}, \mathcal D_{\infty}, X_{\infty}\} $. A half-edge $h_{j}$
  labeled by $j$ is incident to a vertex labeled by $j$.
\item Each half-edge $h$ and each leg $l$ has an element $m(h)$ or $m(l)$ in $C
  \!\times\! \C^{*} \!\times\! \C^{*} $.
\item The legs are labeled with the numbers $\{1, \ldots, m\}$ and each leg is
  incident to a unique vertex.
\end{itemize}

For $j\in \{X_{0}, \mathcal D_{\infty}, X_{\infty}\}$, we will use the symbol
$\mathfrak R_{j}$ to mean the space $j$. By the ``valence" of a vertex $v$,
denoted $\text{val}(v)$, we mean the total number of incident half-edges and
legs. For each $\C^{*}-$fixed stable map $f: (C; q_1, \ldots, q_m)\rightarrow \mathfrak R $ in $\mathcal
K_{0,\overrightarrow{m}}(\mathfrak R,(\beta_{\star},\frac{\delta}{r},d))$, we
can associate a decorated graph $\Gamma$ in the following way.
\begin{itemize}
\item Each edge $e$ corresponds to a genus-zero component $C_e$ of restricted
  stale-map degree $(0,\frac{\delta(e)}{r},d(e))$. Thus $f$ maps
  $C_{e}$ constantly to the base. There are two distinguished points(we will
  also call ramification points) $q_{j}$ and $q_{j'}$ on $C_{e}$ satisfying that
  $q_{j}$ maps to $\mathfrak R_{j}$ and $q_{j'}$ maps to $\mathcal R_{j'}$,
  respectively. There are two corresponding half-edges $h_{j}$ and $h_{j'}$
  associated to $q_{j}$ and $q_{j'}$ respectively.
\item Each vertex $v$ for which $j(v)$ (with unstable exceptional cases noted
  below) corresponds to a maximal sub-curve $C_v$ of $C$ which maps totally into
  $\mathcal X_{0}$, then the restriction of $f$ to $C_{v}$ defines a twisted stable map in
  $$\mathcal
  K_{0,val(v)}(\mathfrak R_{j(v)},\beta(v))\ .$$ The label $\beta(v)$ denotes
  the degree coming from the composition $\mathrm{pr}_{r,s}\circ f|_{C_{v}}$ of
  the restriction and projection to the base.
\item A vertex $v$ is {\it unstable} if stable twisted maps of the type
  described above do not exist. In this case, $v$ corresponds to a single point
  of the component $C_e$ for each incident edge $e$, which may be a node at
  which $C_e$ meets another edge curve $C_{e'}$, a marked point of $C_e$, or an
  unmarked point.
\item The index $m(l)$ on a leg $l$ indicates the rigidified inertia stack
  component $\bar{I}_{m(l)}\mathfrak R$ of $\mathfrak R$ on which the gerby
  marked point corresponding to the leg $l$ is evaluated.
\item A half-edge $h$ of an edge $e$ corresponds a ramification point $q \in
  C_{e}$ such that $f(q)\in \mathfrak R_{j}$ if $h$ is labeled by $j$. Then
  $m(h)$ indicates the rigidified inertia component $\bar{I}_{m(h)}\mathfrak R$
  of $\mathfrak R$ on which the ramification point $q$ associated with $h$ is
  evaluated.
\end{itemize}

Moreover, the decorated graph from a fixed stable map should satisfy the
following:
\begin{enumerate}\label{prop:edge2}
\item(Monodromy constraint for edge) For each edge
  $e=\{h_{0}:=h_{x_{0}},h_{\infty}:=h_{\mathcal D_{\infty}}\}$ linking $X_{0}$
  and $\mathcal D_{\infty}$, we have that $m(h_{0})=(c^{-1},
  \mathrm{e}^{\frac{\delta(e)}{s}},1)$ and $m(h_{\infty})=(c, 1,
  \mathrm{e}^{\frac{\delta(e)}{r}})$ for some $c\in C$. Moreover the choice of
  $c$ should satisfy that $$\delta(e)= age_{c}(L)\; \text{mod}\; \mathbb Z\ ,$$
  which follows from the fact that we have an isomorphism of line bundles
  $f|^{*}_{C_{e}}\mathcal O(sX_{0})\cong f|^{*}_{C_{e}}(\mathcal
  O(rP_{Y})\otimes \mathrm{pr}_{r,s}^{*}L^{\vee})$ on $C_{e}$ and then the ages
  of both line bundles at $q_{\infty}$ should be equal.
\item(Monodromy constraint for vertex) For each vertex $v$, let $I_{v}\subset
  [m]$ be the set of incident legs, and $H_{v}$ be the set of incident
  half-edges. For each $i\in I_{v}$, write $m(l_{i})=(c_{i},a_{i},b_{i})$, and
  for each $h\in H_{v}$, write $m(h)=(c_{h},a_{h},b_{h})$. If $v$ is labeled by
  $0$, we have that $b_{i}$, $b_{h}$ are equal to $1$, and
$$ \mathrm{e}^{\frac{\beta(v)(L)}{s}}\times \prod_{i\in I_{v}} a_{i} \times \prod_{h\in H_{v}}
a_{h}^{-1}=1\ .$$ If $v$ is labeled by $\infty$, then we have $a_{i}$ and
$a_{h}$ are all equal to $1$ and
 $$ \mathrm{e}^{-\frac{\beta(v)(L)}{r}}\times \prod_{i\in I_{v}} b_{i} \times \prod_{h\in E_{v}}b^{-1}_{h}=1\ .$$  
\end{enumerate}

In particular, we note that the decorations at each stable vertex $v$ yield a
tuple
$$\overrightarrow{val}(v) \in (C \times \C^{*}\times \C^{*} )^{\text{val}(v)}$$
recording the multiplicities at every special point of $C_v$. Then the
restriction gives a stable map in
$$ \mathcal
K_{0,\overrightarrow{val}(v)}(\mathfrak R_{j(v)},\beta(v))\ .$$

Assume that $r,s$ are sufficiently large primes. We will do a similar localization analysis as in \S \ref{sec:lc1}.

\subsubsection{Vertex contribution}\label{subsubsec:ver-cntr2}
Assume that $v$ is a stable vertex, the localization contributions for $X_{0}$
and $\mathcal D_{\infty}$ repeat the discussions for $\mathcal D_{0}$ and
$\mathcal D_{\infty}$ in \ref{subsubsec:ver-cntr1} with the replacement of the
words $\mathcal D_{0}$ and $\mathcal D_{\infty}$ by $X_{0}$ and $\mathcal
D_{\infty}$ respectively, and one additional change that the \emph{inverse of
  the Euler class} of the virtual normal bundle for $\mathcal D_{\infty}$ should
be equal to
$$\frac{1}{\lambda}\cdot \sum_{d\geq 0}c_{d}(-R^{\bullet}\pi_{*}f^{*}(L)^{\frac{1}{r}})(\frac{-\lambda}{r})^{|E(v)|-1-d} \ .$$

If $v$ is a stable vertex labeled by $X_{\infty}$, then the vertex moduli $\mathcal
M_{v}$ is given by $\mathcal K_{0,\overrightarrow{val}(v)}(X_{\infty},\beta(v))$
and the fixed part of perfect obstruction theory gives rise to $[\mathcal
K_{0,\overrightarrow{val}(v)}(X_{\infty},\beta(v))]^{\mathrm{vir}}$. The movable
part of the perfect obstruction theory yields the \emph{inverse of the Euler
  class} of the virtual normal bundle which is equal to $\frac{-1}{\lambda}$.

When $v$ is an unstable vertex
over $0$ (resp. $\infty$), let $h$ be the half-edge incident to $v$, the vertex
moduli $\mathcal M_{v}$ is defined to be $\bar{I}_{m(h)^{-1}}\mathcal
D_{0}$(resp. $\bar{I}_{m(h)^{-1}}\mathcal D_{\infty}$) with $[\mathcal
M_{v}]^{\mathrm{vir}}=[\mathcal M_{v}]$ and Euler class of virtual normal
bundle equal to $1$ (resp. $\frac{1}{\lambda}$).

\subsubsection{Edge contribution} \label{subsubsec:edge-cntr2} Here we will only
consider the type of edges linking the vertex labeled by $X_{0}$ and the vertex
labeled by $\mathcal D_{\infty}$ as it is only the type of edge which appears in
the localizaiton contribution of \ref{eq:mainintegralvar}, see Lemma
\ref{lem:edge-exlusion} for more details.

Let $e=\{h_{0}:=h_{X_{0}},h_{\infty}:=h_{\mathcal D_{\infty}}\}$ be an edge in
$\Gamma$ with decorated degree $\delta(e)\in \mathbb Q$ and decorated
multiplicities $m(h_{0})$ and $m(h_{\infty})$. Write $m(h_{\infty})=(c, 1,
\mathrm{e}^{\frac{\delta(e)}{r}})$, let $a_{e}:=a(c)$ be the number as
Assumption \ref{assump:inertiaindex}, we will also write $a$ as $a_{e}$ for
simplicity if there is no confusion. Then we define the edge moduli $\mathcal
M_{e}$ to be the root stack $ \sqrt[as\delta(e)]{ L^{\vee}/I_{c}Y }$. Then the
virtual cycle $[\mathcal M_{e}]^{vir}$ coming from the fix part of the
obstruction theory is the fundamental class $[\mathcal M_{e}]$ and \emph{the
  inverse of the Euler class} of the virtual normal bundle is $\prod_{i=1}^{-
  1-\lfloor {-\delta(e)} \rfloor}\big(\lambda+\frac{i}{\delta(e)}(
c_{1}(L)-\lambda ) \big)$.\footnote{Here, for any rational number $q$, we denote
  $\lfloor q \rfloor$ to be the largest integer no larger than $q$.} In
practice, we will use an equivalent form
$$\frac{1}{\lambda}\prod_{0\leq
  j<\delta(e)}\big(c_{1}(L)+(\delta(e)-j)\frac{\lambda-c_1(L)}{\delta(e)}\big)\
.$$
We note that $\mathcal M_{e}$ allows a finite \'etale map into the corresponding fixed-loci in $\mathcal
K_{0,2}(\mathfrak R, (0,\frac{\delta(e)}{r},0)$ of degree $\frac{1}{as}$. See
appendix \ref{sec:app-edge} for more details.

\subsubsection{Node contributions}\label{subsubsec:node-cntr2}
The deformations in $\mathcal K_{0,\overrightarrow{m}}(\mathfrak R,
(\beta,\frac{\delta}{r},0))$ smoothing a node contribute to the Euler class of
the virtual normal bundle as the first Chern class of the tensor product of the
two cotangent line bundles at the branches of the node. For nodes at which a
component $C_e$ meets a component $C_v$ over the vertex $X_{0}$, this
contribution is
\begin{equation*}
  \frac{\lambda - c_{1}(L)}{a_{e}s\delta(e)} - \frac{\bar{\psi}_v}{a_{e}s}.
\end{equation*}

For nodes at which a component $C_e$ meets a component $C_v$ at the vertex over
$\mathcal D_{\infty}$, this contribution is
\begin{equation*}
  \frac{-\lambda + c_{1}(L)}{a_{e}r\delta(e)} -\frac{\bar{ \psi}_v}{a_{e}r}\ .
\end{equation*}
There is also one dimensional piece from deformations of maps (corresponding to
the normal direction of $E$ in $\mathfrak R$) at the node over
$\mathcal D_{\infty}$, which contributes $\frac{1}{\lambda}$ as the Euler
class of the virtual normal bundle to the localization.

We will not need node contributions from other types of nodes as the above types
suffice the need in this paper, see Lemma \ref{lem:edge-exlusion} for more
details.

\subsubsection{Total contribution}
Here we will only consider two types of graph: (1) all edges in $\Gamma$
connects $X_{0}$ and $\mathcal D_{\infty}$; (2) $\Gamma$ has only one vertex
labeled by $X_{\infty}$ and no edges. They are all we need in the later
localization analysis(see Lemma \ref{lem:edge-exlusion} for more details).

The localization contribution from the decorated graph of type (2) is
$$ -\frac{[\mathcal M_{v}]^{\mathrm{vir}}}{\lambda}\ .$$

For any decorated graph $\Gamma$ of type (1), we define $F_{\Gamma}$ to be
$$\prod_{v:j(v)=0}\mathcal M_{v}\times_{\bar{I}_{\mu}X_{0}}\prod_{e\in E}\mathcal
M_{e}\times_{\bar{I}_{\mu}\mathcal D_{\infty}} \prod_{v:j(v)=\infty} \mathcal
M_{v}
$$ of the following
diagram:
$$\xymatrix{
  F_{\Gamma}\ar[r]\ar[d] &\prod\limits_{v:j(v)=0} \mathcal M_{v}\times
  \prod\limits_{e\in E}\mathcal M_{e}\times \prod\limits_{v:j(v)=\infty}
  \mathcal M_{v}
  \ar[d]^-{ev_{nodes}}\\
  \prod\limits_{E}(\bar{I}_{\mu}X_{0} \times \bar{I}_{\mu}\mathcal
  D_{\infty})\ar[r]^-{(\Delta^{0}\times \Delta^{\infty})^{|E|}}
  &\prod\limits_{E}(\bar{I}_{\mu} X_{0})^{2}\times (\bar{I}_{\mu}\mathcal
  D_{\infty})^{2}\ , }
$$
where $\Delta^{0}=(id,\iota)$(resp. $\Delta^{\infty}=(id,\iota)$) is the
diagonal map of $\bar{I}_{\mu}X_{0}$ (resp. $\bar{I}_{\mu}\mathcal D_{\infty}$).
Here the right-hand vertical map is the product of
evaluation maps at the two branches of each gluing node.

We define $[F_{\Gamma}]^{\mathrm{vir}}$ to be:
\begin{equation*}
  \begin{split}
    \prod_{v:j(v)=0}[\mathcal M_{v}]^{\mathrm{vir}}
    \times_{\bar{I}_{\mu}X_{0}}\prod_{e\in E}[\mathcal
    M_{e}]^{\mathrm{vir}}\times_{\bar{I}_{\mu}\mathcal D_{\infty}}
    \prod_{v:j(v)=\infty}[\mathcal M_{v}]^{\mathrm{vir}}\ .
  \end{split}
\end{equation*}
Then the contribution of decorated graph $\Gamma$ to the virtual localization
is:
\begin{equation}\label{eq:auto-loc2}
  Cont_{\Gamma}=\frac{\prod_{e\in
      E}sa_{e}}{|\text{Aut}(\Gamma)|}(\iota_{\Gamma})_{*}\left(\frac{[F_{\Gamma}]^{\mathrm{vir}}}{e^{\C^{*}}(N^{\mathrm{vir}}_{\Gamma})}\right)\ .
\end{equation}
Here $\iota_{F}:F_{\Gamma}\rightarrow \mathcal K_{0,\overrightarrow{m}}(
\mathfrak R ,(\beta,\frac{\delta}{r},d))$ is a finite \'etale map of degree
$\frac{|\text{Aut}(\Gamma)|}{\prod_{e\in E}sa_{e}}$ into the corresponding
$\C^{*}$-fixed loci. The virtual normal bundle
$e^{\C^{*}}(N^{\mathrm{vir}}_{\Gamma})$ is the product of virtual normal bundles
from vertex contributions, edge contributions and node contributions.
\subsection{Auxiliary cycle}
Take $L_{1}=L$ and $r=1$ in \ref{def:adm-series} and write $\vec{w}$ to be the
weight. Let
$\mu^{X}(z)=\sum_{\beta,\vec{k},c}q^{\beta}\mathbf{t}^{\overrightarrow{k}}\mu^{X}_{\beta,\vec{k},c}(z)$ be an admissible series in
$H^{*}_{CR}(X,\mathbb C)[z][\![t_{1},\cdots,t_{N}]\!][\![\mathrm{Eff}(X)]\!]$.
For any triple $(\beta,\vec{k},c)\in \mathrm{Eff}(X)\times \mathbb Z_{\geq
  0}^{N}\times C$, we define $J^{X,tw}_{\beta,\vec{k},c}(\mu^{X},z)$ to be
\begin{equation}\label{eq:twcomp}
  \begin{split}
    &\bigg(\mu^X_{\beta,\vec{k},c}(z)+\mathbf{Coeff}_{\mathbf{t}^{\overrightarrow{k}}}\bigg[\sum_{m\geq
      0} \sum_{\substack{(\beta_{j}, \overrightarrow{k_{j}},c_{j})_{j=1}^{m}\in \mathrm{Adm}^{m}\; \beta_{\star}\in \mathrm{Eff}(X)\\ \beta_{1}+\cdots+\beta_{m}+\beta_{\star}=\beta\\
      }} \frac{1}{m!}\phi^{\alpha}\langle
    \mathbf{\mathbf{t}^{\overrightarrow{k_{1}}}\mu^{X}_{\beta_{1},\overrightarrow{k_{1}},c_{1}}}(-\bar{\psi}_{1}),\\
    &\cdots,\mathbf{\mathbf{t}^{\overrightarrow{k_{m}}}\mu^{X}_{\beta_{m},\overrightarrow{k_{m}},c_{m}}}(-\bar{\psi}_{m}),\frac{\phi_{\alpha}}{z-\bar{\psi}_{\star}}\rangle^{X}_{0,\vec{m}\cup\star,\beta_{\star}}\bigg]\bigg)
    \prod_{0\leq j< \beta(L)+(\vec{w},\vec{k})}\big(c_{1}(L) +(\beta(L)+(\vec{w},\vec{k}) -j)z
    \big)\ .
  \end{split}
\end{equation}
where $\overrightarrow{m}\cup\star =\big(c_{1}^{-1}\cdots,c_{m}^{-1},c\big)\in
C^{m+1}$. We note that $J^{X,tw}_{0,\vec{0},c}=0$, then the formal series
$$J^{X,tw}(q,\mu^{X},z)=z+\sum_{\beta,\vec{k},c}q^{\beta}\mathbf{t}^{\overrightarrow{k}}J_{\beta,\vec{k},c}^{X,tw}(\mu^{X},z)$$ associated with the input $\mu^{X}(z)$ above is
an admissible series and near $z$.

For any admissible pair $(\beta,\vec{k},c)$ with $\beta(L)+(\vec{w},\vec{k})>0$,
denote $\delta=\beta(L)+ (\vec{w},\vec{k})$. Assume that $r,s$ are sufficiently
large primes. For any nonnegative integer $b$, we will also compare
\eqref{eq:mainintegral1} to the following auxiliary cycle:
\begin{equation}\label{eq:mainintegralvar}
  \begin{split}
    &\sum_{m=0}^{\infty}\sum_{\substack{
        (\beta_{i},\overrightarrow{k_{i}},c_{i})\in \mathrm{Adm}^{m}, \beta_{\star}\in \mathrm{Eff}(X)\\
        \beta_{\star}+\sum_{i=1}^{m}\beta_{i}=\beta \\
        \overrightarrow{k_{1}}+\cdots+\overrightarrow{k_{m}}=\overrightarrow{k}
      }}\frac{1}{m!}(\widetilde{EV_{\star}})_{*}\bigg(
    \prod_{i=1}^{m}ev_{i}^{*}\big(\mathrm{pr}_{r,s}^{*}(\mathbf{t}^{\overrightarrow{k_{i}}}\mu^{X}_{\beta_{i},\vec{k_i},c_{i}}(-\bar{\psi}_{i}))\big)\\
    &\cap \bar{\psi}_{\star}^{b}\cap [\mathcal K_{0,\overrightarrow{m}\cup
      \star}( \mathfrak R
    ,(\beta_{\star},\frac{\delta}{r},0))]^{\mathrm{vir}}\bigg)
  \end{split}
\end{equation}
Here an explanation of the notations is in order:
\begin{enumerate}
\item For degrees $\beta_{\star},\beta_{1},\cdots, \beta_{m}$ in
  $\mathrm{Eff}(X)$ and tuples
  $\overrightarrow{k_{1}},\cdots,\overrightarrow{k_{m}}$ in $\mathbb Z_{\geq
    0}^{N}$ with $\sum_{i=1}^{m}\beta_{i}+\beta_{\star}=\beta$ and $
  \overrightarrow{k_{1}}+\cdots+\overrightarrow{k_{m}}=\overrightarrow{k}$, let
  $(c_{1},\cdots,c_{m})\in C^{m}$ such that
  $(\beta_{i},\overrightarrow{k_i},c_{i})$ are admissible pairs. Write
  $\delta_{i}= \beta_{i}(L)+(\vec{w},\overrightarrow{k_{i}})$, we define
  $\vec{m}\cup \star$ to be the $m+1$tuple
$$\big((c^{-1}_{1},\mathrm{e}^{\frac{\delta_i}{s}},1),\cdots,(c^{-1}_{m},\mathrm{e}^{\frac{\delta_m}{s}},1),(c,1,\mathrm{e}^{\frac{\delta}{r}})\big)\ . $$ 

\item Note $\bar{I}_{m_{\star}}\mathfrak R$ is canonically isomorphic to
  $\bar{I}_{c}Y\times \P^{1}$, it has a natural projection map $p:
  \bar{I}_{m_{\star}}\mathfrak R\rightarrow \bar{I}_{c}Y $. Let $EV_{\star}$ be
  the morphism which is a composition of the following maps:
  $$\xymatrix{ \mathcal K_{0,\overrightarrow{m}\cup \star}(\mathfrak R,(\beta_{\star},\frac{\delta}{r},0))\ar[r]^-{ev_{\star}}
    &\bar{I}_{m_{\star}}\mathfrak R \ar[r]^{p} &\bar{I}_{c}Y\ .}$$ Then we
  define $(\widetilde{EV_{\star}})_{*}$ to be
  $$\iota_{*}(r_{\star}(EV_{\star})_{*})$$
  as in \ref{tilde-ev}. Note here $r_{\star}$ is the order of the band from the
  gerbe structure of $\bar{I}_{\mu}Y$ but not $\bar{I}_{\mu}\mathfrak R$.
\end{enumerate}

Applying $\C^{*}-$localization to the \ref{eq:mainintegralvar}, we will prove
the following:
\begin{proposition}\label{prop:charavar}
  With the notation as above, for $b\geq 0$, we have the following recursive
  relation:
  \begin{equation}\label{eq:recvar}
    \begin{split}
      & [i^{*}J^{X,tw}_{\beta,\vec{k},c}(\mu^{X},z)]_{z^{-b-1}}\\
      &=\bigg[\sum_{m=0}^{\infty}\sum_{n=0}^{\infty} \sum_{\substack{\Gamma\in
          \Lambda_{\beta,\overrightarrow{k},c,m,n}\\ \Gamma\; \textrm{is
            stable}}}\frac{1}{m!n!}(\widetilde {ev_{\star}})_{*}\bigg(\sum
      _{d=0}^{\infty}\epsilon_{*}\big(c_{d}(-R^{\bullet}\pi_{*}f^{*}L^{\frac{1}{r}})(\frac{\lambda}{r})^{-1+m-d}(-1)^{d}\\
      &\cap [\mathcal K_{0,\overrightarrow{m+n}\cup
        \star}(\sqrt[r]{L/Y},\beta_{\star})]^{\mathrm{vir}}\big) \cap
      \prod_{i=1}^{m}\frac{ev_{i}^{*}(\frac{1}{\delta_i} i^{*}\mathbf{t}^{\overrightarrow{k_{i}}}J^{X,tw}_{\beta_{i},\vec{k_i},c_{i}}(\mu^{X},z)|_{z=\frac{\lambda-c_1(L)}{\delta_{i}}})}{\frac{\lambda-ev^{*}_{i}c_1(L)}{r\delta_i}+\frac{\bar{\psi_{i}}}{r}}\\
      &\cap
      \prod_{i=m+1}^{m+n}ev_{i}^{*}(\mathbf{t}^{\overrightarrow{k_i}}\mu^{X}_{\beta_{i},\vec{k_i},c_{i}}(-\bar{\psi}_{i}))
      \cap \bar{\psi}_{\star}^{b}\bigg)\bigg]_{\mathbf{t}^{\overrightarrow{k}}\lambda^{-1}}\ .
    \end{split}
  \end{equation}
  Here $\delta_{i}=\beta_{i}(L)+ (\vec{w}, \overrightarrow{k_{i}})$, and
  $\epsilon:\mathcal K_{0,\overrightarrow{m+n}\cup
    \star}(\sqrt[r]{L/Y},\beta_{\star})\rightarrow \mathcal
  K_{0,\overrightarrow{m+n}\cup \star}(Y,\beta_{\star})$ is the natural
  structural morphism by forgetting root-gerbe structure of $\sqrt[r]{L/Y}$
  (c.f.\cite{tang2021quantum}), where we choose the
  tuple $\overrightarrow{m+n}\cup\star$ for $\mathcal
  K_{0,\overrightarrow{m+n}\cup \star}(\sqrt[r]{L/Y},\beta_{\star})$ to be
  $$\big((c_{1}^{-1},1,\mathrm{e}^{\frac{-\delta_1}{r}}),\cdots, (c_{m+n}^{-1},1,\mathrm{e}^{\frac{-\delta_{m+n}}{r}}), (c,1,\mathrm{e}^{\frac{\delta}{r}})\big)\in
  (C \!\times\! \C^{*} \!\times\! \C^{*})^{m+n+1}\ ,$$ and choose the tuple
  $\overrightarrow{m+n}\cup\star$ for $\mathcal K_{0,\overrightarrow{m+n}\cup
    \star}(Y,\beta_{\star})$ to be
    $$(c_{1}^{-1},\cdots, c_{m+n}^{-1},c)\in C^{m+n+1}\ .$$ 
  \end{proposition}

  Our main strategy to prove this is to apply $\C^{*}-$localization to the
  integral \eqref{eq:mainintegralvar} and then extract the
  $\lambda^{-2}$-coefficient(or equivalently $\mathbf{t}^{\overrightarrow{k}}\lambda^{-2}$-coefficient). Using the polynomiality of
  \eqref{eq:mainintegralvar}, we have that the $\lambda^{-2}$ coefficient must
  be zero, which yields the desired result. To simplify the localization
  calculation, we need the following two lemmas which limit the types of
  localization graphs one need to consider.

\begin{lemma}\label{lem:edge-exlusion}
  Let $\Gamma$ be a decorated graph for $\mathcal
  K_{0,\overrightarrow{m}\cup\star}(\mathfrak R, (\beta,\frac{\delta}{r},0))$.
  If the corresponding fixed loci $F_{\Gamma}$ is non-empty, then there is no edge of $\Gamma$ linking
  $\mathcal D_{\infty}$ and $X_{\infty}$ and no edge of $\Gamma$ linking $X_{0}$
  and $X_{\infty}$.
\end{lemma}
\begin{proof}
  Let $f:C\rightarrow \mathfrak R$ be a $\C^{*}-$fixed stable map in $\mathcal
  K_{0,\overrightarrow{m}\cup\star}(\mathfrak R, (\beta,\frac{\delta}{r},0))$,
  we have the pairing $(f_{*}([C]), [X_{\infty}])=0$. For each stable vertex
  $v$, obviously, we have the pairing $(f_{*}([C_{v}]), [X_{\infty}])=0$. If $e$
  is an edge, there three possibilities:
  \begin{enumerate}
  \item The edge $e$ links a vertex labeled by $\mathcal D_{\infty}$ and
    $X_{\infty}$.
  \item The edge $e$ links a vertex labeled by $\mathcal D_{\infty}$ and
    $X_{0}$.
  \item The edge $e$ links a vertex labeled by $ X_{0}$ and $X_{\infty}$.
  \end{enumerate}
  In the first two cases, we have the pairing $(f_{*}([C_{e}]),
  [X_{\infty}])>0$, while in the third case, the pairing $(f_{*}([C_{e}]),
  [X_{\infty}])=0$. As the total degree of $f^{*}\mathcal O(X_{\infty})$ on $C$
  is zero, we see that the degree of the restriction of $f^{*}\mathcal
  O(X_{\infty})$ to each component $C_{v}$ or $C_{e}$ must be zero, which
  implies that there is no edge linking $\mathcal D_{\infty}$ and $X_{\infty}$
  and no edge linking $X_{0}$ and $X_{\infty}$.
\end{proof}
\begin{remark}\label{rmk:edge-exclusion}
The above lemma also works for the case when $\delta=0$, and the proof is verbatim.
\end{remark}

By the above the discussion, they will only be two types of graphs which will
contribute to the integral \eqref{eq:mainintegralvar}: (1) all edges in $\Gamma$
connects $X_{0}$ and $\mathcal D_{\infty}$; (2) $\Gamma$ has only one vertex
labeled by $X_{\infty}$ and no edges. Now let's assume the graph $\Gamma$ is of
type 1, we have the following lemma by using a similar argument as in Lemma
\ref{lem:vanish1} by using the line bundle $N:=f^{*}(\mathcal O(-P_{Y}))$.
\begin{lemma}\label{lem:vanish2}
  Assume $r,s$ is sufficiently large. If localization graph $\Gamma$ has more
  than one vertex labeled by $\infty$, then the corresponding $\C^{*}-$fixed
  loci moduli $F_{\Gamma}$ is empty, therefore it will contribute zero to
  \eqref{eq:mainintegralvar}.
\end{lemma}

Now we are ready to prove Proposition \ref{prop:charavar}.
\begin{proof}
  If a localization graph $\Gamma$ has a vertex labeled by $X_{\infty}$, by
  Lemma \ref{lem:edge-exlusion}, there will be no edges in $\Gamma$. Then
  $\Gamma$ is made of a single vertex. Such type of graph will contribute
  \begin{equation}\label{eq:typeextra}
    \begin{split}
      \frac{-1}{\lambda}\sum_{m\geq 0}
      \sum_{\substack{(\beta_{j},\overrightarrow{k_j},c_{j})_{j=1}^{m}\in
          \mathrm{Adm}^{m}\\ \beta_{\star}\in \mathrm{Eff}(X)\\
          \beta_{1}+\cdots+\beta_{m}+\beta_{\star}=\beta\\
          \overrightarrow{k_{1}}+\cdots+\overrightarrow{k_{m}}=\overrightarrow{k}
        }} \frac{1}{m!}\phi^{\alpha}\langle{
        \mathbf{t}^{\overrightarrow{k_{1}}}\mathbf{\mu^{X}_{\beta_{1},\overrightarrow{k_{1}},c_{1}}}(-\bar{\psi}_{1}),\cdots,\mathbf{t}^{\overrightarrow{k_{m}}}\mathbf{\mu^{X}_{\beta_{m},\overrightarrow{k_{m}},c_{m}}}(-\bar{\psi}_{m}),\frac{\phi_{\alpha}}{z-\bar{\psi}_{\star}}}\rangle^{X_{\infty}}_{0,\vec{m}\cup\star,\beta_{\star}}
    \end{split}
  \end{equation}
  to \ref{eq:mainintegralvar}.

  The other type of graph which contributes to \eqref{eq:recvar} must satisfy
  that all the vertexes are labeled by either $0$ or $\infty$ and all the edges
  only link $X_{0}$ and $\mathcal D_{\infty}$. Using Lemma \ref{lem:vanish2}, we
  can apply the same strategy of the Proposition \ref{prop:deg-spe} to do the
  localization computation. Denote $v_{\star}$ to be the only vertex of $\Gamma$
  labeled by $\infty$, then the remaining types of localization graph is a
  star-shaped graph introduced in the proof of \ref{prop:deg-spe} and can be
  sorted as below:
  \begin{enumerate}
  \item(Type I) The vertex $v_{\star}$ is unstable and there is only one edge
    incident to $v_{\star}$ in $\Gamma$, the vertex $v$ labeled by $0$ is
    unstable;
  \item(Type II) The vertex $v_{\star}$ is unstable and there is only one edge
    incident to $v_{\star}$ in $\Gamma$, the vertex $v$ labeled by $0$ is
    stable;
  \item(Type III) The vertex $v_{\star}$ is stable.
  \end{enumerate}
  The graph of type I will contribute
  \begin{equation}\label{eq:type1new}
    \frac{1}{\delta\lambda}\mathbf{t}^{\overrightarrow{k}}\prod_{0\leq j<\delta}\big(c_{1}(L)+(\delta-j)\frac{\lambda-c_1(L)}{\delta}\big)i^{*}\mu_{\beta,\vec{k},c}^{X}(\frac{\lambda-c_1(L)}{\delta})\cdot(\frac{\lambda-c_1(L)}{\delta})^{b} 
  \end{equation}
  to \ref{eq:mainintegralvar}.
  

  The graph of type II will contribute
  \begin{equation}\label{eq:type2new}
    \begin{split}
      &\frac{1}{\delta\lambda}\prod_{0\leq
        j<\delta}\big(c_{1}(L)+(\delta-j)\frac{\lambda-c_1(L)}{\delta}\big)
      \sum_{m\geq
        0} \sum_{\substack{(\beta_{i},\overrightarrow{k_i},c_{i})_{i=1}^{m}\in \mathrm{Adm}^{m}\\ \beta_{\star}\in \mathrm{Eff}(X),\; \sum_{i}\beta_{i}+\beta_{\star}=\beta\\ \sum_{i}\overrightarrow{k_{i}}=\overrightarrow{k}  }} \frac{1}{m!}\\
      &\cdot i^{*}\phi^{\alpha}\bigg(\int_{ [\mathcal
        K_{0,\vec{m}\cup\star}(X,\beta_{\star})]^{\mathrm{vir}}}
      \prod_{i=1}^{m}ev_{i}^{*}(\mathbf{t}^{\overrightarrow{k_{i}}}\mu^{X}_{\beta_{i},\overrightarrow{k_{i}},c_{i}}(-\bar{\psi}_{i}))
      \cup \frac{ev^{*}_{\star}\phi_{\alpha}
        (\frac{\lambda-ev_{\star}^{*}c_1(L)}{\delta})^{b}
      }{\frac{\lambda-ev_{\star}^{*}c_1(L)}{\delta}-\bar{\psi_{\star}}} \bigg)\
      .
    \end{split}
  \end{equation}
  to \ref{eq:mainintegralvar}.


  The graph of type III will contribute
  \begin{equation}\label{eq:type3}
    \begin{split}
      &\frac{1}{\lambda}\sum_{m=0}^{\infty}\sum_{n=0}^{\infty}
      \sum_{\substack{\Gamma\in \Lambda_{\beta,\overrightarrow{k},c,m,n}\\
          \Gamma\; \textrm{is
            stable}}}\frac{1}{m!}\frac{1}{n!}(\widetilde{ev_{\star}})_{*}\bigg(\sum
      _{d=0}^{\infty}\epsilon_{*}\big(c_{d}(-R^{\bullet}\pi_{*}f^{*}L^{\frac{1}{r}})(\frac{-\lambda}{r})^{-1+m-d}\\
      &\cap [\mathcal K_{0,\overrightarrow{m+n}\cup
        \star}(\sqrt[r]{L/Y},\beta_{\star})]^{\mathrm{vir}}\big) \cap
      \prod_{i=1}^{m}\frac{ev_{i}^{*}(\frac{1}{\delta_i}i^{*}\mathbf{t}^{\overrightarrow{k_{i}}}J^{X,tw}_{\beta_{i},\vec{k_i},c_{i}}(z)|_{z=\frac{\lambda-c_1(L)}{\delta_i}})}{-\frac{\lambda-ev^{*}_{i}c_1(L)}{r\delta_i}-\frac{\bar{\psi_{i}}}{r}}\\
      &\cap
      \prod_{i=m+1}^{m+n}ev_{i}^{*}(\mathbf{t}^{\overrightarrow{k_{i}}}\mu^{X}_{\beta_{i},\vec{k_i},c_{i}}(-\bar{\psi}_{i}))\cap
      \bar{\psi}^{b}_{\star}\bigg)\ .
    \end{split}
  \end{equation}
  to \ref{eq:mainintegralvar}.

  Now our auxiliary cycle \ref{eq:mainintegralvar} is equal to the sum of
  \ref{eq:typeextra}, \ref{eq:type1new}, \ref{eq:type2new} and \ref{eq:type3}.
  By the polynomiality of \ref{eq:mainintegralvar}, the
  $\mathbf{t}^{\overrightarrow{k}}\lambda^{-2}-$coefficient of the auxiliary cycle is zero, which implies that
  the sum of $\mathbf{t}^{\overrightarrow{k}}\lambda^{-2}-$coefficient of \ref{eq:typeextra},
  \ref{eq:type1new}, \ref{eq:type2new} and \ref{eq:type3} is zero, which proves
  the proposition \ref{prop:charavar}.
\end{proof}

\subsection{Statement of the main theorem and proof}\label{subsec:mainthm}



With the notations in this section, now we state our main theorem:
\begin{theorem}[Main theorem]\label{thm:main1}
  Let $X$ be a smooth Deligne-Mumford stack with projective coarse moduli. Let
  $E:=\oplus_{j=1}^{r}L_{j}$ be a direct sum of \emph{semi-positive} line
  bundles over $X$ with a regular section of $E$ which cuts off a smooth
  complete intersection $Y\subset X$. Let
  $\mu^{X}(z):=\sum_{(\beta,\overrightarrow{k},c)\in
    \mathrm{Adm}}q^{\beta}\mathbf{t}^{\overrightarrow{k}}\mu^{X}_{\beta,\vec{k},c}(z)$ be an admissible series
  and $J^{X,tw}(q,\mu^{X},z)$ be the hypergeometric modification of
  $J^{X}(q,\mu^{X},z)$ as in \ref{eq:tw}. Let $i:\bar{I}_{\mu}Y\rightarrow
  \bar{I}_{\mu}X$ be the inclusion, then the series
  $i^{*}J^{X,tw}(q,\mu^{X},-z)$ is a point on the Lagrange cone of $Y$. More
  precisely, for any admissible pair $(\beta,\overrightarrow{k},c)$, if we
  define
 $$\mu^{X,tw}_{\beta,\vec{k},c}(z):= \big[J^{X}_{\beta,\vec{k},c}(\mu^{X},z)\prod_{j=1}^{r} \prod_{0\leq
   m<\beta(L_{j})+(\vec{w_{j}},\vec{k})}\big(c_{1}(L_{j}) +(\beta(L_{j})+(\vec{w_{j}},\vec{k})-m)z
 \big)\big]_{+} $$ to be the truncation in nonnegative $z-$powers. Here
 $J^{X}_{\beta,\vec{k},c}$ is defined in \ref{eq:jcomp}. Write
$$ \mu^{X,tw}(q,z):= \sum_{(\beta,\overrightarrow{k},c)\in
  \mathrm{Adm} }q^{\beta}\mathbf{t}^{\vec{k}}\mu^{X,tw}_{\beta,\vec{k},c}(z)\ .$$ Then for any
integer $b\geq 0$, we have the following relation:
\begin{equation}\label{eq:recmain}
  \begin{split}
    [J^{Y}_{\beta,\vec{k},c}(i^{*}\mu^{X,tw},z)]_{z^{-b-1}}=[i^{*}J^{X,tw}_{\beta,\vec{k},c}(\mu^{X},z)]_{z^{-b-1}}\
    .
  \end{split}
\end{equation}
\end{theorem}

As $\mu^{X,tw}$ is also an admissible series in $H^{*}(\bar{I}_{\mu}X,\mathbb
C)[z][\![t_{1},\cdots,t_{N}]\!][\![\mathrm{Eff}(X)]\!]$, by induction on the
number of line bundles, it's sufficient to prove the case when $Y$ is a
\emph{hypersurface}. Now let's assume that $Y$ is a hypersurface and use the
notations in previous sections.

We will first prove \ref{eq:recmain} for any admissible pair
$(\beta,\overrightarrow{k},c)$ with $\beta(L)+(\vec{w},\vec{k})=0$. Note in this
case we have $\mu^{X,tw}_{\beta,\vec{k},c}(z)=\mu^{X}_{\beta,\vec{k},c}(z)$.

We will adopt the convention that for any substack $Z$ of $\mathfrak R$ and
$c\in C$, we will write $\bar{I}_{c}Z := \bar{I}_{\mu}Z\cap \bar{I}_{(c,1,1)}Z$.
Now fix an admissible pair $(\beta,\overrightarrow{k},c)$ such that
$\beta(L)+(\vec{w},\vec{k})=0$, we have $age_{c}L=0$. Then $\bar{I}_{c}Y$ is
proper hypersurface of $\bar{I}_{c}X $ with normal bundle $L$.\footnote{Here the
  normal bundle is a pull-back of the line bundle $L$ along the natural morphism
  from $\bar{I}_{\mu}X$ to $X$ by sending $(x,g)$ to $x$.} Recall that the
divisor $E$ of $\mathfrak R$ is isomorphic to $\P Y_{r,s}$, Denote by
$\mathrm{pr}^{E}_{r,s}: E\rightarrow Y$ the projection to the base, which
induces a morphism on the corresponding rigidified inertia stacks, which we
still denote to be $\mathrm{pr}^{E}_{r,s}$. We have $\bar{I}_{c^{-1}}E\cong
(\mathrm{pr}^{E}_{r,s})^{-1}(\bar{I}_{c^{-1}}Y ) $ and $\bar{I}_{c^{-1}}E$ is a
hypersurface of $\bar{I}_{c^{-1}}\mathfrak R$.

Denote by $i_{E}$ the inclusion from $\bar{I}_{c^{-1}}E $ to $\bar{I}_{c^{-1}}
\mathfrak R$. For any admissible pair $(\beta,\overrightarrow{k},c)$ with
$\beta(L)+(\vec{w},\vec{k})=0$, consider the following class
\begin{equation}\label{eq:int-semipos}
  \begin{split}
    &\sum_{m=0}^{\infty} \sum_{\substack{ (\beta_{i},\overrightarrow{k_i},c_{i})_{i=1}^{m}\in \mathrm{Adm}^{m}\\ \beta_{\star}+\sum_{i=1}^{m}\beta_{i}=\beta,\; \beta_{\star}\in \mathrm{Eff}(X)\\ \overrightarrow{k_{1}}+\cdots+\overrightarrow{k_{m}}=\overrightarrow{k}}}\frac{1}{m!}(\mathrm{pr}^{E}_{r,s})_{*}\circ i^{*}_{E} \circ (\widetilde{ev_{\star}})_{*}\\
    &\bigg(\mathcal K_{0,\overrightarrow{m}\cup \star}(\mathfrak R,
    (\beta_{\star},0,0))\cap \bar{\psi}_{\star}^{b} \cap \prod_{i=1}^{m}
    \mathrm{pr}^{*}_{r,s}(\mathbf{t}^{\overrightarrow{k_i}}\mu^{X}_{\beta_{i},\vec{k_i},c_{i}}(-\bar{\psi}_{i}))\bigg)\ .
  \end{split}
\end{equation}
Here for any $\beta_{\star},\beta_{1},\cdots, \beta_{m}$,
$\overrightarrow{k_{1}}\cdots,\overrightarrow{k_{m}}$ and $c_{1},\cdots,c_{m}$
such that $(\beta_{i},\overrightarrow{k_i},c_{i})_{i=1}^{m}\in \mathrm{Adm}^{m}$
with $\sum_{i=1}^{m}\beta_{i}+\beta_{\star}=\beta$ and
$\overrightarrow{k_{1}}+\cdots+\overrightarrow{k_{m}}=\overrightarrow{k}$, we
associated an element in $(C \!\times\! \C^{*} \!\times\! \C^{*})^{m+1}$
$$\overrightarrow{m}\cup
\star:=\big((c_{1}^{-1},1,1),\cdots,(c_{m}^{-1},1,1,),(c,1,1)\big)\ .$$
$i_{E}^{*}$ is the Gysin pull-back from the (equivariant) cohomology of
$\bar{I}_{(c^{-1},1,1)} \mathfrak R$ to the (equivariant) cohomology of
$\bar{I}_{(c^{-1},1,1)}E$. Note that $\beta_{i}(L)+(\vec{w},\vec{k_{i}})=0$ for
all $1\leq i\leq m$,

\begin{lemma}\label{lem:deg0main}
  For any admissible pair $(\beta,\overrightarrow{k},c)$ with
  $\beta(L)+(\vec{w},\vec{k})=0$, let $b$ be an nonnegative integer, we have
  $$[J^{Y}_{\beta,\vec{k},c}(q,i^{*}\mu^{X},z)]_{z^{-b-1}}=
  i^{*}[J^{X}_{\beta,\vec{k},c}(q,\mu^{X},z)]_{z^{-b-1}}\ ,$$ where $i^{*}$ is
  induced from the inclusion from $\bar{I}_{\mu}Y$ to $\bar{I}_{\mu}X$.
\end{lemma}

\begin{proof}
  We will apply $\C^{*}-$localization to \ref{eq:int-semipos}. First we claim
  there are only two types of localization graphs contributing to the integral
  \eqref{eq:int-semipos}: 1, the localization $\Gamma_{1}$ has only one vertex
  and it's labeled by $X_{0}$; 2, the localization graph $\Gamma_{2}$ has only
  one vertex and it's labeled by $\mathcal D_{\infty}$. Indeed, there is no
  vertex labeled by $X_{\infty}$ as the marking $q_{\star}$ must go to the
  divisor $E$ to make nonzero localization contribution to \ref{eq:int-semipos}.
  Furthermore, as $deg(f^{*}\mathcal O(P_{Y}))=0$, we claim there is no edge
  linking $X_{0}$ and $\mathcal D_{\infty}$. Indeed, recall that we only need to
  consider edge curve $C_{e}$ linking $X_{0}$ and $\mathcal D_{\infty}$ by Lemma
  \ref{lem:edge-exlusion} (see also Remark \ref{rmk:edge-exclusion}). then the claim comes from that the line bundle
  $\mathcal O(P_{Y})$ is of nonnegative degree restricting to all vertex curves
  $C_{v}$ and of positive degree restricting to all edge curves $C_{e}$. The
  localization contribution from $\Gamma_{1}$ is equal to
  \begin{equation}\label{eq:deg0amb}
    \begin{split}
      &(\mathrm{pr}^{E}_{r,s})_{*}\circ i^{*}_{E}\circ
      (\iota_{X_{0}})_{*}\bigg[\sum_{m=0}^{\infty} \sum_{\substack{
          (\beta_{i},\overrightarrow{k_i},c_{i})_{i=1}^{m}\in
          \mathrm{\mathrm{Adm}}^{m}\\
          \beta_{\star}+\sum_{i=1}^{m}\beta_{i}=\beta,\; \beta_{\star}\in
          \mathrm{Eff}(X)\\
          \overrightarrow{k_{1}}+\cdots+\overrightarrow{k_{m}}=\overrightarrow{k}}}\frac{1}{m!}(\widetilde
      {ev_{\star}})_{*}\bigg(\sum
      _{d=0}^{\infty}\big(c_{d}(-R^{\bullet}\pi_{*}f^{*}(L^{\vee})^{\frac{1}{s}})(\frac{\lambda}{s})^{-1-d}\\
      &\cap [\mathcal K_{0,\vec{m}\cup
        \star}(X_{0},\beta_{\star})]^{\mathrm{vir}}\big) \cap
      \prod_{i=1}^{m}ev_{i}^{*}(\mathrm{pr}^{*}_{r,s}\mathbf{t}^{\overrightarrow{k_{i}}}\mu^{X}_{\beta_{i},\vec{k_i},c_{i}}(-\bar{\psi}_{i}))
      \cap \bar{\psi}_{\star}^{b}\bigg)\bigg].
    \end{split}
  \end{equation}
  Here $\iota_{X_{0}}:\bar{I}_{c^{-1}}X_{0}\rightarrow \bar{I}_{c^{-1}}\mathfrak
  R$ is the natural inclusion and $(\widetilde{ev_{\star}})_{*}$ is a morphism
  from $H^{*}_{\C^{*}}(\mathcal K_{0,\vec{m}\cup \star}(X_{0},\beta_{\star}))$
  to $H^{*}_{\C^{*}}(\bar{I}_{c^{-1}}X_{0})$ as defined in \ref{tilde-ev}.

  We have the following commutative diagram where the upper square and outer
  square are Cartersian:
$$\xymatrix{
  \bar{I}_{c^{-1}}\mathcal D_{0}\ar[r]^{i_{\mathcal
      D_{0}}}\ar[d]^{\iota_{\mathcal D_{0}}} &\bar{I}_{c^{-1}}X_{0} \ar[d]^{\iota_{X_{0}}}\\
  \bar{I}_{c^{-1}}E\ar[r]^{i_{E}}\ar[d]^{\mathrm{pr}^{E}_{r,s}}
  &\bar{I}_{c^{-1}}\mathfrak R\ar[d]^{\mathrm{pr}_{r,s}} \\
  \bar{I}_{c^{-1}}Y\ar[r]^{i} & \bar{I}_{c^{-1}}X\ . }$$ Here $i$ is induced
from the inclusion from $\mathcal D_{0}$ to $X_{0}$ and $\iota_{\mathcal D_{0}}$
is induced from the inclusion from $\mathcal D_{0}$ to $E$. As the age
$age_{c}(L)=0$, all rows are regular embeddings of codimension one. Then by the
commutativity of the proper push-forward and Gysin pullback for the upper
square, \ref{eq:deg0amb} is equal to
\begin{equation}\label{eq:deg0amb1}
  \begin{split}
    &(\mathrm{pr}^{E}_{r,s}\circ \iota_{\mathcal
      D_{0}})_{*}i^{*}_{E}(\iota_{\mathcal D_{0}})_{*}\big[\cdots\big].
  \end{split}
\end{equation}
Here the dots inside the big square bracket copies the stuff in the big square
bracket of \ref{eq:deg0amb}. Using the outer square, \ref{eq:deg0amb1} is equal
to
\begin{equation}\label{eq:deg0amb2}
  \begin{split}
    & i^{*}_{Y}(\mathrm{pr}_{r,s}\circ \iota_{X_{0}})_{*}\big[\cdots\big].
  \end{split}
\end{equation}

As the composition $\mathrm{pr}_{r,s}\circ \iota_{X_{0}}$ is induced from
natural de-root stackfication of $X_{0}\cong \sqrt[s]{L^{\vee}/X} $,
$(\mathrm{pr}_{r,s}\circ \iota_{X_{0}})_{*}$ is the identity map on cohomology
groups as it induces an identity map on their coarse moduli spaces and the cohomology group of an orbifold is
identified with the cohomology group of its coarse moduli space
(see~\cite{Dan_Abramovich_2008} for more details). Therefore \ref{eq:deg0amb2}
is equal to
\begin{equation}\label{eq:deg0amb3}
  \begin{split}
    &i^{*}\bigg[\sum_{m=0}^{\infty} \sum_{\substack{
        (\beta_{i},\overrightarrow{k_i},c_{i})_{i=1}^{m}\in
        \mathrm{\mathrm{Adm}}^{m}\\
        \beta_{\star}+\sum_{i=1}^{m}\beta_{i}=\beta,\; \beta_{\star}\in
        \mathrm{Eff}(X) \\
        \overrightarrow{k_{1}}+\cdots+\overrightarrow{k_{m}}=\overrightarrow{k}}}\frac{1}{m!}(\widetilde
    {ev_{\star}})_{*}\bigg(\sum
    _{d=0}^{\infty}\epsilon_{*}\big(c_{d}(-R^{\bullet}\pi_{*}f^{*}(L^{\vee})^{\frac{1}{s}})(\frac{\lambda}{s})^{-1-d}\\
    &\cap [\mathcal K_{0,\vec{m}\cup
      \star}(X_{0},\beta_{\star})]^{\mathrm{vir}}\big) \cap
    \prod_{i=1}^{m}ev_{i}^{*}(\mathbf{t}^{\overrightarrow{k_{i}}}\mu^{X}_{\beta_{i},\vec{k_i},c_{i}}(-\bar{\psi}_{i}))
    \cap \bar{\psi}_{\star}^{b}\bigg)\bigg],
  \end{split}
\end{equation}
where $\epsilon:\mathcal K_{0,\vec{m}\cup\star}(X_{0},\beta_{\star})\rightarrow
\mathcal K_{0,\vec{m}\cup\star}(X,\beta_{\star})$ is induced from the natural structural map
forgetting the root of $X_{0}\cong \sqrt[s]{L^{\vee}/X}$.

The localization contribution from $\Gamma_{2}$ is equal to
\begin{equation}\label{eq:deg0hyper}
  \begin{split}
    & \frac{1}{\lambda}(\mathrm{pr}^{E}_{r,s})_{*}\circ i^{*}_{E}\circ
    (i_{\mathcal D_{\infty}})_{*}\bigg[\sum_{m=0}^{\infty} \sum_{\substack{
        (\beta_{i},\overrightarrow{k_i},c_{i})_{i=1}^{m}\in \mathrm{Adm}^{m}\\
        \beta_{\star}\in \mathrm{Eff}(X),\;
        \beta_{\star}+\sum_{i=1}^{m}\beta_{i}=\beta\\
        \overrightarrow{k_{1}}+\cdots+\overrightarrow{k_{m}}=\overrightarrow{k}}}\frac{1}{m!}(\widetilde
    {ev_{\star}})_{*}\bigg(\sum
    _{d=0}^{\infty}\big(c_{d}(-R^{\bullet}\pi_{*}f^{*}L^{\frac{1}{r}})(\frac{-\lambda}{r})^{-1-d}\\
    &\cap [\mathcal K_{0,\vec{m}\cup
      \star}(\sqrt[r]{L/Y},\beta_{\star})]^{\mathrm{vir}}\big) \cap
    \prod_{i=1}^{m}ev_{i}^{*}(\mathrm{pr}_{r,s}^{*}\mathbf{t}^{\overrightarrow{k_{i}}}\mu^{X}_{\beta_{i},\vec{k_i},c_{i}}(-\bar{\psi}_{i}))
    \cap \bar{\psi}_{\star}^{b} \bigg)\bigg]
  \end{split}
\end{equation}
Here $i_{\mathcal D_{\infty}}:\bar{I}_{c^{-1}}\mathcal D_{\infty}\rightarrow
\bar{I}_{c^{-1}}\mathfrak R$ is the inclusion and $(\widetilde{ev_{\star}})_{*}$
is the morphism from $H^{*}(\mathcal K_{0,\vec{m}\cup \star}(\mathcal
D_{\infty},\beta_{\star}))$ to $H^{*}(\bar{I}_{c^{-1}}\mathcal D_{\infty})$ as
defined in \ref{tilde-ev}. The inclusion $i_{\mathcal D_{\infty}}$ can be
rewritten as the composition of two inclusions $\iota_{\mathcal
  D_{\infty}}:\bar{I}_{c^{-1}}\mathcal D_{\infty}\hookrightarrow
\bar{I}_{c^{-1}}E$ and $i_{E}: \bar{I}_{c^{-1}}E\hookrightarrow
\bar{I}_{c^{-1}}\mathfrak R $. Then we have
\begin{equation}\label{eq:deg0hyper1}
  \begin{split}
    &\frac{1}{\lambda}(\mathrm{pr}^{E}_{r,s})_{*}\circ i^{*}_{E}\circ (i_{\mathcal D_{\infty}})_{*}\bigg[\cdots\bigg]\\
    &=\frac{1}{\lambda}(\mathrm{pr}^{E}_{r,s})_{*}\circ i^{*}_{E}\circ (i_{E})_{*}\circ (\iota_{\mathcal D_{\infty}})_{*}\bigg[\cdots\bigg]\\
    &=(\mathrm{pr}^{E}_{r,s})_{*}(\iota_{\mathcal D_{\infty}})_{*}\bigg[\cdots\bigg]\\
    &=\bigg[\cdots \bigg]\ .
  \end{split}
\end{equation}
Here the dots inside the big square bracket copies the stuff in the big square
bracket of \ref{eq:deg0hyper}. The second equality in \ref{eq:deg0hyper1} uses
the fact that the morphism
$i^{*}_{E}(i_{E})_{*}:H^{*}_{\C^{*}}(\bar{I}_{c^{-1}}E)\rightarrow
H^{*}_{\C^{*}}(\bar{I}_{c^{-1}}E)$ is equal to the multiplication map by
$e^{\C^{*}}(N_{\bar{I}_{c^{-1}}Y/\bar{I}_{c^{-1}}\mathfrak R })$ using excess
intersection formula, which implies that $i^{*}_{E}(i_{E})_{*}$ acts on the
space $(\iota_{\mathcal D_{\infty}})_{*}H^{*}_{\C^{*}}(\mathcal D_{\infty})$ by
multiplication by $\lambda$ as
$e^{\C^{*}}(N_{\bar{I}_{c^{-1}}Y/\bar{I}_{c^{-1}}\mathfrak R
})|_{|_{\bar{I}_{c^{-1}}\mathcal D_{\infty}}}=\lambda$. Finally
\ref{eq:deg0hyper} becomes
\begin{equation}\label{eq:deg0hyper2}
  \begin{split}
    &\sum_{m=0}^{\infty} \sum_{\substack{
        (\beta_{i},\overrightarrow{k_i},c_{i})_{i=1}^{m}\in
        \mathrm{Adm}^{m},\beta_{\star}\in \mathrm{Eff}(X)\\
        \beta_{\star}+\sum_{i=1}^{m}\beta_{i}=\beta\\
        \overrightarrow{k_{1}}+\cdots+\overrightarrow{k_{m}}=\overrightarrow{k}}}\frac{1}{m!}(\widetilde
    {ev_{\star}})_{*}\bigg(\sum
    _{d=0}^{\infty}\epsilon_{*}\big(c_{d}(-R^{\bullet}\pi_{*}f^{*}L^{\frac{1}{r}})(\frac{-\lambda}{r})^{-1-d}\\
    &\cap [\mathcal K_{0,\vec{m}\cup
      \star}(\sqrt[r]{L/Y},\beta_{\star})]^{\mathrm{vir}}\big) \cap
    \prod_{i=1}^{m}ev_{i}^{*}(i^{*}\mathbf{t}^{\overrightarrow{k_{i}}}\mu^{X}_{\beta_{i},\vec{k_i},c_{i}}(-\bar{\psi}_{i}))
    \cap \bar{\psi}_{\star}^{b} \bigg)
  \end{split}
\end{equation}
using the morphism $\epsilon:\mathcal
K_{0,\vec{m}\cup\star}(\mathcal D_{\infty},\beta_{\star})\rightarrow \mathcal
K_{0,\vec{m}\cup\star}(Y,\beta_{\star})$ induced from forgetting the root of
$\mathcal D_{\infty}\cong \sqrt[s]{L/Y}$.

Now taking $\lambda^{-1}$ coefficients of the contributions \ref{eq:deg0amb3}
and \ref{eq:deg0hyper2} from $\Gamma_{1}$ and $\Gamma_{2}$, their sum is zero by
polynomiality of \ref{eq:int-semipos}, this finishes the proof of
\eqref{eq:recmain} with the help of the fact
$$ \epsilon_{*}\big( [\mathcal
K_{0,\vec{m}\cup\{\star\}}(\sqrt[s]{L^{\vee}/X},\beta_{\star})]^{\mathrm{vir}}\big)=\frac{1}{s}
[\mathcal K_{0,\vec{m}\cup\{\star\}}(X,\beta_{\star})]^{\mathrm{vir}}$$
$$ \epsilon_{*}\big( [\mathcal
K_{0,\vec{m}\cup\{\star\}}(\sqrt[r]{L/Y},\beta_{\star})]^{\mathrm{vir}}\big)=\frac{1}{r}
[\mathcal K_{0,\vec{m}\cup\{\star\}}(Y,\beta_{\star})]^{\mathrm{vir}} $$
proved in~\cite{tang2021quantum}.
\end{proof}

Now we are well prepared to prove the main theorem \ref{thm:main1}:
\begin{proof}[Proof of the main theorem \ref{thm:main1}]
  Since $X$ has projective coarse moduli, there exists a \emph{positive} line
  bundle $M$ on $X$, i.e., $\beta(M)>0$ for all nonzero degree $\beta\in
  \mathrm{Eff}(X)$. Let's prove the relation \ref{eq:recmain} by induction on
  the nonnegative rational number $e:=\beta(M)+ |\vec{k}|$.\footnote{Here
    $|\vec{k}|=\sum_{i}k_{i}$.} When $e=0$, we have
  $(\beta,\vec{k})=(0,\vec{0})$, then
  $J^{Y}_{0,\vec{0},c}(i^{*}\mu^{X,tw},z)=J^{X,tw}_{0,\vec{0},c}(\mu^{X},z)=0$.
  Now suppose the relation \ref{eq:recmain} holds for all admissible pair
  $(\beta',\overrightarrow{k'},c')$ with $\beta'(M)+ |\overrightarrow{k'}|<e$
  for some positive rational number $e$, it remains to show the relation holds
  for any admissible pair $(\beta,\overrightarrow{k},c)$ with
  $\beta(M)+|\overrightarrow{k}|=e$. If $\beta(L)+
  (\vec{w},\overrightarrow{k})=0$, by Proposition \ref{lem:deg0main}, the
  relation holds (note in this case, we have
  $J^{X,tw}_{\beta,\vec{k},c}(\mu,z)=J^{X}_{\beta,\vec{k},c}(\mu,z)$).
  Otherwise, we use the recursive relations \ref{eq:recvar} and
  \ref{eq:rec1-spe}; in which the terms corresponding to the element $\Gamma\in
  \Lambda_{\beta,\overrightarrow{k},c,m,n}$ with $\beta_{\star}\neq 0$ are the
  same as their inputs are associated with the extended degrees
  $(\beta_{i},\overrightarrow{k_{i}})$ with $\beta_{i}(M)+
  |\overrightarrow{k_{i}}|<e$, which are the same by induction assumption. We
  are left to show that the terms in \ref{eq:recvar} and \ref{eq:rec1-spe}
  corresponding to the element $\Gamma\in
  \Lambda_{\beta,\overrightarrow{k},c,m,n}$ with $\beta_{\star}=0$ are the same.
  Note $m+n\geq 2$ as $\Gamma$ is stable and
  $J^{Y}_{0,\vec{0},c}(i^{*}\mu^{X,tw},z)=J^{X,tw}_{0,\vec{0},c}(\mu^{X},z)=0$,
  then we can still use inductive assumption. This completes the proof.
\end{proof}

\begin{remark}
  The above proof can also be extended to the equivariant case, where we assume
  that $X$ has an algebraic torus $T-$action which preserves $Y$. Since we don't
  have a particular application in mind, we will leave the details to the
  interested readers.
\end{remark}

Let $\{\phi_{\alpha}\}$ be a \emph{graded} basis of $H^{*}(\bar{I}_{\mu}X,\mathbb C)$ such
that $\phi_{\alpha}\in H^{*}(\bar{I}_{c}X)$ for some $c\in C$. Let $\{t_{\alpha}\}$
be a coordinate system of $H^{*}(\bar{I}_{\mu}X,\mathbb C)$ corresponding to the
basis $\{\phi_{\alpha}\}$. We also put $\mathbb Z_{2}$-grading on $\{t_{\alpha}\}$ according to
$\phi_{\alpha}$. Then we have the following:
\begin{corollary}
  Let $i:\bar{I}_{\mu}Y \rightarrow \bar{I}_{\mu}X$ be the inclusion of
  rigidified inertia stacks and $\mathbf{t}=\sum_{\alpha}t_{\alpha}\phi_{\alpha}$. Then the series $i^{*} J^{X,tw}(q,\mathbf{t},-z)$
  lies on the Lagrangian cone $\mathcal L_{Y}$ of $Y$.
\end{corollary}

\addtocontents{toc}{\protect\setcounter{tocdepth}{1}}

\begin{appendix}
  \section{Edge contribution}
  In the proof of our main theorem, we need to apply torus localization formula
  to the moduli of stable maps to $\P Y_{r,s}$ (or to $\mathfrak R$ which
  factors through $E$), which is an orbi-$\P^{1}$ bundle
  over a stack, for which we don't find a suitable reference of fitting in the
  discussion of the
  explicit localization contribution from the torus-fixed $1$-dimensional orbits
  (named edge contribution) in our case. In this appendix, we will explain how to compute
  the edge contributions in \ref{subsubsec:edge-cntr1} and
  \ref{subsubsec:edge-cntr2}. To achieve this, we will construct a space called
  $\mathcal M_{e}$ with a family of $\C^{*}-$fixed stable maps to $\P Y_{r,s}$
  corresponding to a decorated edge. This space has the property that it allows
  a finite \'etale map to the corresponding substack of $\C^{*}-$fixed loci in
  $\mathcal K_{0,2}(\P Y_{r,s},(0,\frac{\delta}{r}))$ (or $\mathcal
  K_{0,2}(\mathfrak R,(0,\frac{\delta}{r},0))$). Then we can use $\mathcal
  M_{e}$ as a substitute for $\C^{*}-$fixed loci in the edge contribution and
  the explicit description of the family map helps us carry out the localization
  analysis concretely. To begin with, we will prove a classification result (see
  Lemma \ref{thm:classify-mor}) about $\C^{*}-$fixed stable maps to $\P Y_{r,s}$
  which map to a fiber of $\P Y_{r,s}$ over $Y$ using the relation between
  fundamental groups and $\C^{*}-$fixed stable maps to an orbi-$\P^{1}$ curve
  from~\cite{liu2020stacky}. Note that we will always assume that $r,s$ are
  sufficiently large primes in this appendix.


\subsection{Local picture}
Write $\mathfrak P:=\P Y_{r,s}$ for simplicity. For any $\mathbb C-$point $y \in
Y$, let $G_{y}$ be the residual isotropy group of $y$ in $Y$. Let
$\rho:G_{y}\rightarrow \C^{*}=Aut(L|_{y})$ be the associated representation.
Denote $\mathfrak P_{y}$ by the fiber product in the following square
$$\xymatrix{
  \mathfrak P_{y}\ar[r]\ar[d] &\mathfrak P\ar[d]^{\mathrm{pr}_{r,s}}\\
  \mathbb B G_{y}\ar[r]^{i} &Y\ , }$$ where $i:\mathbb B G_{y}\rightarrow Y$ is
the inclusion. Let $U:=\mathbb C^{2}\backslash \{0\}$, the fiber curve
$\mathfrak P_{y}$ can be represented as the quotient stack
$$[ \C^{*} \!\times\!  U/(G_{y} \!\times\!  (\C^{*})^{2} )]$$
via the action
$$ (m, x_{1},x_{2}) (g,t_{1},t_{2})=(\rho^{-1}(g)t_{1}^{-s}t_{2}^{r}m, t_{1}x_{1},t_{2}x_{2})\ ,$$
where $(m,x_{1},x_{2})\in \C^{*}\!\times\! U$, $(g,t_{1},t_{2})\in G_{y}
\!\times\! (\C^{*})^{2}$.

The $\C^{*}-$action on $\mathfrak P$ defined in \ref{subsec:space1} restricts to
be a $\C^{*}-$action on $\mathfrak P_{y}$ by scaling $x_{1}$ with weight one.
There are two $\C^{*}-$fixed points of $\mathfrak P_{y}$: (1) The point $0$
corresponding to $x_{1}=0$; (2) the point $\infty$ corresponding to $x_{2}=0$.
Now let's describe the isotropy groups of points $0$ and $\infty$, for which we
need the following definition. Let $G$ be an arbitrary finite group and $\eta$
be a character of $G$. Assume that image of $\eta$ is a finite cyclic group
$\boldsymbol{\mu_{t}}\subset \C^{*}$ of order $t$. In other words, we have the
following exact sequence from the action
 $$\xymatrix{
   1\ar[r] &\mathrm{ker}(\eta) \ar[r] &G\ar[r]^{\eta}
   &\boldsymbol{\mu}_{t}\ar[r] &1\ . }$$

 For any positive integer $i$, we define the central extension $G(\eta,i)$ of
 $G$ by the cyclic group $\boldsymbol{\mu_{i}}$:
$$ G(\eta,i):= \{(g,b)\in G\times \C^{*} | \eta(g)=b^{i} \}\ .$$

We note $G(\eta,i)$ fits in the following Cartesian digram with all rows and all
columns exact:
$$\xymatrix{
  & &\boldsymbol{\mu}_{i}\ar[r]^{=}\ar@{|->}[d] &\boldsymbol{\mu}_{i}\ar@{|->}[d]\\
  1\ar[r] &\mathrm{ker}(\eta)\ar[d]^{=} \ar[r] &G(\eta,i)
  \ar[r]^{\eta_{i}}\ar@{->>}[d]
  &\boldsymbol{\mu}_{ti}\ar[r]\ar@{->>}[d]^{a\mapsto a^{i}}
  &1\\
  1\ar[r] &\mathrm{ker}(\eta) \ar[r] &G\ar[r]^{\eta} &\boldsymbol{\mu}_{t}\ar[r]
  &1\ ,\\
}$$ where the morphism $\eta_{i}$ is induced from the projection from $G
\!\times\! \C^{*}$ to the second factor. Then the isotropy group $G_{0}$(resp.
$G_{\infty}$) of $0$(resp. $\infty$) is isomorphic to $G_{y}(\rho^{-1},s)$(resp.
$G_{y}(\rho,r)$). Denote
$$\mathfrak U_{0}:=\mathfrak P_{y} \backslash \{\infty\}\cong [\C/G_{y}(\rho^{-1},s)],\;
\text{and}\; \mathfrak U_{\infty}:=\mathfrak P_{y} \backslash \{0\}\cong
[\C/G_{y}(\rho,r)]\ , $$ and
$$\mathfrak
o_{e}:=\mathfrak P_{y}\backslash \{0,\infty \}\ . $$ Now we quote the result
from~\cite[\S 6]{liu2020stacky}(see also~\cite{noohi2006uniformization} ).
Assume that the image of $\rho$ is $\boldsymbol{\mu}_{t}$. Then the curve
$\mathfrak P_{y}$ is a $G_{e}:=\mathrm{Ker}(\rho)$-gerbe over $\P^{1}_{tr,ts}$.
More explicitly, note $\mathfrak o_{e}$ is isomorphic to the quotient stack
$[\C^{*}/G_{y}]$\footnote{One can check that $\mathfrak o_{e}$ can be also
  represented as the quotient stack $[\C^{*}/G_{y}(\rho^{-1},s)]$ or
  $[\C^{*}/G_{y}(\rho,r)]$ via $\rho^{-1}_{s}$ (resp. $\rho_{r}$).} then the
morphism $\rho:G_{y}\rightarrow \boldsymbol{\mu}_{t}$ induces that $\mathfrak
o_{e}$ is a $G_{e}-$gerbe over its coarse moduli $o_{e}\cong [\C^{*}
/\boldsymbol{\mu}_{t}]=\C^{*}$. Then the morphism $\rho^{-1}_{s}$ induces a
morphism from $\mathfrak U_{0}\cong [\mathbb C/G_{y}(\rho^{-1},s)]$
(resp.$\mathfrak U_{\infty}\cong [\mathbb C/G_{y}(\rho,r)]$ ) to $[\mathbb
C/\boldsymbol{\mu}_{ts}]$(resp. $[\C/\boldsymbol{\mu}_{tr}]$). The above
discussion fits into the following commutative digram
$$\xymatrix{
  \mathfrak U_{0} \ar[r] &[\mathbb C/\boldsymbol{\mu}_{ts}]\\
  \mathfrak o_{e}\ar[d]\ar[r]\ar[u] &o_{e}\cong[\mathbb
  C^{*}/\boldsymbol{\mu}_{ts}]\cong[\C^{*}/\boldsymbol{\mu_{t}}]\cong[\C^{*}/\boldsymbol{\mu}_{tr}]\ar[d]\ar[u]\\
  \mathfrak U_{\infty} \ar[r] &[\mathbb C/\boldsymbol{\mu}_{tr}]\ . }$$ Let
$\mathfrak p_{e}\cong \mathbb B G_{e}$ be a point of $\mathfrak o_{e}$. Define
$$H_{e}:=\pi_{1}(\mathfrak o_{e},\mathfrak p_{e})\ ,$$ then the coarsification map
$\mathfrak o_{e}\rightarrow o_{e}$ (with the corresponding base points
$\mathfrak p_{e}\rightarrow p_{e}$) induces a surjection
$$\phi_{e}:H_{e}=\pi_{1}(\mathfrak o_{e},\mathfrak p_{e})\rightarrow \pi_{1}(o_{e},p_{e})\cong \mathbb
Z \ ,$$ whose kernel is isomorphic to $G_{e}$ as shown in \cite[\S
6]{liu2020stacky}.

Choose a path from $\mathfrak p_{e}$ to $0$ and a path from $\mathfrak p_{e}$ to
$\infty$, the open embeddings $ \mathfrak o_{e}\hookrightarrow \mathfrak U_{0}$
and $ \mathfrak o_{e} \hookrightarrow \mathfrak U_{\infty}$ induce surjective
group homomorphisms

$$ H_{e}\!=\! \pi_{1}(\mathfrak o_{e},\mathfrak p_{e})\stackrel{\pi_{e,0}}{\longrightarrow}\pi_{1}(\mathfrak
U_{0},0)\cong G_{y}(\rho^{-1},s), \quad H_{e}\!=\! \pi_{1}(\mathfrak
o_{e},\mathfrak
p_{e})\stackrel{\pi_{e,\infty}}{\longrightarrow}\pi_{1}(\mathfrak
U_{\infty},\infty)\cong G_{y}(\rho,r)\ .$$

Applying $\pi_{1}$, the above discussion will be summarized into the following
commutative diagram with all rows exact:
\begin{equation}\label{eq:cd1}
  \xymatrix{
    1\ar[r] & G_{e} \ar[r] & \pi_{1}(\mathfrak U_{0})\cong
    G_{y}(\rho^{-1},s)\ar[r]^-{\rho^{-1}_{s}}
    &\boldsymbol{\mu}_{ts}\ar[r] &1\\
    1\ar[r] & G_{e} \ar[r]\ar@{=}[u]\ar@{=}[d] &\pi_{1}(\mathfrak o_{e}) \cong
    H_{e} \ar[r]^-{\phi_{e}}\ar[u]^{\pi_{e,0}}\ar[d]_{\pi_{e,\infty}}
    &\mathbb Z=\pi_{1}(o_{e})\ar[r]\ar[d]\ar[u] &1\\
    1\ar[r] & G_{e} \ar[r] & \pi_{1}(\mathfrak U_{\infty})\cong
    G_{y}(\rho,r)\ar[r]^-{\rho_{r}}
    &\boldsymbol{\mu}_{tr}\ar[r] &1\ ,\\
  }
\end{equation}
where in the third column $\mathbb Z\rightarrow \boldsymbol{\mu}_{tr}$ and
$\mathbb Z\rightarrow \boldsymbol{\mu}_{ts}$ are given by $d\mapsto
e^{\frac{d}{rt}}$ and $d\mapsto e^{\frac{d}{st}}$ respectively. By chasing
diagram, we can show that $H_{e}$ is isomorphic to the fiber product
$G_{y}(\rho^{-1},s) \!\times\!_{\boldsymbol{\mu}_{ts}}\mathbb Z$ (or
$G_{y}(\rho,r) \!\times\!_{\boldsymbol{\mu}_{tr}}\mathbb Z$ as well) using the
above diagram. Furthermore using the fact $G_{y}(\rho^{-1},s)$ isomorphic to the
fiber product $G_{y}\!\times\!_{\boldsymbol{\mu}_{t}}\boldsymbol{\mu}_{ts}$, we
have $H_{e}$ is also isomorphic to the fiber product $G_{y}
\!\times\!_{\boldsymbol{\mu}_{t}} \mathbb Z$.

Let $\P^{1}_{m,n}$ be the unique smooth DM stack with coarse moduli $\P^{1}$ and
trivial generic stablizer and only two special points $q_{\infty}\cong \mathbb
B\boldsymbol{\mu}_{m}$ and $q_{0}\cong \mathbb B\boldsymbol{\mu}_{n}$. Here we
can choose a $\C^{*}-$action on $\mathbb P^{1}_{m,n}$ such that $q_{0}$,
$q_{\infty}$ are the only $\C^{*}-$fixed points and the $\C^{*}-$weight of the
normal space $N_{q_{0}}$(resp. $N_{q_{\infty}}$) is $\frac{1}{m}$ (resp.
$\frac{-1}{n}$). Let $f:\mathbb P^{1}_{m,n}\rightarrow \mathfrak P$ be a
$\C^{*}-$fixed stable map of degree $(0,\frac{\delta}{r})$ which factors through
$\mathfrak P_{y}$. Assume that $f(q_{0})=0\in \mathfrak P_{y}$ and
$f(q_{\infty})=\infty \in \mathfrak P_{y}$. Moreover, if the
multiplicity\footnote{Here we can replace the index set $C$ for
  $\bar{I}_{\mu}\mathfrak P_{y}$
  by the set $\mathrm{Conj}(G_{y})$ of conjugacy classes of $G_{y}$.} at $q_{0}$
is $([g^{-1}],\mathrm{e}^{\frac{\delta}{s}},1)\in \mathrm{Conj}(G_{y})\!\times\!
\C^{*} \!\times\! \C^{*}$ such that $(x,[g])\in \bar{I}_{c}Y$ where $g\in
G_{y}$, then the multiplicity at $q_{\infty}$ is $([g],1,
\mathrm{e}^{\frac{\delta}{r}}) \in \mathrm{Conj}(G_{y})\!\times\! \C^{*}
\!\times\! \C^{*}$, and vice versa. Let $a$ be the order of $g$, when $r,s$ are
sufficiently large primes, as $f$ is representable, we have $m=ar,n=as$. Assume
that $age_{g}(L|_{y} )=\frac{h}{a}$ for some integer $h$ with $0\leq h<a$. As
$\delta=r\cdot deg (f^{*}\mathcal O(\mathcal D_{\infty}))$, using orbifold
Riemann-Roch for the $f^{*}\mathcal O(r \mathcal D_{\infty})$, the condition
$$ \delta -\frac{h}{a}\in \mathbb Z\ ,$$
is a necessary condition to ensure the map $f$ to exist. Then $a\delta$ must be
an integer.

Let $T_{ar,as}$ be the subgroup of $(\C^{*})^{2}$ defined by the equation
$t_{1}^{as}=t_{2}^{ar}$. Let $F:U\rightarrow \C^{*} \!\times\! U$ be the morphism sending $(x,y)$ to
$(1,x^{a\delta},y^{a\delta})$. Then $F$ is equivariant with respect to the group
homomorphism from $T_{ar,as}$ to $G_{y}\times (\C^{*})^{2}$:
$$(t_{1},t_{2})\mapsto
(\tau(t_{1}^{-s}t^{r}_{2}),t_{1}^{a\delta},t_{2}^{a\delta}) \ ,$$ where
$\tau:\boldsymbol{\mu}_{a}\rightarrow G_{y}$ is the group morphism sending
$\mu_{a}$ to $g$. We can check $F$ descents to be a morphism $\widetilde{F}$
from $\mathbb P^{1}_{ar,as}$ to $\mathfrak P_{y}$ of degree
$(0,\frac{\delta}{r})$\footnote{We use the degree notation by viewing
  $\widetilde{F}$ as a morphism to $\mathfrak P$.} with the multiplicity
decoration $([g], 1,\mathrm{e}^{\frac{\delta}{r}})$ at $q_{\infty}$. Moreover,
if we change $g$ to any conjugate of $g$ in $G_{y}$, the descended morphisms are
all isomorphic to each other. Conversely, we will show the following.
\begin{lemma}\label{thm:classify-mor}
  With the notations as above, let $f$ be a $\C^{*}-$fixed morphism of degree
  $(0,\frac{\delta}{r})$ to $\mathfrak P$ factor through $\mathfrak P_{y}$ such
  that the multiplicity is given by $([g],1, \mathrm{e}^{\frac{\delta}{s}})$ at
  $q_{\infty}$. Then the morphism $f$ must be isomorphic to $\widetilde {F}$ (up
  to a unique $2-$isomorphism) constructed as above.
\end{lemma}
\begin{proof}
  Let $f_{1}$, $f_{2}$ be two stable maps with same associated decorated graph
  as above, if we assume they are isomorphic when restricting to $O_{e}: =
  \P^{1}_{ar,as}\backslash \{q_{0},q_{\infty}\}$, then $f_{1}$ is isomorphic to
  $f_{2}$ up to a unique $2-$isomorphism(see e.g.,\cite[Lemma 2.5]{Cheong_2015}
  for a proof). Thus we only need to show the restriction $f|_{\mathcal O_{e}}$
  is uniquely determined by the degree data and multiplicity data at
  $q_{\infty}$. Note that $f|_{O_{e}}:O_{e}\rightarrow \mathfrak o_{e}$ is an
  \'etale covering map. Indeed, as $O_{e}$ and the coarse moduli $o_{e}$ of
  $\mathfrak o_{e}$ are both isomorphic to $\C^{*}$ and the coarsing map from
  $\mathfrak o_{e} $ to $o_{e}$ is \'etale, then the \'etaleness of $f|_{O_{e}}$
  follows from the \'etaleness of the composition
  $$O_{e}\rightarrow \mathfrak o_{e}\rightarrow o_{e} $$
  as any $\C^{*}-$fixed morphism from $\C^{*}$ to $\C^{*}$ must be a cyclic
  covering. By \cite[Theorem 18.19]{noohi05_found_topol_stack_i}, we have a
  bijection between conjugacy classes of subgroup
  $(f|_{O_{e}})_{*}(\pi_{1}(O_{e},P_{e}))$ ($P_{e}$ is the inverse image of the
  point $\mathfrak p_{e}$ under $O_{e}\rightarrow \mathfrak o_{e}$) of
  $\pi_{1}(\mathfrak o_{e},\mathfrak p_{e})$ and isomorphism classes of covering
  map of $\mathfrak o_{e}$ (not necessarily base point preserving), it's enough
  to show the conjugacy class of subgroup
  $(f|_{O_{e}})_{*}(\pi_{1}(O_{e},P_{e}))$ is uniquely determined by the degree
  data and the monodromy data at $q_{\infty}$.

  We will drop the base points in the notation of fundmental group in the sequel
  if it's clear in the context. Let
  $$U_{0}:= \P^{1}_{ar,as} \backslash \{q_{\infty}\}\cong [\C/\boldsymbol{\mu}_{as}],\quad
  \text{and}\quad U_{\infty}:= \P^{1}_{ar,as} \backslash \{q_{0}\}\cong
  [\C/\boldsymbol{\mu}_{ar} ]\ . $$Choose a path from $P_{e}$ to $q_{0}$ in
  $U_{0}$ and a path from $P_{e}$ to $q_{\infty}$ in $U_{\infty}$, then the open
  embeddings $ O_{e}\hookrightarrow U_{0}$ and $ O_{e} \hookrightarrow
  U_{\infty}$ induces surjective group homomorphisms
$$\xymatrix{
  & \mathbb Z=\pi_{1}(O_{e})\ar[r]^-{\pi_{q_{0}}}&\pi_{1}(U_{0})\cong
  \boldsymbol{\mu}_{as} &\mathbb Z=
  \pi_{1}(O_{e})\ar[r]^-{\pi_{q_{\infty}}}&\pi_{1}(U_{\infty})\cong\boldsymbol{\mu_{ar}}\
  . }$$ Here we choose the identification between $\pi_{1}(O_{e})$ and $\mathbb
Z$ such that $\pi_{q_{0}}(1)=e^{\frac{-1}{as}}$ and
$\pi_{q_{\infty}}(1)=e^{\frac{1}{ar}}$.\footnote{Different choices of path from
  $P_{e}$ to $q_{0}$(resp. $q_{\infty}$) will only affect $\pi_{q_{0}}$(resp.
  $\pi_{q_{\infty}}$) up to conjugacy, then the image $\pi_{q_{0}}(1)$ (resp.
  $\pi_{q_{\infty}}(1)$) doesn't depend on the choice of the path.} Then the
morphism
 $$\xymatrix{
   \pi_{1}(O_{e})\ar[r]
   &\pi_{1}(U_{\infty})\ar[r]^{(f|_{U_{\infty}})_{*}}&\pi_{*}(\mathfrak
   U_{\infty}) }$$ sends the generator $1$ of $\pi_{1}(O_{e})$ to the conjugacy
 class of $(g,e^{\frac{\delta}{r}})\in G_{y}(\rho,r)\cong \pi_{1}(\mathfrak
 U_{\infty})$ by the multiplicity data at $q_{\infty}$. We have a commutative
 digram with all vertical lines are inclusions:
 $$\xymatrix{
   &O_{e}\ar[r]^{f|_{O_{e}}}\ar[d]&\mathfrak o_{e}\ar[r]\ar[d] &o_{e}\ar[d]\\
   &U_{\infty} \ar[r] &\mathfrak U_{\infty}\ar[r] & [\mathbb
   C/\boldsymbol{\mu}_{tr}]\ . }$$ Using \ref{eq:cd1}, the above diagram induces
 the following commutative diagram
 \begin{equation}\label{diagram-dec}
   \xymatrix{
     &\mathbb Z=\pi_{1}(O_{e}) \ar[r]^{(f|_{O_{e}})_{*}}\ar[d] &\pi_{1}(\mathfrak o_{e}) \cong
     H_{e} \ar@{->>}[r]^-{\phi_{e}}\ar[d]^{\pi_{e,\infty}}
     &\pi_{1}(o_{e})\cong\mathbb Z\ar[d]\\
     &\boldsymbol{\mu}_{ar}=\pi_{1}(U_{\infty}) \ar[r] & \pi_{1}(\mathfrak
     U_{\infty})\cong G_{y}(\rho,r)\ar@{->>}[r]^-{\rho_{r}}
     &\boldsymbol{\mu}_{tr}\ , }
 \end{equation}
 where in the third column $\mathbb Z\rightarrow \boldsymbol{\mu}_{tr}$ is given
 by $d\mapsto \mathrm{e}^{\frac{d}{tr}}$. Then we see that the composition
 $$\xymatrix{
   \pi_{1}(O_{e})\cong \mathbb Z\ar[r]^{(f|_{O_{e}})_{*}} &\pi_{1}(\mathfrak
   o_{e})=H_{e}\ar[r]^{\phi_{e}}&\pi_{1}(o_{e})\cong \mathbb Z\ . }$$ from
 $\mathbb Z$ to $\mathbb Z$ is multiplication by $t\delta$ as the generator $1$
 of $\pi_{1}(O_{e})$ sends to $e^{\frac{\delta}{r}}$ in the group
 $\boldsymbol{\mu}_{tr}$ in whatever path from the up-left corner to the
 bottom-right corner in the above commutative diagram and $r$ is sufficiently
 large.

 Let $\gamma_{e}\in H_{e}$ be the image of generator $1\in \pi_{1}(O_{e})$ in
 $H_{e}$. By the previous discussion, $\gamma_{e}$ maps to the conjugacy
 class of $(g,e^{\frac{\delta}{r}})$ via $\pi_{e,\infty}$, and maps to $t\delta$
 in $\pi_{1}(o_{e})$. Conversely, using the fact $H_{e}$ is isomorphic to fiber
 product $G_{y}(\rho,r)\times_{\boldsymbol{\mu}_{tr}}\mathbb Z$ using the right
 square of \ref{diagram-dec} , we see that the conjugacy class of $\gamma_{e}$,
 is uniquely determinded by the conjugacy class of the image of $\gamma_{e}$ in
 $\pi_{1}(\mathfrak U_{\infty})$ and the image of $\gamma_{e}$ in
 $\pi_{1}(o_{e})$, which are determined by the decorated data. Hence the
 conjugacy class of the subgroup $(f|_{O_{e}})_{*}(\pi_{1}(O_{e}))$ is also
 uniquely determined by the decorated data. This completes the proof.
\end{proof}

Next we describe the automorphism group $Aut(f)$ of $f$.
\begin{proposition}\label{prop:aut-order}
  Use the notation in the above proof, we have $Aut(f)\cong
  C_{H_{e}}(\gamma_{e})/\langle{\gamma_{e}}\rangle$, and the order of
  automorphism group is equal to
  $|Aut(f)|=|C_{H_{e}}(\gamma_{e})/\langle{\gamma_{e}}
  \rangle|=\delta|C_{G_{y}}(g)|$. Here $C_{H_{e}}(\gamma)$ (resp.
  $C_{G_{y}}(g)$) is the centralizer of $\gamma_{e}$(resp. $g$) in the group
  $H_{e}$(resp. $G_{y}$).
\end{proposition}
\begin{proof}
  The automorphism group of $f$ is isomorphic to the automorphism group of the
  restriction $f|_{O_{e}}$, then the first claim follows from standard fact
  about covering deck transformation(not necessarily basepoint preserving) as
  $\gamma_{e}$ generates the subgroup $(f|_{O_{e}})_{*}(\pi_{1}(O_{e}))$ in
  $\pi_{1}(\mathfrak o_{e})$.
  
  Assume that the image $\phi_{e}(C_{H_{e}}(\gamma_{e}))$ in $\mathbb Z$ is
  generated by the positive integer $u$, as $H_{e}$ is isomorphic to the fiber
  product $G_{y}\times_{\boldsymbol{\mu}_{t}}\mathbb Z$.\footnote{We have already
    seen that $H_{e}$ is isomorphic to
    $G_{y}(\rho^{-1},s)\times_{\boldsymbol{\mu}_{st}}\mathbb Z$, then we use the
    isomorphism $G_{y}(\rho^{-1},s)\cong
    G_{y}\times_{\boldsymbol{\mu}_{t}}\boldsymbol{\mu}_{st}$ } More precisely,
 $$ H_{e}\cong \{(g,d)\in G_{y} \!\times\! \mathbb Z|
 \rho(g)=\mathrm{e}^{\frac{d}{t}}\}\ .$$ Then $\gamma_{e}=(g,t\delta)\in H_{e}$.
 We can show that $t$ and $t\delta$ are both divided by $u$. Thus the image of
 $C_{G_{y}}(g)$ under $\rho$ is generated by $\mathrm{e}^{\frac{u}{t}}$.

 The group $C_{H_{e}}(\gamma_{e})$ can be represented as
  $$ \{(h,d)\in C_{G_{y}}(g) \!\times\! \mathbb Z| \rho(h)=\mathrm{e}^{\frac{d}{t}}\} $$
  in $H_{e}$. Then the claim about $|Aut(f)|$ follows from the following two
  exact sequences:
$$\xymatrix{
  1\ar[r] &G_{e}\cap C_{G_{y}}(g)\ar[r]^{\theta} &C_{H_{e}}(\gamma_{e})/\langle
  {\gamma_{e}} \rangle\ar[r] & u \mathbb Z / t\delta\mathbb Z \ar[r] &1\ , }$$
and
 $$\xymatrix{
   1\ar[r] &G_{e}\cap C_{G_{y}}(g)\ar[r] &C_{G_{y}}(g)
   \ar[r]^-{\rho}&\boldsymbol{\mu}_{\frac{t}{u}}\ar[r]&1\ . }$$ Here $u \mathbb
 Z$(resp. $t\delta \mathbb Z$) is the subgroup of $\mathbb Z$ generated by the
 integer $u$(resp. $t\delta$), the morphism $\theta$ is the induced from the
 restriction to $G_{e}\cap C_{G_{y}}(g)$ of the inclusion from $G_{e}$ to
 $H_{e}$.
\end{proof}

We can make the above proposition more concretely. Let $(h,d)\in
C_{H_{e}}(\gamma_{e})$ be an element, we associate $(h,d)$ an element in the
automorphism group $Aut(f)$ as follows: define
$$\theta_{h,d}:U\rightarrow U:\quad (x,y)\mapsto (x,e^{\frac{d}{art\delta}}y)
\ ,$$ which descents to be an isomorphism of $\P^{1}_{ar,as}$ which we still
denote to be $\theta_{h,d}$. Recall that $F$ is a morphism from $U$
$\rightarrow \C^{*}\times U$ defined as in Lemma \ref{thm:classify-mor}, then
$F\circ \theta_{h,d}$ differ from $F$ by the action of
the element $(h,1,e^{\frac{d}{rt}})$ on the target via the action of
$G_{y}\times (\C^{*})^{2}$ on $\C^{*}\times U$ which also commutes the action of
$T_{ar,as}$ and $G_{y}\times (\C^{*} )^{2}$ on the source and target
respectively. Then by descent, this defines an $2-$isomorphism
$\alpha_{h,d}:F\rightarrow F\circ \theta_{h,d}$. Thus
the pair $(\theta_{h,d},\alpha_{h,d})$ defines an automorphism in $Aut(f)$.

For two elements $(h_{1},d_{1})$ and $(h_{2},d_{2})$ in $C_{H_{2}}(\gamma_{e})$,
the pairs $(\theta_{h_{1},d_{1}},\alpha_{h_{1},d_{1}})$ and
$(\theta_{h_{2},d_{2}},\alpha_{h_{2},d_{2}})$ are isomorphic to each other if
there exists some integer $n$ such that $h_{2}=h_{1}g^{n}$ and
$d_{2}=d_{1}+nt\delta$. Indeed, let $m_{(1,e^{\frac{nt\delta}{art\delta}})}$ be
the automorphism of $U$ acted upon by $(1,e^{\frac{nt\delta}{art\delta}})$, we
have $\theta_{h_{2},d_{2}}= m_{ (1,e^{\frac{nt\delta}{art\delta}})} \circ
\theta_{h_{1},d_{1}}$, which defines a two-isomorphism
$\beta_{n}:\theta_{h_{1},d_{1}} \rightarrow \theta_{h_{2},d_{2}}$ in
$\mathrm{Hom}(\P^{1}_{ar,as},\P^{1}_{ar,as})$, then the composition $f\circ
\beta_{n}$ defines a two-isomorphism from $F\circ
\theta_{h_{1},d_{1}}$ to $F\circ \theta_{h_{2},d_{2}}$ such that
$\alpha_{h_{2},d_{2}}=f\circ \beta_{n}\circ \alpha_{h_{1},d_{1}}$.

The morphism $f:\P^{1}_{ar,as}\rightarrow \mathfrak P_{y}$ defines a $\mathbb
C-$point of $\mathcal K_{0,2}(\P Y_{r,s},(0,\frac{\delta}{r}))$, where we denote
by $q_{0}$(resp. $q_{\infty}$) by the first (resp. second) marking, then the
evaluation map $$ ev_{q_{\infty}}: \mathcal K_{0,2}(\P
Y_{r,s},(0,\frac{\delta}{r}))\rightarrow \bar{I}_{\mu}\P Y_{r,s}\ ,$$ induces a
morphism from the isotropy group $Aut(f)$ to the isotropy group $Aut(
(y,([g],1,e^{\frac{\delta}{r}})))$, which, from the concrete description of
$Aut(f)$ above, coincides the morphism from $C_{H_{e}}(\gamma_{e})/\langle
{\gamma_{e}} \rangle$ to $C_{G_{y}}(g)/\langle{g}\rangle$ induced from the
projection of $G_{y}\times \mathbb Z$ to the first factor. From this concrete
description, we also know the kernel of the morphism of isotropy groups:

\begin{proposition}\label{prop:ker-aut}
  The kernel of the morphism $Aut(f)\rightarrow Aut\big(
  (y,([g],1,e^{\frac{\delta}{r}}))\big)$ is isomorphic to $\mathbb Z/a\delta
  \mathbb Z$.
\end{proposition}
\begin{proof}
  This follows from by applying the snake lemma to the first two exact rows of
  the following diagram:
$$\xymatrix{
  1\ar[r] &\mathbb Z \cong \langle{(1,at\delta)}\rangle \ar@{^{(}->}[d]\ar[r]
  &\mathbb Z\cong \langle{\gamma_{e}} \rangle
  \ar[r]^{pr_{1}}\ar@{^{(}->}[d] &\mathbb Z/a \mathbb Z \cong \langle{g} \rangle \ar[r]\ar@{^{(}->}[d] &1\\
  1\ar[r] &\mathbb Z\cong \langle{(1,t)} \rangle \ar[r]\ar@{->>}[d]
  &C_{H_{e}}(\gamma_{e})
  \ar[r]^{pr_{1}}\ar@{->>}[d] & C_{G_{y}}(g) \ar[r]\ar@{->>}[d] &1\\
  1\ar[r] & \mathbb Z/a\delta \mathbb Z \ar[r] &C_{H_{e}}(\gamma_{e})/\langle
  {\gamma_{e}} \rangle\ar[r]^{pr_{1}}\ar@{=}[d] & C_{G_{y}}(g)/\langle{g}\rangle
  \ar[r]\ar@{=}[d]
  &1\\
  & &Aut(f) &Aut\big( (y,([g],1,e^{\frac{\delta}{r}}))\big) \ , }$$ where
$pr_{1}$ are induced from the projection from $G_{y}\times \mathbb Z$ to the
first factor.
\end{proof}

\subsection{Construction of family stable map}\label{sec:app-edge}
Let $e$ be an edge with decorated degree $\frac{\delta}{r}$ and decorated
multiplicity $(c,1,\mathrm{e}^{\frac{\delta}{r}})$ at the half-edge $h_{\infty}$
as in \S \ref{sec:lc1}. Let $a:=a(c)$ be the integer associated to $c$ as in the
Assumption \ref{assump:inertiaindex}. Define the space $\mathcal M_{e}$ to be
the root gerbe $\sqrt[as\delta]{ L^{\vee}/I_{c}Y }$. We will generalize the
construction in Lemma \ref{thm:classify-mor} to a family version; we will
construct a family of $\C^{*}-$fixed stable map over $\mathcal M_{e}$ to $\P
Y_{r,s}$ with the associated decorated degree. Then we get a
morphism $$g:\mathcal M_{e}\rightarrow \mathcal K_{0,2}(\P
Y_{r,s},(0,\frac{\delta}{r})) \ .$$ We will show the image of $g$ is a closed
and open substack of the $\C^{*}-$fixed part of $\mathcal K_{0,2}(\P
Y_{r,s},(0,\frac{\delta}{r}))$, and $g$ is a finite \'etale map into the image
of $g$. Then we can use $\mathcal M_{e}$ to do an explicit computation of the
edge contribution in the virtual localization formula.

By an abuse of notation, denote $L^{\vee}$ to be the line bundle over $\mathcal
M_{e}$ via a pull-back of the line bundle $L^{\vee}$ over $Y$ along the natural
map $\mathcal M_{e}\rightarrow Y$, which is a composition of maps $\mathcal
M_{e}\rightarrow \bar{I}_{c}Y \rightarrow Y$. Denote by $(L^{\vee})^{0}$ the
open set of $L^{\vee}$ with zero section removed.

Let $U:=\mathbb C^{2}\backslash \{0\}$. Let $\mathcal C_{e}$ be the quotient
stack $[(L^{\vee})^{0}\times U / \C^{*}\times T_{ar,as}]$ defined by the (right)
action
$$(l,x,y)(t,t_{1},t_{2})=(t^{-as\delta}l, tt_{1}x,t_{2}y)\ , $$
where $(l,x,y)\in (L^{\vee})^{0}\times U $ and $(t,t_{1},t_{2})\in \C^{*} \times
T_{ar,as}$. Here $T_{ar,as}$ is the subgroup of $(\C^{*})^{2}$ defined
by$\{(t_{1},t_{2})\in (\C^{*})^{2}|\; t_{1}^{as}=t_{2}^{ar} \}$. Then $\mathcal C_{e}$ is a family of orbi-$\P^{1}$ bundle over $\mathcal M_{e}$
such that all fibers are isomorphic to $\P^{1}_{ar,as}$.

\begin{remark}\label{rmk:psi-edge}
  Let $\bar{\psi}_{\star}$ be the psi-class associated with the gerbe marking
  corresponding to the zero loci of $y$ of $\mathcal C_{e}$. We see that
  $\bar{\psi}_{\star}$ is equal to $\frac{\lambda-c_1(L)}{\delta}$.
\end{remark}

We will define a family map $f$ using the categorical description of $\P
Y_{r,s}$ and $C_{e}$:
\begin{equation*}
  \xymatrix{
    \mathcal C_e  \ar^-{f}[r]\ar_{\pi}[d] & \P Y_{r,s} \\
    \mathcal M_{e}:=\sqrt[as\delta]{ L^{\vee}/I_{c}Y }.
  }
\end{equation*}

\subsubsection{A categorical description about stacks $\P Y_{r,s}$ and $\mathcal C_{e}$}
The space $\P Y_{r,s}$ can be represented as the quotient $[(L^{\vee})^{0}\times
U / (\C^{*})^{2}]$ defined by the action
$$(l,x,y)(t_{1},t_{2})=(t_{1}^{-s}t_{2}^{r}l, t_{1}x,t_{2}y) \ .$$
We see that, for any $\mathbb C-$scheme $S$, an object in $ \P
Y_{r,s}(S)$ is given by a tuple $(w,M_{1},M_{2},\rho, s_{1},s_{2})$
(c.f.~\cite[\S 10.2.7,\S 10.2.8]{MR3495343}), where
\begin{enumerate}
\item $w:S\rightarrow Y$ is an object\footnote{Here we view the morphism
    $S\rightarrow Y$ as an equivalent definition of objects in $Y(S)$ using the
    2-Yoneda Lemma.} in $Y(S)$,
\item $M_{1},M_{2}$ are line bundles over $S$,
\item $\rho: M_{1}^{\otimes s} \otimes M_{2}^{-\otimes r}\cong w^{*}L^{\vee} $
  is an isomorphism,
\item $ s_{i}\in \Gamma(S,M_{i}), i\in\{1,2\}$ are sections such that $
  s_{1},s_{2}$ don't vanish simultaneously.
\end{enumerate}

Let $(w' ,M' _{1},M' _{2},\rho' , s' _{1},s' _{2})$ be another object in $\P
Y_{r,s}(S')$, given a morphism $g:S'\rightarrow S$, then a morphism from $(w'
,M' _{1},M' _{2},\rho' , s' _{1},s' _{2})$ to $(w,M_{1},M_{2},\rho,
s_{1},s_{2})$ over $g$ is a tuple $(g^b,\alpha_{1},\alpha_{2})$. Here $ g^{b}:
w'\rightarrow w$ is a morphism over $g$, $\alpha_{i}: g^{*}M_{i}\rightarrow
M'_{i} $ is an isomorphism which maps $g^{*}s_{i}$ to $s_{i}'$,. Furthermore, we
require a commutative square
$$\xymatrix{
  &g^{*}M_{1}^{\otimes s}\otimes g^{*}M_{2}^{-\otimes r}\ar[d]_{g^{*}\rho}
  \ar[rr]^{\alpha_{1}^{\otimes s}\otimes \alpha^{-\otimes r}_{2}} &
  &M_{1}'^{\otimes s}\otimes M_{2}'^{-\otimes r}\ar[d]^{\rho'}\\
  & g^{*}w^{*}L^{\vee} \ar[rr]^{g^b} & & (w')^{*}L^{\vee}\ , }$$ and
$\alpha_{i}$ should satisfy obvious commutative laws for composition of two maps
(c.f.~\cite[\S 9.1.10(i)]{MR3495343}).

\begin{remark}
  Let $w:S\rightarrow Y$ be any object in $Y(S)$, denote by $\mathbb
  P_{r,s}(w^{*}L^{\vee}\oplus \mathbb C )$ the stack quotient
  $$ [(w^{*}L^{\vee})^{0}\times U /
  (\C^{*})^{2}]\ ,$$ where the action is defined in the same way as $\P
  Y_{r,s}$. One can verify that we have a Cartesian diagram:
   $$\xymatrix{
     \mathbb P_{r,s}(w^{*}L^{\vee}\oplus \mathbb C )\ar[r]\ar[d] &\P
     Y_{r,s}\ar[d]^{\mathrm{pr_{r,s}}}\\
     S\ar[r]^{w} &Y\ . }$$ This property can also lead to our categorical
   description of $\P Y_{r,s}$, see a closed-related discussion in~\cite[\S
   10.2.8]{MR3495343}.
 \end{remark}




 We also give a categorical description of $\mathcal C_{e}$. Recall that
 $\mathcal C_{e}$ can be represented as the quotient stack $[(L^{\vee})^{0}
 \times U / ( T_{ar,as}\!\times\! \C^{*})]$ given by the action
$$(l,x,y)(t_{1},t_{2},t_{3})=(t_{3}^{-as\delta}l, t_{3}t_{1}x,t_{2}y) $$
where $(l,x,y)\in (L^{\vee})^{0}\times U $ and $(t_{1},t_{2},t_{3})\in T_{ar,as}
\!\times\! \C^{*}$. For any scheme $S$ over $\mathrm{Spec}(\mathbb C)$, an
object in $ \mathcal C_{e}(S)$ is given by a tuple $\tau:=(w,\theta,
M_{1},M_{2},M_{3},\rho_{1},\rho_{2}, s_{1},s_{2})$, where
\begin{enumerate}
\item $w$ is a object in $Y(S)$, $\theta$ is an automorphism in $Aut_{S}(w)$
  such that $(w,\theta)\in I_{c}Y$,
\item $M_{1},M_{2},M_{3}$ are line bundles over $S$,
\item $\rho_{1}:M_{3}^{\otimes as\delta}\cong w^{*}L^{\vee}$ and
  $\rho_{2}:M_{1}^{\otimes as} \otimes M_{2}^{-\otimes ar}\cong \mathcal O_{S}$
  are isomorphisms,
\item $ s_{1}\in \Gamma(S,M_{1}\otimes M_{3})$ and $s_{2}\in \Gamma(S,M_{2})$
  are sections such that $ s_{1},s_{2}$ don't vanish simultaneously.
\end{enumerate}

Given another $\mathbb C-$scheme $S'$, let $\tau'=(w' ,\theta',M' _{1},M'
_{2},M_{3}', \rho'_{1},\rho'_{2} , s' _{1},s' _{2})$ be an object in $\mathcal
C_{e}(S')$, given a morphism $g:S'\rightarrow S$, then a morphism from $(w'
,\theta',M' _{1},M' _{2}, M_{3}',\rho'_{1},\rho_{2} , s' _{1},s' _{2})$ to
$(w,\theta, M_{1},M_{2},M_{3},\rho_{1},\rho_{2}, s_{1},s_{2})$ over $g$ is a
tuple $(g^b,\alpha_{1},\alpha_{2},\alpha_{3})$, where $g^{b}:
(w',\theta')\rightarrow (w,\theta)$ is a morphism in $I_{c}Y$ over the morphism
$g$, $\alpha_{i}:g^{*}M_{i}\rightarrow M_{i}'$ are isomorphisms such that
$\alpha_{1}\otimes \alpha_{3}: g^{*}M_{1}\otimes g^{*}M_{3}\rightarrow
M'_{1}\otimes M_{3}' $ maps $g^{*}s_{1}$ to $s_{1}'$, and $\alpha_{2}:
g^{*}M_{2}\cong M'_{2} $ maps $g^{*}s_{2}$ to $s_{2}'$. Furthermore, these data
are require to satisfy the following conditions:
\begin{enumerate}
\item We have two commutative squares
$$\xymatrix{
  &g^{*}M_{3}^{\otimes as\delta}\ar[d]_{g^{*}\rho_{1}}
  \ar[rr]^{\alpha_{3}^{\otimes as\delta}} & &M_{3}'^{\otimes as\delta}\ar[d]^{\rho'_{1}}\\
  & g^{*}w^{*}L^{\vee}  \ar[rr]^{g^b} & & (w')^{*}L^{\vee}\\
}$$ and
 $$\xymatrix{
   &g^{*}M_{1}^{\otimes as}\otimes g^{*}M_{2}^{-\otimes
     ar}\ar[d]_{g^{*}\rho_{2}}\ar[rr]^{\alpha_{1}^{\otimes as}\otimes
     \alpha^{-\otimes ar}_{2}} & &M_{1}'^{\otimes as}\otimes M_{2}'^{-\otimes
     ar}
   \ar[d]^{\rho'_{2}}\\
   & g^{*}\mathcal O_{S}\ar[rr]^{\cong}& & \mathcal O_{S'} \ . }$$
\item $\alpha_{1},\alpha_{2}$ should satisfy commutative laws for composition of
  maps.
\end{enumerate}

For each object $\tau:=(w,\theta, M_{1},M_{2},M_{3},\rho_{1},\rho_{2},
s_{1},s_{2}) \in \mathcal C_{e}(S)$, we give the following definition:
\begin{definition}\label{def:induced-mor}
  The collection of line bundles $M_{1}^{-\otimes s}\otimes M_{2}^{\otimes r}$
  and morphisms $\alpha_{1}^{\otimes s}\otimes \alpha^{-\otimes r}_{2}$ defines
  a line bundle on $\mathcal C_{e}$, which we denote to be $\mathcal
  L$.\footnote{This line bundle is also isomorphic to $L_{-s\chi_{1}}\otimes
    L_{r\chi_{2}}$ defined in \S \ref{subsec:computeedge}.} Using $\rho_{2}$, we
  see that the line bundle $M_{1}^{-\otimes s}\otimes M_{2}^{\otimes r}$
  represents a torsion element in $\mathrm{Pic}(S)$ of which the order is
  divided by $a$. This defines a morphism $\alpha:\mathcal C_{e}\rightarrow
  \mathbb B\boldsymbol{\mu}_{a}$ satisfying that: let $\mathbb L=[\mathbb
  C/\boldsymbol{\mu}_{a}]$ be the line bundle over $\mathbb B
  \boldsymbol{\mu}_{a}$ associated with the standard representation of
  $\boldsymbol{\mu}_{a}$, we have $\alpha^{*}\mathbb L\cong \mathcal L$. Then we
  can naturally associate each object $\tau:S\rightarrow \mathcal C_{e}$ a
  morphism $\bar{f}_{1}:=\alpha\circ \tau$ and an isomorphism $\alpha_{\tau}:
  M_{1}^{-\otimes s}\otimes M_{2}^{\otimes r} \cong \bar{f}_{1}^{*}\mathbb L$.
  Given another object $\tau'\in \mathcal C_{e}(S')$ and a morphism
  $u:=(g^b,\alpha_{1},\alpha_{2},\alpha_{3})$ in $Hom_{\mathcal
    C_{e}}(\tau',\tau)$ over $g:S'\rightarrow S$, we have

\begin{equation}\label{alpha}
  \xymatrixcolsep{4pc}\xymatrix{
    &g^{*}(M_{1}^{\otimes s}\otimes
    M_{2}^{\otimes -r})\ar[d]^{g^{*}\alpha_{\tau}}  \ar[r]^{\alpha_{1}^{-\otimes s}\otimes \alpha_{2}^{\otimes r}
    } &M_{1}^{'-\otimes s}\otimes
    M_{2}^{'\otimes r}\ar[d]^{\alpha_{\tau'}} \\
    &g^{*}\bar{f}_{1}^{*}\mathbb L \ar[r]^{\alpha(u)} &(\bar{f}'_{1})^{*}\mathbb L\ .
  }
\end{equation}

Note $\bar{f}_{1}$ induces a $S-$morphism $f_{1}:S\rightarrow S \times \mathbb B
\boldsymbol{\mu}_{a}$ such that $\mathrm{pr}_{2}\circ f_{1}=\bar{f}_{1}$, where
$\mathrm{pr}_{2}$ is the projection from $S\times \mathbb B
\boldsymbol{\mu}_{a}$ to the second factor.

Recall that there is an isomorphism of the stack $\mathrm{RepHom}(\mathbb B
\mathbf{\mu}_{a},Y)$ with the $I_{\mu_{a}}Y$ from ~\cite[\S
3.2]{Dan_Abramovich_2008}, the data $(w,\theta)$ naturally defines a morphism
$f_{2}:S\times \mathbb B \boldsymbol{\mu}_{a}\rightarrow Y $ such that the
generator of the isotropy group $Aut_{S}(i_{S})= (\boldsymbol{\mu_{a}})_{S}$
sends to the automorphism $\theta\in Aut_{S}(w)$, where $i_{S}$ is a
$S-$morphism $i_{S}:S\rightarrow S\times \mathbb B \boldsymbol{\mu}_{a}$
corresponding to the trivial $\boldsymbol{\mu}_{a}$ principal bundle over $S$.
Moreover we have $w=f_{2}\circ i_{S}$.
\end{definition}

Let $a$ be the integer associated to $c$ as in \ref{assump:inertiaindex}. Write
$age_{c}(L)=\frac{h}{a}$ for an integer $h$ with $0\leq h<a$. Let $\mathbb
L_{a,-h}$ be the line bundle over $\mathbb B \boldsymbol{\mu}_{a}$ associated to
the character of $\boldsymbol{\mu}_{a}$ by sending $\mu_{a}^{i}$ to
$\mu_{a}^{-ih}$. Then $\mathbb L_{a,-h}=\mathbb L^{-\otimes h}$ is defined in
\ref{def:induced-mor}.

\begin{proposition}\label{prop:induced-mor}
  We have a natural isomorphism $\epsilon_{f_{2}}:w^{*}L^{\vee}\boxtimes \mathbb
  L_{a,-h}\rightarrow f_{2}^{*}L$, and we have a natural isomorphism
  $\bar{\rho}_{\tau}:w^{*}L^{\vee}\otimes (M_{1}^{\otimes sh}\otimes
  M_{2}^{\otimes -rh})\rightarrow (f_{2}\circ f_{1})^{*}L^{\vee}$.
\end{proposition}
\begin{proof}
  We will prove $\epsilon_{f_{2}}$ by \'etale descent. First note $i_{S}:S
  \rightarrow S\times \mathbb B \mu_{a}$ is an \'etale covering with Galois
  group $\boldsymbol{\mu}_{a}$, then $f_{2}$ can be recovered form the descent
  data given by the pair $(w,\theta)$(c.f.,~\cite{MR2223406}). Then the
  $f^{*}_{2}L^{\vee}$ can be recovered descent data given by the pair
  $(w^{*}L^{\vee}, \widetilde {\theta}:w^{*}L^{\vee}\rightarrow w^{*}L^{\vee})$,
  where $\widetilde{\theta}$ is an automorphism of $w^{*}L^{\vee}$ induced by
  $\theta$. As $age_{c}(L)=\frac{h}{a}$, $\widetilde{\theta}$ is nothing else
  but the scalar multiplication by $\mathrm{e}^{\frac{-h}{a}}$. By descent, we
  see that the line bundle $f_{1}^{*}L^{\vee}$ is isomorphic to
  $w^{*}L^{\vee}\boxtimes \mathbb L_{a,-h}$.

  Finally the claim about $(f_{2}\circ f_{1})^{*}L^{\vee}$ follows the next
  lemma.
\end{proof}

\begin{lemma}\label{lem:induced-mor}
  Let $W,U_{1},U_{2}$ be stacks and $L_{1}$ and $L_{2}$ be two line bundles over
  $U_{1}$ and $U_{2}$ respectively. Let $t_{1}:W\rightarrow U_{1}$ and
  $t_{2}:W\rightarrow U_{2}$ be two morphisms and $t_{1}\times
  t_{2}:W\rightarrow U_{1}\times U_{2}$ is the induced morphism using the fiber
  product. Then we have that a natural isomorphism $(t_{1}\times
  t_{2})^{*}(L_{1}\boxtimes L_{2})\cong t_{1}^{*}L_{1}\otimes t_{2}^{*}L_{2}$.
\end{lemma}

\begin{definition}\label{def:rhotau}
  Recall that we have $a\delta\equiv h$ mod $a$, then $\rho^{\otimes
    \frac{a\delta-h}{a}}_{2}$ induces a natural isomorphism
  $w^{*}L^{\vee}\otimes (M_{1}^{\otimes as\delta}\otimes M_{2}^{\otimes
    -ar\delta})\cong M_{1}^{\otimes sh}\otimes M_{2}^{\otimes -rh} $ which we
  still denote to be $\rho^{\otimes \frac{a\delta-h}{a}}_{2}$. We construct an
  natural isomorphism $$\rho_{\tau}: (M_{1}\otimes M_{3})^{\otimes
    as\delta}\otimes M_{2}^{-\otimes ar\delta}\rightarrow (f_{2}\circ
  f_{1})^{*}L^{\vee} $$ as the composition
  \begin{equation}
    \begin{aligned}
      &(M_{1}\otimes M_{3})^{\otimes as\delta}\otimes M_{2}^{-\otimes
        ar\delta}\stackrel{\cong}{\longrightarrow} M_{3}^{\otimes
        as\delta}\otimes (M_{1}^{\otimes as\delta}\otimes M_{2}^{-\otimes
        ar\delta})
      \\
      &\xrightarrow{ \rho_{1}\otimes \rho^{\otimes
          \frac{a\delta-h}{a}}_{2}}w^{*}L^{\vee}\otimes(M_{1}^{\otimes
        sh}\otimes M_{2}^{\otimes -rh}) \xrightarrow{
        \bar{\rho}_{\tau}}(f_{2}\circ f_{1})^{*}L^{\vee}\ .
    \end{aligned}
  \end{equation}
\end{definition}

\subsubsection{Defining the family map}\label{subsec:defing-mor}
We first collect some facts about line bundles over stacks. By a line bundle
$\mathcal L$ over a stack $\mathcal X$, we adopt the following definition: for
any scheme $S$ over $\mathbb C$ and any object $x$ of $\mathcal X(S)$, we
associate a line bundle over $S$ and we denote it by $x^{*}\mathcal L$, this
notation coincides with the notation of pullback of line bundles when we view an object of
$\mathcal X$ as a morphism $x:S\rightarrow \mathcal X$ by the $2-$Yoneda lemma. Let
$S'$ be another scheme and $x'\in \mathcal X(S')$. Let $g^{b}:x'\rightarrow x$
be a morphism over $g:S'\rightarrow S$. Then there is an isomorphism
$g^{*}x^{*}\mathcal L\rightarrow x^{'*}\mathcal L $, which we denote to be
$g^{b}$ as well, the collection of line bundles and morphisms between them
should satisfy obvious compatibility with composition of morphisms in $\mathcal
X$.

Let $f:\mathcal Y\rightarrow \mathcal X$ be a morphism between stacks, then the
pullback $f^{*}\mathcal L$ of the line bundle $\mathcal L$ along $f$ can be defined as follows: to each object $y\in
Ob(\mathcal Y)$, we associate the line bundle $f(y)^{*}\mathcal L$; the
morphisms between line bundles are defined in an obvious way. Let $f':\mathcal
Y\rightarrow \mathcal X$ be another morphism and $\alpha$ be a 2-isomorphism
from $f'$ to $f$, then we can define an isomorphism from $f^{*}\mathcal L$ to
$f^{'*}\mathcal L$ as follows, view $\mathcal X$ and $\mathcal Y$ as its
underlying category, then $\alpha$ is a natural transformation from $f'$ to $f$,
which are viewed as functors. For any object $y\in Ob(\mathcal Y)$ over a scheme
$S$, there is an associated morphism $\alpha_{y}:f^{'}(y)\rightarrow f(y)$ over
$S$, which determines a morphism $\alpha_{y}:f(y)^{*}\mathcal L\rightarrow
f^{'}(y)^{*}\mathcal L$ of line bundles over $S$. This defines the isomorphism
$f^{*}\mathcal L\rightarrow f^{'*}\mathcal L$, which we will still denote by
$\alpha$ by an abuse of notation. Let $g:\mathcal Z\rightarrow \mathcal Y$ and
$h:\mathcal X\rightarrow \mathcal W$ be two morphisms between stacks, we will
write $\alpha\circ g$ to mean the $2-$isomorphism from
$f'\circ g$ to $f\circ g$ and write $h\circ \alpha $ to mean the $2-$isomorphism from
$h\circ f'$ to $h\circ f$.


Now we define the family map $f:\mathcal C_{e}\rightarrow \P Y_{r,s} $ over
$\mathcal M_{e}$. To any object $$\tau=(w,\theta,
M_{1},M_{2},M_{3},\rho_{1},\rho_{2},s_{1},s_{2})\in \mathcal C_{e}(S)\ ,$$ we
associate an object $\phi=( f_{2} \circ f_{1}, (M_{1}\otimes M_{3})^{\otimes
  a\delta},M^{\otimes a\delta}_{2}, \rho_{\tau}, s_{1}^{a\delta},s_{2}^{a\delta}
)\in \mathcal \P Y_{r,s}(S)$, where $f_{1},f_{2},\rho_{\tau}$ (and
$\bar{f}_{1}$, $\alpha_{\tau}$, $\epsilon_{f_{2}}$, $\bar{\rho}_{\tau}$) are
defined in Definition \ref{def:induced-mor}, Definition \ref{def:rhotau} (and
Proposition \ref{prop:induced-mor}). Given another scheme $S'$ and another
object $\tau'\in \mathcal C_{e}(S')$, let $\phi'$ be the associated object in
$\P Y_{r,s}$, similarly for $f_{1}',f_{2}',\rho_{\tau'}$ (and $\bar{f'}_{1}$,
$\alpha_{\tau'}$, $\epsilon_{f'_{2}}$, $\bar{\rho}_{\tau'}$). Let
$u:=(g^b,\alpha_{1},\alpha_{2},\alpha_{3})$ be a morphism in $Hom_{\mathcal
  C_{e}}(\tau',\tau)$ over $g:S'\rightarrow S$, , we need some preparation
before we define a morphism from $\phi'$ to $\phi$ over $g$.

Using the notations in \ref{def:induced-mor}, by 2-Yoneda lemma, the morphism
$\alpha(u)$ induces a unique 2-isomorphism $\bar{\eta}_{1}:\bar{f}_{1}'
\rightarrow \bar{f}_{1}\circ g$ marking the following triangle commute:
$$\xymatrix{
  &S' \ar[rd]_{\bar{f}_{1}'} \ar[rr]^{g}& &S\ar[ld]^{\bar{f}_{1}} \\
  & &\mathbb B\boldsymbol{\mu_{a}} &\ . }$$ Let $\bar{\eta}_{1}:
(\bar{f}_{1}\circ g)^{*}\mathbb L \rightarrow (\bar{f}'_{1})^{*}\mathbb L$ be the
isomorphism induced from the 2-isomorphism $\bar{\eta}_{1}$. Then we have the
following commutative triangle

\begin{equation}\label{alpha*}
  \xymatrixcolsep{4pc}\xymatrix{
    &g^{*}\bar{f}_{1}^{*}\mathbb L\ar[d]^{\cong}\ar[rd]^{\alpha(u)}   \\
    &(\bar{f}_{1}\circ
    g)^{*}\mathbb L \ar[r]^{\bar{\eta}_{1}} &(\bar{f}'_{1})^{*}\mathbb L\ .
  }
\end{equation}

In the following, to simplify the presentation, we will treat
$g^{*}\bar{f}_{1}^{*}\mathbb L$ and $(\bar{f}_{1}\circ g)^{*}\mathbb L$ are the
same via the vertical isomorphism in \ref{alpha*} defined in an obvious way,
thus $\bar{\eta}_{1}=\alpha(u)$ above, and $\bar{\eta}_{1}$ can be also viewed as
a map from $g^{*}\bar{f}_{1} ^{*}\mathbb L$. The same rule applies to other
pullbacks of a line bundle along a composition of maps.

Note that $\bar{\eta}_{1}$ induces a 2-isomorphism $\eta_{1}: (g\times id)\circ
f_{1}'\rightarrow f_{1}\circ g$ satisfying that
$\bar{\eta}_{1}=\mathrm{pr}_{2}\circ \eta_{1}$ and making the following square
commute:
$$\xymatrix{
  &S' \ar[d]^{f_{1}'} \ar[r]^{g} &S\ar[d]^{f_{1}} \\
  &S'\times \mathbb B\boldsymbol{\mu_{a}} \ar[r]^{g\times id} &S\times \mathbb
  B\boldsymbol{\mu_{a}}\ . }$$ For any line bundle $M$ over $S$ and any integer
$n$, the $2-$isomorphism $\eta_{1}$ induces an isomorphism $\eta_{1}:
(f_{1}\circ g)^{*}(M\boxtimes \mathbb L^{\otimes n}) \rightarrow ((g\times
id)\circ f'_{1})^{*}(M\boxtimes \mathbb L^{\otimes n})$ such that we have the
following commutative square using \ref{alpha} and Lemma \ref{lem:induced-mor}

\begin{equation}\label{eta1}
  \xymatrixcolsep{4pc}\xymatrix{
    &g^{*}M\otimes g^{*}(M_{1}^{\otimes sn}\otimes
    M_{2}^{\otimes -rn})\ar[d]^{\cong}  \ar[rr]^{id \otimes(\alpha_{1}^{-\otimes sn}\otimes \alpha_{2}^{\otimes rn})
    }& &g^{*}M\otimes (M_{1}^{'-\otimes sn}\otimes
    M_{2}^{'\otimes rn})\ar[d]_{\cong} \\
    &g^{*}f_{1}^{*}(M\boxtimes \mathbb L^{\otimes n}) 
    \ar[rr]^{\eta_{1}}& &((g\times
    id)\circ  f'_{1})^{*}(M\boxtimes \mathbb L^{\otimes n}) \ .
  }
\end{equation}
Here the left vertical map is induced from Lemma \ref{lem:induced-mor} by taking
$t_{1}\times t_{2}=f_{1}$ and pullback along the map $g$, while the right
vertical map is induced from Lemma \ref{lem:induced-mor} by taking $t_{1}\times
t_{2}=(g\times id)\circ f'_{1}$.

Using the equivalence of the stack $\mathrm{RepHom}(\mathbb B
\mathbf{\mu}_{a},Y)$ with the $I_{\mu_{a}}Y$ from ~\cite[\S
3.2]{Dan_Abramovich_2008}, the morphism $g^b$ between $(w',\theta')\in
I_{\mu_{a}}Y(S')$ and $(w,\theta)\in I_{\mu_{a}}Y(S)$ induces a 2-isomorphism $
\eta_{2}:f_{2}' \rightarrow f_{2}\circ(g\times id ) $ marking the following
triangle commute
$$\xymatrix{
  &S'\times \mathbb B\boldsymbol{\mu_{a}}\ar[rd]_{f'_{2}} \ar[rr]^{g\times id} & &S\times \mathbb B\boldsymbol{\mu_{a}}\ar[ld]^{f_{2}}\\
  & &Y & \ .}$$ The 2-isomorphism $\eta_{2}$ induces a morphism from
$(f_{2}\!\circ\!(g\times id))^{*}L^{\vee}$ to  m,$(f'_{2})^{*}L^{\vee}$ which makes
the following diagram commutes:
\begin{equation}\label{eta2}
  \xymatrix{
    &g^{*}w^{*}L^{\vee}\boxtimes \mathbb L_{a,-h}\ar[d]_{(g\times id)^{*}\epsilon_{f_{2}}}\ar[rr]^{g^b\boxtimes id
    }& &(w')^{*}L^{\vee}\boxtimes \mathbb L_{a,-h}\ar[d]^{\epsilon_{f'_{2}}}\\
    & (f_{2}\!\circ\!(g\!\times\!
    id))^{*}L^{\vee} \ar[rr]^-{\eta_{2}} & & (f'_{2})^{*}L^{\vee}  \ .
  }
\end{equation}
Indeed, using the notations in Proposition \ref{prop:induced-mor}, the
2-isomorphism $i_{S}\circ \eta_{2}$ induces a morphism $g^{*}w^{*}L^{\vee}
\rightarrow (w')^{*}L^{\vee}$, which is exactly induced by $g^{b}$. As $g^{b}$
can be view a morphism from $w'$ to $w$ which commutes with $\theta'$ and
$\theta$, then by descent, we get the above diagram.

Then we make the composition $\eta_{2}\star \eta_{1}:= (\eta_{1}\circ f_{2})
\circ (f_{1}'\circ \eta_{2})$,\footnote{This is the so-called Godement product.}
which is a 2-isomorphism from $f'_{2}\circ f'_{1}$ to $f_{2}\circ f_{1}\circ g$.
When we view $f_{2}'\circ f_{1}'$ (resp. $f_{2}\circ f_{1}$) as a point in
$Y(S')$ (resp. $Y(S)$), then $\eta_{2}*\eta_{1}$ can be viewed as a morphism
over $g$ in $Y$ from the point $f_{2}'\circ f_{1}'$ to the point $f_{2}\circ
f_{1}$.

Now we can define the family map:
\begin{definition}[\textbf{Family map}]\label{def:familymap}
  With the notations as above, for any object
 $$\tau=(w,\theta, M_{1},M_{2},M_{3},\rho_{1},\rho_{2},s_{1},s_{2})\in
 \mathcal C_{e}(S) \ ,$$ we define the image $f(\tau)$ under the family map $f$
 to be
 $$\phi=( f_{2} \circ
 f_{1}, (M_{1}\otimes M_{3})^{\otimes a\delta},M^{\otimes a\delta}_{2},
 \rho_{\tau}, s_{1}^{a\delta},s_{2}^{a\delta} )\in \mathcal \P Y_{r,s}(S)\ .$$
 Given a morphism $u:\tau'\rightarrow \tau$ in $\mathcal C_{e}$, we define the
 morphism $f(u): \phi' \rightarrow \phi$ to be $(\eta_{2}\star
 \eta_{1},(\alpha_{1}\otimes \alpha_{3})^{\otimes a\delta},(\alpha_{2})^{\otimes
   a\delta})$.
\end{definition}

We should check that $f(u)$ is indeed a well-defined morphism, for which we only
need to check the following diagram commutes
\begin{equation}\label{eq:comm-verify}
  \xymatrixcolsep{4pc}\xymatrix{
    &g^{*}(M_{1}\!\otimes \!M_{3})^{\otimes as\delta}\!\otimes\!
    g^{*}M_{2}^{-\otimes ar\delta}\ar[d]_{g^{*}\rho_{\tau}}
    \ar[rr]^{(\alpha_{1}\otimes \alpha_{3})^{\otimes as\delta}\otimes \alpha^{-\otimes ar\delta}_{2}} &
    &(M_{1}'\!\otimes\! M_{3}')^{\otimes as\delta}\!\otimes\! M_{2}'^{-\otimes ar\delta}\ar[d]^{\rho_{\tau'}}\\
    & g^{*}(f_{2}\circ f_{1})^{*}L^{\vee}  \ar[rr]^{\eta_{2}\star\eta_{1}} & &(f'_{2}\circ f'_{1})^{*}L^{\vee}\ .  
  }
\end{equation}
Note that we have the following diagram
 $$\xymatrixcolsep{4pc}\xymatrix{
   &g^{*}(M_{1}\!\otimes \!M_{3})^{\otimes as\delta}\!\otimes\!
   g^{*}M_{2}^{-\otimes ar\delta}\ar[d]_{\cong} \ar[rr]^{(\alpha_{1}\otimes
     \alpha_{3})^{\otimes as\delta}\otimes \alpha^{-\otimes ar\delta}_{2}} &
   &(M_{1}'\!\otimes\! M_{3}')^{\otimes as\delta}\!\otimes\! M_{2}^{'-\otimes ar\delta}\ar[d]^{\cong}\\
   & g^{*}M_{3}^{\otimes as\delta}\!\otimes\! g^{*}(M_{1}^{\otimes
     as\delta}\otimes M_{2}^{-\otimes ar\delta})\ar[d]_{g^{*}\rho_{1}\otimes
     g^{*}\rho_{2}^{\frac{a\delta-h}{a}}} \ar[rr]^{( \alpha_{3})^{\otimes
       as\delta}\otimes \alpha_{1}^{as\delta}\otimes\alpha^{-\otimes
       ar\delta}_{2}} & &M_{3}^{'\otimes as\delta}\!\otimes\! (M_{1}^{'\otimes
     as\delta}\otimes M_{2}^{'-\otimes
     ar\delta})\ar[d]^{\rho'_{1}\otimes (\rho_{2}')^{\frac{a\delta-h}{a}}}\\
   &g^{*}w^{*}L^{\vee}\otimes g^{*}(M_{1}^{\otimes sh}\otimes M_{2}^{\otimes
     -rh})\ar[rr]^{g^b \otimes (\alpha_{1}^{\otimes sh}\otimes
     \alpha_{2}^{-\otimes rh})}& &(w')^{*}L^{\vee}\otimes (M_{1}^{'\otimes
     sh}\otimes M_{2}^{'\otimes -rh}) \ , }$$ by the definition of
 $\rho_{\tau}$, we see that showing that \ref{eq:comm-verify} is a commutative
 diagram is equivalent to show that the following diagram commutes:
 \begin{equation}\label{eq:comm-hor}
   \xymatrixcolsep{4pc}\xymatrix{
     &g^{*}w^{*}L^{\vee}\otimes g^{*}(M_{1}^{\otimes sh}\otimes
     M_{2}^{\otimes -rh})\ar[d]_{g^{*}\bar{\rho}_{\tau}} \ar[rr]^{g^b \otimes (\alpha_{1}^{\otimes sh}\otimes \alpha_{2}^{-\otimes rh})}& &(w')^{*}L^{\vee}\otimes (M_{1}^{'\otimes sh}\otimes
     M_{2}^{'\otimes -rh})\ar[d]^{\bar{\rho}_{\tau'}}\\
     &g^{*}(f_{2}\circ f_{1})^{*}L^{\vee}\ar[rr]^{\eta_{2}\star\eta_{1}} & &(f'_{2}\circ f'_{1})^{*}L^{\vee} \ .
   }
 \end{equation}

 Denote $\bar{\rho}_{\tau,\tau'}:=f_{1}^{*}(g\times id)^{*}\epsilon_{f_{2}}$,
 then $\bar{\rho}_{\tau,\tau'} $ induces an isomorphism from
 $g^{*}w^{*}L^{\vee}\otimes (M_{1}^{'\otimes sh}\otimes M_{2}^{'\otimes -rh})$
 to $(f_{2}\circ(g\times id)\circ f'_{1})^{*}L^{\vee}$. The 2-isomorphism
 $\eta_{1}\circ f_{2}$ induces a morphism $\eta_{1}\circ f_{2}: g^{*}(f_{2}\circ
 f_{1})^{*}L^{\vee}\rightarrow (f_{2}\circ(g\times id)\circ f'_{1})^{*}L^{\vee}
 $. Take $M=w^{*}L^{\vee}$, $n=-h$ in \ref{eta1}, using Proposition
 \ref{prop:induced-mor}, we get the following commutative diagram
$$\xymatrixcolsep{4pc}\xymatrix{
  &g^{*}w^{*}L^{\vee}\otimes g^{*}(M_{1}^{\otimes sh}\otimes M_{2}^{\otimes
    -rh})\ar[d]^{g^{*}\bar{\rho}_{\tau'}} \ar[rr]^{id
    \otimes(\alpha_{1}^{\otimes sh}\otimes \alpha_{2}^{-\otimes rh}) }&
  &g^{*}w^{*}L^{\vee}\otimes (M_{1}^{'\otimes sh}\otimes M_{2}^{'\otimes
    -rh})\ar[d]_{\bar{\rho}_{\tau',\tau}}
  \\
  & g^{*}(f_{2}\circ f_{1})^{*}L^{\vee} \ar[rr]^{\eta_{1}\circ f_{2}} &
  &(f_{2}\circ(g\times id)\circ f'_{1})^{*}L^{\vee} \ . }$$

The 2-isomorphism $f_{1}'\circ \eta_{2}$ induces a morphism $f_{1}'\circ
\eta_{2}: (f_{2}\circ(g\times id)\circ f'_{1})^{*}L^{\vee} \rightarrow
(f'_{2}\circ f'_{1})^{*}L^{\vee}$. Pull back the square \ref{eta2} via
$f_{1}^{'*}$, we get the commutative square:
$$\xymatrix{
  &g^{*}w^{*}L^{\vee}\otimes (M_{1}^{'\otimes sh}\otimes M_{2}^{'\otimes
    -rh})\ar[d]_{\bar{\rho}_{\tau',\tau}} \ar[rr]^{g^b\otimes id }&
  &(w')^{*}L^{\vee}\otimes (M_{1}^{'\otimes sh}\otimes
  M_{2}^{'\otimes -rh})\ar[d]^{\bar{\rho}_{\tau'}}\\
  &(f_{2}\circ(g\times id)\circ f'_{1})^{*}L^{\vee} \ar[rr]^{f_{1}'\circ
    \eta_{2}} & & (f'_{2}\circ f'_{1})^{*}L^{\vee} \ . }$$ Join the above two
squares horizontally, this verifies the commutativity of \ref{eq:comm-hor}.

\begin{remark}\label{rmk:ex-toric}
  When $Y$ is a complete intersection in a toric stack, our construction of the
  family map is equivalent to the construction in~\cite[\S
  5.3.2]{wang19_mirror_theor_gromov_witten_theor_without_convex}.
\end{remark}

\subsection{Computation of edge contribution}\label{subsec:computeedge}
Recall that in \S \ref{sec:prop-Lagrangian}, we define a $\C^{*}-$action on $\P
Y_{r,s} $ which is induced from the $\C^{*}-$action on $(L^{\vee})^{0}
\!\times\! U$:
$$ t\cdot (l,x,y)=(tl, x,y)\ . $$
In the language of functor of points, for any $t\in \C^{*}$ and any scheme $S$,
this action will change any object $(w,M_{1},M_{2},\rho, s_{1},s_{2})$ of $Ob(\P
Y_{r,s} (S)) $ to $(w, M_{1},M_{2},t\cdot \rho, s_{1},s_{2})$ of $Ob(\P
Y_{r,s}(S))$. where $t\cdot \rho$ means that we compose $\rho$ by the
automorphism of $w^{*}L^{\vee}$ obtained by scaling by $t$.

Now we define a $\C^{*}-$action on $\mathcal C_{e}$ such that the morphism $f$
is $\C^{*}-$equivariant. Recall that $\mathcal C_{e}$ can be represented as
$[(L^{\vee})^{0}\times U / T_{ar,as} \!\times\! \C^{*} ]$ given by the action
$$(l,x,y)(t_{1},t_{2},t_{3})= (t_{3}^{-as\delta}l, t_{3}t_{1}x,t_{2}y) \ .$$
We define the $\C^{*}-$action on $(L^{\vee})^{0}\times U$ in the way that
$$t\cdot (l,x,y)= (t\cdot l, x, y)\ .$$
One can see this action descends to be $\C^{*}-$action on $\mathcal C_{e}$. In
the language of functor of points, for any $t\in \C^{*} $ and any scheme $S$
this $\C^{*}-$action will change any object $(w,\theta,
M_{1},M_{2},M_{3},\rho_{1},\rho_{2}, s_{1},s_{2})\in \mathcal C_{e}(S)$ to
$(w,\theta, M_{1},M_{2},M_{3},t\cdot \rho_{1},\rho_{2}, s_{1},s_{2})\in \mathcal
C_{e}(S)$, where $t\cdot \rho_{1}$ means that we compose $\rho_{1}$ by the
automorphism of $w^{*}L^{\vee}$ obtained by scaling by $t$. Then this action
will change the associated $\rho_{\tau}$ by $t\rho_{\tau}$, which implies that
$f$ is $\C^{*}-$equivariant.

Let $\chi_{1}$ (resp. $\chi_{2}$) be the character associated to $T_{ar,as}$ by
sending $(t_{1},t_{2})\in T_{ar,as}$ to $t_{1}$ (resp. $t_{2}$). Then we can
construct the line bundle $L_{\chi_{1}}$(resp. $L_{\chi_{2}}$) over $\mathcal
C_{e}$ as the stack quotient $[(L^{\vee})^{0}\times U \times \mathbb C/
T_{ar,as} \!\times\! \C^{*}]$ using the Borel construction:
$$(l,x,y,z)\sim (t_{3}^{-as\delta}l,t_{3}t_{1}x,t_{2}y,t_{1}z)\; (\mathrm{resp}.\;  (t_{3}^{-as\delta}l,t_{3}t_{1}x,t_{2}y,t_{2}z))$$
for any $(l,x,y,z)\in (L^{\vee})^{0}\times U \times \mathbb C$ and
$(t_{1},t_{2},t_{3})\in T_{ar,as} \!\times\! \C^{*}$. Then for any scheme $S$
and any object $\tau:= (w,\theta, M_{1},M_{2},M_{3},\rho_{1},\rho_{2},
s_{1},s_{2})$ of $Ob(\mathcal C_{e}(S))$, we have $\tau^{*}L_{\chi_{1}}=M_{1}$
and $\tau^{*}L_{\chi_{2}}=M_{2}$. We also note $M_{3}$ is isomorphic to the line
bundle $\tau^{*}\pi^{*}L^{\vee}$. Using the stack quotient representation, the
line bundles $L_{\chi_{1}}$ and $L_{\chi_{2}}$ also carry a $\C^{*}-$action such
that the fiber of $L_{\chi_{1}}$ over $\infty-$section(corresponding to the
coordinate $y=0$) is of weight $\frac{-1}{as\delta}$, the fiber of
$L_{\chi_{2}}$ over $D_{\infty}$ is of weight $0$. Then we have the following:
\begin{proposition}\label{prop:equiv-sec}
  The coordinate $x$ of $(L^{\vee})^{0}\times U$ descents to be a
  $\C^{*}-$equivariant section of the $\C^{*}-$equivariant line bundle
  $L_{\chi_{1}}\otimes (\pi^{*}L^{\vee})^{\frac{1}{as\delta}}\otimes \mathbb
  C_{\frac{\lambda}{as\delta}}$ over $\mathcal C_{e}$, and the coordinate $y$ of
  $(L^{\vee})^{0}\times U$ descents to be a $\C^{*}-$equivariant section of the
  $\C^{*}-$equivariant line bundle $L_{\chi_{2}}$ over $\mathcal C_{e}$.
\end{proposition}

Fix the edge $e=\{h_{0},h_{\infty}\}$ with the decoration as above, we will
denote
$$\mathcal K_{0,h_{0}\sqcup h_{\infty}}(\P
Y_{r,s}, (0,\frac{\delta}{r})):= \mathcal K_{0,2}(\P Y_{r,s}, (0,
\frac{\delta}{r}))\cap
ev^{-1}_{1}(\bar{I}_{(c^{-1},\mathrm{e}^{\frac{\delta}{s}},1)}\P Y_{r,s}) \cap
ev^{-1}_{2}(\bar{I}_{(c,1,\mathrm{e}^{\frac{\delta}{r}})}\P Y_{r,s})\ .$$ Here
we label the first marking by $h_{0}$ and the second marking by $h_{\infty}$.
Let $\mathcal K_{0,h_{0}\sqcup h_{\infty}}(\P Y_{r,s},
(0,\frac{\delta}{r}))^{\C^{*}}$ be the $\C^{*}-$fixed loci of $\mathcal
K_{0,h_{0}\sqcup h_{\infty}}(\P Y_{r,s}, (0,\frac{\delta}{r}))$, write $\mathcal
K:=\mathcal K_{0,h_{0}\sqcup h_{\infty}}(\P Y_{r,s},
(0,\frac{\delta}{r}))^{\C^{*}}$ for short, we will first show the following:
\begin{proposition}\label{prop:etale-ev}
  Let $[\mathcal K]^{vir}$ be the virtual cycle associated to the $\C^{*}-$fixed
  part of the perfect obstruction theory of $\mathcal K_{0,h_{0}\sqcup
    h_{\infty}}(\P Y_{r,s}, (0,\frac{\delta}{r}))$ to $\mathcal K$. We have
  that $[\mathcal K]^{\mathrm{vir}}=[\mathcal K] $ in $A_{*}(\mathcal K)$ and
  the evaluation map
$$ev_{h_{\infty}}: \mathcal K\rightarrow \bar{I}_{(c,1,e^{\frac{\delta}{r}})}\P Y_{r,s}\, $$
is \'etale and finite of degree $\frac{1}{a\delta}$.
\end{proposition}
\begin{proof}
  Let $\mathbb E$ be the restriction of the perfect obstruction
  theory\footnote{Here a perfect obstruction is a morphism $\mathbb T\rightarrow
    \mathbb E$ in~\cite{behrend97_intrin_normal_cone}, where $\mathbb T= \mathbb
    L^{\vee}$ is the derived dual of the cotangent complex.} for $\mathcal
  K_{0,h_{0}\sqcup h_{\infty}}(\P Y_{r,s}, (0,\frac{\delta}{r}))$ to $\mathcal
  K$, then the $\C^{*}-$fixed part of $\mathbb E$, which we denote to be
  $\mathbb E^{fix}$, is a perfect obstruction theory for $\mathcal K$
  by~\cite{Graber_1999}. We will show that $\mathbb E= \mathbb E^{fix}$ and it's
  isomorphic to the locally free sheaf $ev^{*}_{h_{\infty}}\mathbb
  T_{\bar{I}_{(c,1,e^{\frac{\delta}{r}})}\P Y_{r,s}}$. Note this implies that $\mathbb
  E^{fix}\cong \mathbb T_{\mathcal K}\cong ev^{*}_{h_{\infty}}\mathbb
  T_{\bar{I}_{(c,1,e^{\frac{\delta}{r}})}\P Y_{r,s}} $. This will implies the
  \'etaleness part, while the finiteness part will follow the properness of
  $ev_{h_{\infty}}$, Lemma \ref{thm:classify-mor} and Proposition
  \ref{prop:ker-aut}.

  It suffices to check $ \mathbb E=\mathbb E^{fix}\cong
  ev^{*}_{h_{\infty}}\mathbb T_{\bar{I}_{(c,1,e^{\frac{\delta}{r}})}\P Y_{r,s}}
  $ locally on every $\mathbb C$-point. For any $\C-$point $[f]$ of $\mathcal K$
  represented by a twisted stable map $f$, the tangent space $T_{[f]}$ and the
  obstruction space $Ob_{[f]}$ of the space $\mathcal K$ at the point $[f]$ fit
  into the following exact sequence of $\C^{*}-$presentations\footnote{All the
    arrows in this proof are $\C^{*}-$equivaraint unless otherwise stated.}
  \begin{equation}\label{eqn:B}
    \begin{aligned}
      0\to H^{0}(C,\mathbb T_{C}(-[q_{0}]-[q_{\infty}]))\to H^{0}(C,f^{*}\mathbb
      T_{\P
        Y_{r,s}})\to T_{[f]}\\
      \to H^{1}(C,\mathbb T_{C}(-[q_{0}]-[q_{\infty}]))\to H^{1}(C,f^{*}\mathbb
      T_{\P Y_{r,s}})\to Ob_{[f]}\to 0\ .
    \end{aligned}
  \end{equation}
  We have the following two exact sequences of vector bundles over $\P Y_{r,s}$
  (which are also $\C^{*}$-equivariant)
  \begin{equation}\label{eq:firstseq}
    \xymatrix{
      0\ar[r] &\mathbb T_{\P Y_{r,s}/ Y} \ar[r] &\mathbb T_{\P Y_{r,s}
      }\ar[r] &\mathrm{pr}_{r,s}^{*}\mathbb T_{Y} \ar[r] &0\ ,
    }
  \end{equation}
  and the Euler exact sequence
  \begin{equation}\label{eq:secondseq}
    \xymatrix{
      0\ar[r] &\mathbb C \ar[r] &\mathcal O(\mathcal D_{0})\oplus \mathcal O(\mathcal D_{\infty})\ar[r] &\mathbb T_{\P Y_{r,s}/Y}\ar[r] &0\ .
    }
  \end{equation}
  Then we have the following two right exact sequences
  \begin{equation}\label{eqn:C}
    \begin{aligned}
      0\to H^{0}(C,f^{*}\mathbb T_{\P Y_{r,s}/Y}) \to H^{0}(C,f^{*}\mathbb T_{\P Y_{r,s}})\to  H^{0}(C,f^{*}\mathrm{pr}_{r,s}^{*}\mathbb T_{Y})\\
      \to H^{1}(C,f^{*}\mathbb T_{\P Y_{r,s}/Y}) \to H^{1}(C,f^{*}\mathbb T_{\P
        Y_{r,s}})\to H^{1}(C,f^{*}\mathrm{pr}_{r,s}^{*}\mathbb T_{Y})\to 0\ .
    \end{aligned}
  \end{equation}
  and
  \begin{equation}\label{eqn:D}
    \begin{aligned}
      0 \to H^{0}(C,\mathcal O) \to H^{0}(C,f^{*} \mathcal O(\mathcal D_{0})\oplus f^{*}\mathcal O(\mathcal D_{\infty}))\to  H^{0}(C,f^{*}T_{\P Y_{r,s} /Y})\\
      \to H^{1}(C,\mathcal O) \to H^{1}(C,f^{*} \mathcal O(\mathcal D_{0})\oplus
      f^{*}\mathcal O(\mathcal D_{\infty}))\to H^{1}(C,f^{*}T_{\P Y_{r,s}
        /Y})\to 0\ .
    \end{aligned}
  \end{equation}
  By Lemma \ref{thm:classify-mor}, we have $$f^{*} \mathcal O(\mathcal
  D_{0})\oplus f^{*}\mathcal O(\mathcal D_{\infty})\cong \mathcal
  O(a\delta[q_{0}]+a\delta[q_{\infty}])$$ over $C\cong \P^{1}_{ar,as}$, as $r,s$
  are sufficiently large primes, then we
  have
$$ H^{1}(C,f^{*} \mathcal O(\mathcal D_{0})\oplus f^{*}\mathcal O(\mathcal
D_{\infty}))=0\ ,$$ which implies that $H^{1}(C,f^{*}\mathbb T_{\P Y_{r,s}
  /Y})=0$ by \ref{eqn:D}. As $H^{1}(C,\mathbb T_{C}(-[q_{0}]-[q_{\infty}]))=0$,
we have $$Ob_{[f]}\cong H^{1}(C,f^{*}\mathrm{pr}^{*}_{r,s}\mathbb T_{Y})$$
and $$T_{[f]}\cong H^{0}(C,f^{*}\mathbb T_{\P Y_{r,s}})/ H^{0}(C,\mathbb
T_{C}(-[q_{0}]-[q_{\infty}])) $$ by \ref{eqn:B} and \ref{eqn:C}.

Moreover, the exact sequence \ref{eqn:D} becomes
\begin{equation}\label{eqn:D*}
  \begin{aligned}
    0 \to \mathbb C \to \mathbb C\oplus \mathbb C \to H^{0}(C,f^{*}\mathbb T_{\P
      Y_{r,s} /Y})\to 0 \to 0\to 0\to 0\ ,
  \end{aligned}
\end{equation} 
which implies that $H^{0}(C,f^{*}\mathbb T_{\P Y_{r,s} /Y})\cong \mathbb C$(of
weight 0). Furthermore, the composition
$$H^{0}(C,T_{C}(-[q_{0}]-[q_{\infty}]))\rightarrow H^{0}(C,f^{*}\mathbb T_{\P Y_{r,s}})\rightarrow H^{0}(C,f^{*}\mathrm{pr}_{r,s}^{*}\mathbb T_{Y})$$
where the arrows are coming from \ref{eqn:B} and \ref{eqn:C}, is equal to zero
as the composition is induced from the composition of sheaves
$$\mathbb T_{C}(-[q_{0}]-[q_{\infty}])\rightarrow f^{*}\mathbb T_{\P Y_{r,s}}\rightarrow
f^{*}\mathrm{pr}_{r,s}^{*}\mathbb T_{Y}\ ,$$ which is zero. Note
$H^{0}(C,\mathbb T_{C}(-[q_{0}]-[q_{\infty}]))$ is one dimensional, which
corresponds to the space of tangent vector field of orbi-$\P^{1}$ vanishing on
the two markings, then the space $H^{0}(C,\mathbb T_{C}(-[q_{0}]-[q_{\infty}]))$
maps isomorphically into the subspace $H^{0}(C,f^{*}\mathbb T_{\P Y_{r,s}/Y})
\subset H^{0}(C,f^{*}\mathbb T_{\P Y_{r,s}})$ via \ref{eqn:C}, which implies
that
\begin{equation*}
  \begin{split}
    &T_{[f]}\cong H^{0}(C,f^{*}\mathbb T_{\P
      Y_{r,s}})/ H^{0}(C,\mathbb T_{C}(-[q_{0}]-[q_{\infty}]))\\
    &\cong H^{0}(C,f^{*}\mathbb T_{\P Y_{r,s}})/ H^{0}(C,f^{*}\mathbb T_{\P
      Y_{r,s}/Y})\cong H^{0}(C,f^{*}\mathrm{pr}^{*}_{r,s}\mathbb T_{\P
      Y_{r,s}})\ .
  \end{split}
\end{equation*}

Let $\pi:C\rightarrow \mathrm{Spec}(\mathbb C)$ be the projection to a point,
next we show that $$R^{\bullet}\pi_{*}f^{*}\mathrm{pr}^{*}_{r,s}T_{\P
  Y_{r,s}}\cong (ev^{*}_{q_{\infty}}\mathbb
T_{\bar{I}_{(g,1,e^{\frac{\delta}{r}})}\P Y_{r,s}})|_{[f]}\ ,$$ which will imply
that $$ T_{[f]}\cong (ev^{*}_{q_{\infty}}\mathbb
T_{\bar{I}_{(g,1,e^{\frac{\delta}{r}})}\P Y_{r,s}})|_{[f]}\ ,\quad Ob_{[f]}=0\ ,
 $$
 in particular $T_{[f]}$ is $\C^{*}-$fixed. Applying
 $R^{\bullet}\pi_{*}f^{*}(-)$ to the first sequence \ref{eq:firstseq}, we get a
 distinguished triangle
$$\xymatrix{
  &R^{\bullet}\pi_{*} f^{*}\mathbb T_{\P Y_{r,s}/ Y} \ar[r]
  &R^{\bullet}\pi_{*}f^{*}\mathbb T_{\P Y_{r,s} }\ar[r]
  &R^{\bullet}\pi_{*}f^{*}\mathrm{pr}_{r,s}^{*}\mathbb T_{Y} \ar[r]^-{+1}
  &R^{\bullet}\pi_{*} f^{*}\mathbb T_{\P Y_{r,s}/ Y}[1]\ . }$$ Let $b_{\infty}:
q_{\infty}\hookrightarrow C$ be the gerbe section of $C$ corresponding to the
half-edge $h_{\infty}$. Using the divisor sequence for the section $b_{\infty}$
$$0\rightarrow \mathcal O(-q_{\infty})\rightarrow \mathcal O_{C}\rightarrow \mathcal O_{q_{\infty}}\rightarrow 0 \ ,$$
tensor with $f^{*}\mathrm{pr}_{r,s}^{*}\mathbb T_{Y}$ and taking
$R^{\bullet}\pi_{*}$, we have
$$ R^{\bullet}\pi_{*}f^{*}\mathrm{pr}_{r,s}^{*}\mathbb
T_{Y}=R^{\bullet}(\pi\circ b_{\infty})_{*}(f\circ
b_{\infty})^{*}(\mathrm{pr}_{r,s}^{*}\mathbb T_{Y})\ ,$$ where we use the fact
that $R^{\bullet}\pi_{*}f^{*}\mathrm{pr}_{r,s}^{*}\mathbb T_{Y}(-q_{\infty})=0$.
Next we show that
$$ R^{\bullet}(\pi\!\circ\!
b_{\infty})_{*}(f\!\circ\! b_{\infty})^{*}(\mathrm{pr}_{r,s}^{*}\mathbb T_{Y})
\cong (ev^{*}_{q_{\infty}}\mathbb T_{\bar{I}_{(c,1,e^{\frac{\delta}{r}})}\P
  Y_{r,s}})|_{[f]} \ ,$$

First by the tangent bundle lemma in~\cite[\S 3.6.1]{Dan_Abramovich_2008}, we
have
$$ (\pi\!\circ\!
b_{\infty})_{*}(f\!\circ\! b_{\infty})^{*}(\mathbb T_{\P Y_{r,s}}) \cong
(ev^{*}_{q_{\infty}}\mathbb T_{\bar{I}_{(c,1,e^{\frac{\delta}{r}})}\P
  Y_{r,s}})|_{[f]} \ ,$$ thus we only need to show that
$$ R^{\bullet}(\pi\!\circ\!
b_{\infty})_{*}(f\!\circ\! b_{\infty})^{*}(\mathrm{pr}_{r,s}^{*}\mathbb
T_{Y})\cong (\pi\!\circ\! b_{\infty})_{*}(f\!\circ\! b_{\infty})^{*}\mathbb
T_{\P Y_{r,s}}\ .$$ As $(\pi\circ b_{\infty})_{*}$ is exact on coherent sheaf,
we have $R^{\bullet}(\pi\circ b_{\infty})_{*}= (\pi\circ b_{\infty})_{*}$ and
the following exact sequence of vector spaces:
\begin{equation}\label{eq:thirdseq}
  \begin{aligned}
    0\to (\pi\!\circ\! b_{\infty})_{*}(f\!\circ\! b_{\infty})^{*}\mathbb T_{\P
      Y_{r,s}/ Y} \to (\pi\!\circ\! b_{\infty})_{*}(f\!\circ\!
    b_{\infty})^{*}\mathbb T_{\P Y_{r,s} }\to (\pi\!\circ\!
    b_{\infty})_{*}(f\!\circ\! b_{\infty})^{*}\mathrm{pr}_{r,s}^{*}\mathbb T_{Y}
    \to 0\ ,
  \end{aligned}
\end{equation}
we are left to show that $(\pi\!\circ\! b_{\infty})_{*}(f\!\circ\!
b_{\infty})^{*}\mathbb T_{\P Y_{r,s}/ Y}=0$, which can be checked by the next
lemma \ref{lem:verify-etale-edge} by pulling back to every $\mathbb C$-point
using the idea of proving the tangent bundle lemma~\cite[\S
3.6.1]{Dan_Abramovich_2008}.
\end{proof}

\begin{lemma}\label{lem:verify-etale-edge}
  Let $x:=[f]$ be a $\mathbb C$-point of $\mathcal K$ represented by a twisted
  stable map $f:C\rightarrow \P Y_{r,s}$ as in Lemma \ref{thm:classify-mor}. Use
  the notation in the proof of Proposition \ref{prop:etale-ev}, we have
 $$(\pi\!\circ\!
 b_{\infty})_{*}(f\!\circ\! b_{\infty})^{*}\mathbb T_{\P Y_{r,s}/ Y}=0\ .$$
\end{lemma}
\begin{proof}
  The gerbe marking $q_{\infty}$ of $C$ is canonically isomorphic to the
  classifying stack $\mathbb B\boldsymbol{\mu_{ar}}$, which gives a unique lift
  of the $\mathbb C-$point $x$ in the gerbe marking $q_{\infty}$, which we
  denote to be $z$. The isotropy group $Aut(z)$ of $z$ in $\mathbb B
  \boldsymbol{\mu_{ar}}$ is canonically isomorphic to the cyclic group
  $\boldsymbol{\mu}_{ar} $, and the generator of $\boldsymbol{\mu}_{ar}$ acts on
  the vector space $f^{*}\mathbb T_{\P Y_{r,s}/Y}|_{z}$ via the multiplication
  by $e^{\frac{\delta}{r}}$ due to the very choice of the multiplicity
  associated to the half-edge $h_{\infty}$. Note the space $(\pi\!\circ\!
  b_{\infty})_{*}(f\!\circ\! b_{\infty})^{*}\mathbb T_{\P Y_{r,s}/ Y}$ is given
  by the $\boldsymbol{\mu}_{ar}-$invariant space $(f^{*}\mathbb T_{\P
    Y_{r,s}/Y}|_{z})^{\boldsymbol{\mu}_{ar}}$, then it is equal to zero.
\end{proof}

The family map $f$ constructed in Section \ref{subsec:defing-mor} induces a
morphism $$g:\mathcal M_{e}\rightarrow \mathcal K$$ by the universal property of
moduli stack $ \mathcal K$. Then $g$ is surjective by Lemma
\ref{thm:classify-mor}. Furthermore,. we have:

\begin{proposition}\label{prop:deg-edge}
  The morphism $g:\mathcal M_{e}\rightarrow \mathcal K$ is finite \'etale of
  degree $\frac{1}{as}$. The evaluation map $ev_{h_{\infty}}$ from $ \mathcal K$
  to $\bar{I}_{(c, 1, \frac{\delta}{r})}\P Y_{r,s} $ is finite \'etale of degree
  $\frac{1}{a\delta}$.
\end{proposition}
\begin{proof}
  Let $\mathrm{pr}_{r,s}:\bar{I}_{(c,1,\frac{\delta}{r})}\P Y_{r,s}\rightarrow
  \bar{I}_{c}Y $ be the morphism induced from the projection of $\P Y_{r,s}$ to
  the base $Y$, we see $\mathrm{pr}_{r,s}$ is an isomorphism. Then the
  composition $\mathrm{pr}_{r,s}\!\circ \!ev_{h_{\infty}}\! \circ \!g$ is finite
  \'etale of degree $\frac{1}{a^{2}s\delta}$ as the composition can be also
  obtained by first forgetting root construction of $\mathcal M_{e}$ and taking
  rigidification afterwards. By two-of-there property for \'etale
  morphisms\cite[Lemmma 100.35.6]{stacks-project}, the \'etaleness of $g$ comes
  from that $ev_{h_{\infty}}$ is finite \'etale, which is proved in Proposition
  \ref{prop:etale-ev}. Finally as the composition $ev_{h_{\infty}}\circ g$ is of
  degree $\frac{1}{a^{2}s\delta}$ and $ev_{h_{\infty}}$ is of degree
  $\frac{1}{a\delta}$ by Proposition \ref{prop:aut-order}, we see that $g$ is of
  degree $\frac{1}{as}$.
\end{proof}

\begin{lemma}\label{lem:edge-contr}
  We have the following:
  \begin{enumerate}
  \item When the edge moduli $\mathcal M_{e}$ arises from the localization
    analysis in \ref{subsubsec:edge-cntr1}, we have that $[\mathcal
    M_{e}]^{vir}=[\mathcal M_{e}]$ and the Euler class of the virtual normal
    bundle is equal to $1$.
  \item When the edge moduli $\mathcal M_{e}$ arises from the localization
    analysis in \ref{subsubsec:edge-cntr2}, we have $[\mathcal
    M_{e}]^{vir}=[\mathcal M_{e}]$ and the inverse of Euler class of the virtual
    normal bundle is equal to
    $$\prod_{i=1}^{-1-\lfloor{-\delta}  \rfloor}\big(\lambda+\frac{i}{\delta}(
    c_{1}(L)-\lambda ) \big)\ .$$
  \end{enumerate}
\end{lemma}
\begin{proof}
  The first case follows from Proposition \ref{prop:etale-ev}. Now we calculate
  the edge contribution from the second case. Recall the divisor $E$ of
  $\mathfrak R$ introduced in Section \ref{subsec:space2} is isomorphic to $\P
  Y_{r,s}$. The normal bundle $N_{E/\mathfrak R}$ is isomorphic to $\mathcal
  O(-s \mathcal D_{0})$, and the fiber of the normal bundle $N_{E/\mathfrak R}$
  over $\mathcal D_{\infty}$ has $\C^{*}-$weight $1$, thus we have
  $f^{*}N_{E/\mathfrak R}= L_{as\delta\chi_{1}}^{\vee}\otimes \pi^{*} L$ as
  $\C^{*}-$equivariant line bundles over $\mathcal C_{e}$. Then we see that $$
  -R^{*}\pi_{*}f^{*}N_{E/\mathfrak R}=R^{1}\pi_{*}(\pi^{*}L\otimes
  L_{as\delta\chi_{1}}^{\vee})$$ in the $K-$group $K^{*}(\mathcal M_{e})$.

  Write $\mathfrak P:=\P Y_{r,s}$ for simplicity. Let $\mathbb E_{\mathfrak R}$
  (resp. $\mathbb E_{\mathfrak P}$) be the pull-back of the perfect obstruction
  theory of $\mathcal K_{0,2}(\mathfrak R,(0,\frac{\delta}{r},0))$ (resp.
  $\mathcal K_{0,2}(\mathfrak P,(0,\frac{\delta}{r}))$) to $\mathcal M_{e}$. As
  $\mathcal M_{e}$ is \'etale over the corresponding $\C^{*}$-fixed loci in
  $\mathcal K_{0,2}(\mathfrak P,(0,\frac{\delta}{r}))$ (or $\mathcal
  K_{0,2}(\mathfrak R,(0,\frac{\delta}{r},0))$) by Proposition
  \ref{prop:deg-edge}, we see that the $\C^{*}-$fixed parts $\mathbb
  E_{\mathfrak P}^{fix}$ and $\mathbb E_{\mathfrak R}^{fix}$ are perfect
  obstruction theory for $\mathcal M_{e}$. Note we have the distinguished
  triangle:
  \begin{equation}\label{eq:edge-relpot}
    \xymatrix{
      \mathbb E_{\mathfrak P} \ar[r] &\mathbb
      E_{\mathfrak R}\ar[r] &R^{\bullet}\pi_{*}f^{*}N_{E/\mathfrak R}\ar[r]^{+1} &\mathbb E_{\mathfrak P}[1]\ .
    }
  \end{equation}
  As shown in Lemma \ref{prop:etale-ev}, $\mathbb E_{\mathfrak P}$ is
  $\C^{*}-$fixed, then the movable part of $\mathbb E_{\mathfrak R}$ is equal to
  the movable part of $R^{\bullet}\pi_{*}f^{*}N_{E/\mathfrak R}$, which we will
  calculate below.

  The restriction of sections $x^{-sa\delta+asi}y^{ari }$($1\leq i \leq
  -1-\lfloor{-\delta} \rfloor$) to each fiber curve $C_{e}$ of $\mathcal C_{e}$
  over $\mathcal M_{e}$ is a basis of the vector space $H^{1}(C_{e},
  f|_{C_{e}}^{*}N_{E/\mathfrak R} )$. However they may not be sections of the
  vector bundle $R^{1}\pi_{*}(f^{*}N_{E/\mathfrak R} )$. Instead, by Proposition
  \ref{prop:equiv-sec}, one can see that $x^{-sa\delta+asi}y^{ari }$ is a
  section of vector bundle $R^{1}\pi_{*}(f^{*}N_{E/\mathfrak R} )\otimes
  L^{-\otimes \frac{i}{\delta}}\otimes \mathbb
  C_{-\frac{\delta-i}{\delta}\lambda} $. Therefore
  $R^{1}\pi_{*}f^{*}N_{E/\mathfrak R} $ is isomorphic to the direct sum of line
  bundles $\bigoplus_{i=1}^{-1-\lfloor{-\delta} \rfloor} L^{\otimes
    \frac{i}{\delta}}\otimes \mathbb C_{\frac{\delta-i}{\delta}\lambda}$. Then
  the inverse of the Euler class of the virtual normal bundle form $\mathbb
  E_{\mathfrak R}$ is equal to $\prod_{i=1}^{ -1-\lfloor{-\delta}
    \rfloor}\big(\lambda+\frac{i}{\delta}( c_{1}(L)-\lambda ) \big)$.

  The calculation above also shows that the $\C^{*}-$fixed part of $\mathbb
  E_{\mathfrak R}$ is the equal to $\mathbb E_{\mathfrak P}$, then we have
  $[\mathcal M_{e}]^{vir}=[\mathcal M_{e}]$ in Section
  \ref{subsubsec:edge-cntr2}.
\end{proof}

\begin{remark}
When $Y$ is of zero dimensional, our computation of the Euler class of the
virtual normal bundle coincides the computation
from~\cite[\S 6]{liu2020stacky}.  
\end{remark}
\end{appendix}

\bibliographystyle{amsxport} \bibliography{references}

\end{document}